\documentstyle[12pt]{article}
\def\m@th{\mathsurround=0pt }
\def\eqalign#1{\null\,\vcenter{\openup\jot \m@th
   \ialign{\strut\hfil$\displaystyle{##}$&$
      \displaystyle{{}##}$\hfil \crcr#1\crcr}}\,}
\newtheorem{theorem}{Theorem}[section]

\newtheorem{corollary}[theorem]{Corollary}
\newtheorem{lemma}[theorem]{Lemma}

\newtheorem{definition}[theorem]{Definition}

\def\ad{{\rm ad }\, }
\def\gkdim{{\rm GKdim }\ }

\def\wt{{\rm wt }}

\def\gr{{\rm gr}}
\def\bideg{{\rm bideg}}
\def\tip{{\rm tip}}
\def\ht{{\rm ht }}
\def\supp{{\rm supp }}
\begin{document}

\begin{center}
{\large{\bf  Coideal Subalgebras and  Quantum Symmetric Pairs } }
\vskip 5mm Gail Letzter\footnote{supported by NSA grant no. MDA
904-99-1-0033.

AMS subject classification 17B37}

Mathematics Department

Virginia Polytechnic Institute

and State University

Blacksburg, VA 24061

\end{center}

\begin{abstract}

Coideal subalgebras of the quantized enveloping 
algebra are surveyed, with selected proofs included.  
The first half of the paper studies generators, Harish-Chandra modules, 
and associated quantum homogeneous 
spaces. The second half   discusses various  well known  quantum coideal subalgebras
and the implications of the abstract theory on these examples.   The focus 
is on 
 the locally finite part of the quantized enveloping algebra, 
analogs of enveloping algebras of nilpotent Lie subalgebras, and coideals 
used to form quantum symmetric pairs.   The last family of examples is 
explored in detail.  Connections are made to the construction of quantum 
symmetric spaces.
\end{abstract} 

\medskip

The introduction of quantum groups  in the early 1980's has had a
tremendous influence on the theory of Hopf algebras. Indeed, quantum
groups provide a source of new and interesting examples.
We shall discuss the reverse impact:  the theory of quantum groups uses the Hopf
structure extensively.  This special structure is often 
  hidden in the classical setting, while it is prominent and 
  fundamental for quantum analogs.

 Let  ${\bf g}$ be a semisimple Lie algebra and write $G$ for 
   the corresponding   connected, simply connected  algebraic
 group. There are two standard types of quantum groups
associated to ${\bf g}$ and $G$. The
first is the quantized enveloping algebra which is a quantum
analog of the enveloping algebra of ${\bf g}$. The second 
is the quantized function algebra which is a quantum analog of the
algebra of regular functions on $G$. We will be focusing on a
particular aspect of the Hopf theory of both types of quantum groups: the study
of (one-sided) coideal subalgebras.

 One of the reasons coideal subalgebras are so important in the study 
of  quantum groups  is that quantum groups
 do not  have ``enough" Hopf subalgebras.  This shortage
of Hopf subalgebras is especially notable for quantized
enveloping algebras. 
Consider a Lie 
subalgebra ${\bf t}$ of the Lie algebra ${\bf g}$.   The enveloping 
algebra $U({\bf t})$ of ${\bf t}$ is a Hopf subalgebra of $U({\bf 
g})$.   However, upon passage to the quantum case, $U_q({\bf t})$, 
even when 
it is defined, is 
often not isomorphic to a Hopf subalgebra of $U_q({\bf g})$.  In many 
cases,
there are   subalgebras of  $U_q({\bf g})$ which are not 
Hopf subalgebras but are still  
good quantum analogs of $U({\bf t})$. Moreover, these subalgebras 
often turn out to be coideals. For example, let ${\bf g}={\bf 
n}^-\oplus {\bf h}\oplus {\bf n}^+$ be the triangular decomposition 
of ${\bf g}$.  There is a natural subalgebra $U^+$ of $U_q({\bf 
g})$ which is an analog of  the subalgebra  $U({\bf n}^+)$ of 
$U({\bf g})$.  This subalgebra $U^+$ is a 
coideal but is not a Hopf subalgebra of the quantizing enveloping algebra. 
Of more interest to us is the 
 fixed Lie subalgebra ${\bf g}^{\theta}$ corresponding to an
involution $\theta$ of ${\bf g}$. In the classical
case, the symmetric pair ${\bf g}^{\theta},{\bf g}$ is used to
form  symmetric spaces.  However, in the quantum case, $U_q({\bf 
g}^{\theta})$ does not usually embed inside of $U_q({\bf g})$.  
Thus it was initially unclear how to develop the  theory of quantum 
symmetric spaces.  In  [K], Koornwinder    constructed two-sided coideal analogs of ${\bf
g}^{\theta}$ in type $A_1$ and used them to produce  quantum 
symmetric spaces. More families of coideal analogs were
discovered in [N], [NS],[DN], and [L1]. In [L2], a uniform
approach was developed in the maximally split case using one-sided
coideal subalgebras of the quantized enveloping algebra. The
one-sided coideal condition turned out to be critical in
characterizing these quantum analogs of $U({\bf g}^{\theta})$.

Quantum symmetric spaces were first defined using the quantized 
function algebra. (See for example [KS, 11.6.3 and 11.6.4].)  Koornwinder's work [K] inspired the development of a
 quantum
symmetric space theory   using analogs of  ${\bf 
g}^{\theta}$  contained in $U_q({\bf g})$.   The axiomatic theory of quantum
symmetric spaces (see [Di]) proceeded more rapidly than the
discovery of a general way to construct  examples.
Indeed, Dijkhuizen [Di, end of Section 3] outlined the desirable properties 
that analogs 
of ${\bf g}^{\theta}$ contained in $U_q({\bf g})$
 should have in order to form ``nice" quantum
symmetric spaces. As in Koornwinder's work [K], one of the key
properties is the coideal condition. Another crucial property of an 
analog is that its finite-dimensional spherical modules be 
characterized in a similar way to the characterization in the 
classical case. This   is obtained in [L3] for the coideal subalgebras of
[L2].  The proof  uses  quantum
Harish-Chandra modules associated to quantum symmetric pairs.  
The coideal condition  plays a
prominent role   in defining and developing the theory of quantum Harish-Chandra 
modules ([L3]).

This paper is based on a talk given at the MSRI Hopf Algebra
Workshop. It offers a panorama of the use of coideal subalgebras in
constructing quantum symmetric pairs, in forming   quantum
Harish-Chandra modules,  and in producing quantum symmetric
spaces. In the first half  of the paper, we present topics in the
general theory of quantum coideal subalgebras.   Section 1 sets 
notation and presents some basic facts about coideal subalgebras 
inside arbitrary Hopf algebras. 
  In
Section 2, we define  Harish-Chandra
modules associated to quantum ``reductive" pairs. We prove a basic 
result: every $U_q({\bf g})$ module contains a large Harish-Chandra module associated to a quantum reductive pair. 
In Section 3, we discuss how coideal subalgebras of the quantized
enveloping algebra can be used in the dual quantum function
algebra setting.  Connections are made to the
  theory of   quantum homogeneous spaces.  
  Section 4 studies
  filtrations  on the quantized
enveloping algebra and their impact on coideal subalgebras. As
a result, we  obtain a nice description of the
generators of a coideal subalgebra under mild restrictions.
 
The final three sections are devoted to specific  coideal
subalgebras of the quantized enveloping algebra.  Section 5
discusses the locally finite part, $F(U)$, of $U_q({\bf g})$.   It is 
well known that the classical enveloping algebra   $U({\bf g})$  can be written as a 
direct sum  of finite-dimensional $\ad {\bf g}$ modules.   This result 
plays an important role in understanding the structure of $U({\bf g})$
and   classifying its primitive ideals. 
Unfortunately, the quantized enveloping algebra contains infinite 
dimensional $U_q({\bf g})$ modules with respect to the adjoint 
action.   Thus it is often necessary to use the locally finite part,   
$F(U)$, 
which is the maximal subalgebra of $U_q({\bf g})$ that can be written 
as a direct sum of finite-dimensional simple ad $U_q({\bf g})$ modules. This 
algebra
$F(U)$     is 
not a Hopf subalgebra of $U_q({\bf g})$, but it is a coideal.  The structure of this
coideal subalgebra is briefly reviewed with some consideration for the 
implications
of the results of Section 4. Certain quantum  Harish-Chandra
modules defined originally in [JL3] using $F(U)$ are  elucidated
in terms 
of the general approach presented in Section 2.
Section 6 considers coideal subalgebra analogs of
 enveloping algebras of nilpotent and parabolic Lie
subalgebras of ${\bf g}$.  Much of the material in this section
 is based on  
[Ke].
The last part, Section 7, is devoted to the theory of quantum
symmetric pairs. This material is largely drawn from
[L2] and [L3]. However, since the papers appeared, we have found 
simpler    approaches which are presented here with many proofs  included. 
We show how to lift a 
maximally split
involution $\theta$ of ${\bf g}$ to the quantum setting.  Exploiting this lift, we 
define a coideal subalgebra $B_{\theta}$ of $U_q({\bf g})$.
 As in [L2],   $B_{\theta}$ is
characterized as the ``unique'' maximal coideal subalgebra of 
$U_q({\bf g}^{\theta})$ which specializes to $U({\bf g}^{\theta})$ 
as $q$ goes to $1$.  
  Using the results of Section 4, we  give a new proof of this
uniqueness theorem   which
  does not involve the intricate specialization arguments found in [L2]. 
  We also take the opportunity to make some corrections in the case 
  work necessary to make the uniqueness tight.  
Results on the Harish-Chandra module and quantum symmetric space theory
associated to these pairs are described.

\bigskip
\noindent
{\bf Acknowledgement.}
Part of this paper was written while the author spent a month as a 
visiting professor at the
University of Rheims. The author would like to thank the mathematics 
department there for their hospitality and J. Alev for his valuable comments.
  The author would  also like to thank the referee for a painstakingly 
careful reading of the first draft and  many useful 
suggestions.   Finally, the author would like to thank Dan Farkas 
whose support transcends multiple revisions.

 \section{Background and Notation}
 Let $H$ be a Hopf algebra over a field $k$ with comultiplication ${\it\Delta}$,
 antipodal map $\sigma$, and counit $\epsilon$.
 Given any $a\in H$,
 write ${\it\Delta}(a)=\sum a_{(1)}\otimes a_{(2)}$ using Sweedler notation.
 A vector subspace $I$ of $H$ is called a
 left   coideal  if $${\it \Delta}(I)\subset H\otimes I.$$
 Similarly, $I$ is called a right coideal if
 ${\it \Delta}(I)\subset I\otimes H.$ In particular, a left
 (resp. right) coideal is a
 left (resp. right) $H$ comodule contained in $H$. If $I$ is both a left
 (resp. right) coideal and a subalgebra of $H$, then we simply say that $I$
is a left
 (resp. right) coideal subalgebra.    There is also a notion of
 two sided coideals but those are  generally not considered  here.
 We will
 usually choose to discuss left coideals and left coideal subalgebras;
 analogous results can be proved for the right-handed versions.

We first present   two general results
about coideals
 inside of an arbitrary Hopf algebra.
First, assume that $H$ contains the group algebra $kG$ of a group
$G$. Choose a vector space complement $Y$ to $kG$ in $H$.  Let
$P$ be the  projection map of $H$ onto $kG$ as vector spaces using
the decomposition $H=kG\oplus Y$.  Assume that $H$ is a left $kG$
comodule where
 the comodule structure comes from the comultiplication and the projection
$P$.  In
 particular, $H$ is the direct sum of vector subspaces $H_g$ where
 $$(P\otimes Id){\it\Delta}(H_g)\subset g\otimes H_g.\leqno{(1.1)}$$  Given any
 left coideal $I$ of $H$, set $I_g=I\cap H_g$.

 \begin{lemma}
A left coideal  $I$ contained in  $H$ is equal to a direct sum of
the
 vector spaces $I_g$.  Thus $I$ is a left $kG$ comodule.
 \end{lemma}

  \noindent
 {\bf Proof:}  Write $a\in I$ as $a=\sum_{g\in G}a_g$ where each
 $a_g\in H_g$.  The lemma follows if we show that each $a_g\in
 I$.  By (1.1), $${\it\Delta}(a)\in\sum_{g\in G}g\otimes a_g+Y\otimes
 H.$$
 The coideal property now ensures that each $a_g\in I$.
 $\Box$

\medskip
Every Hopf algebra $H$ comes equipped with a (left) adjoint
action. Using this adjoint action, $H$ becomes an $(\ad H)$
module. In particular, given $a,b\in H$, $$(\ad
a)b=\sum a_{(1)}b\sigma(a_{(2)}).\leqno{(1.2)}$$
 In the quantum case,
it is often interesting to consider ad-invariant coideals. The
following result (which is basically [Jo, Lemma 1.3.5]) is
particularly useful.

\begin{lemma}
Let $I$ be a left coideal in $H$ and let $M$ be a Hopf subalgebra
of $H$.   Then $(\ad M)I$ is an ad $M$ invariant left coideal of $H$.
\end{lemma}

\noindent {\bf Proof:}
First note that
 ([Jo, 1.1.10])
$${\it\Delta}(\sigma(a))=\sum \sigma(a_{(2)})\otimes
\sigma(a_{(1)}).\leqno{(1.3)}$$  Hence $${\it\Delta}((\ad a) b)=
{\it\Delta}(\sum a_{(1)}b\sigma(a_{(2)}))
=\sum (a_{(1)}b_{(1)}\sigma(a_{(4)}))\otimes
(a_{(2)}b_{(2)}\sigma(a_{(3)})).$$  The result follows from the
fact that ${\it\Delta}(a_{(2)})=\sum a_{(2)}\otimes a_{(3)}$. $\Box$
\medskip

Before defining the quantized enveloping algebra, we recall some 
basic facts about semisimple Lie algebras.
Denote the set of nonnegative integers by ${\bf N}$, the 
complex numbers by ${\bf C}$, and the real numbers 
 by ${\bf R}$.
 Let ${\bf g}$ be a complex semisimple Lie algebra
 with triangular decomposition ${\bf n}^-\oplus{\bf h}\oplus{\bf
 n}^+$. Write $h_1,\dots, h_n$ for a basis of ${\bf h}^*$.
  Let $\Delta$ denote the root system of ${\bf g}$ and write 
  $\Delta^+$ for the set of positive roots.  Recall that $\Delta$ is 
  a subset of   ${\bf h}^*$.  Furthermore,
    ${\bf n}^+$
 (resp. ${\bf n}^-$) has 
 a basis of root vectors $\{e_{\beta}|\beta\in {\Delta^+}\}$ (resp. 
 $\{f_{-\beta}|\beta\in {\Delta^+}\}$). These root vectors are common 
  eigenvectors, called weight vectors,
   for the adjoint action of ${\bf h}$ on ${\bf g}$. 
   In particular, $(\ad h_i)e_{\beta}=[h_i,e_{\beta}]=\beta(h_i)e_{\beta}$
   and  $(\ad h_i)f_{-\beta}=[h_i,f_{-\beta}]=-\beta(h_i)f_{-\beta}$
   for each $\beta\in \Delta^+$. We further assume that 
   $\{e_{\beta},f_{-\beta}|\beta\in {\Delta^+}\}\cup\{h_1,\dots, h_n\}$
   is a Chevalley basis for ${\bf g}$ ([H, Theorem 25.2]).
  Let  $\pi=\{\alpha_1,\dots,\alpha_n\}$ denote the
 set of  simple roots in ${\Delta^+}$ and  $(\ ,\ )$ denote the Cartan inner product
 on ${\bf h}^*$. Recall further that $(\ ,\ )$ is positive definite on 
 the real vector space spanned by $\pi$.  The set $\pi$ is a basis for   ${\bf h}^*$. Given 
 $\alpha_i\in \pi$, we write $e_i$ (resp. $f_i$) for $e_{\alpha_i}$
 (resp. $f_{-\alpha_i}$). The Cartan matrix associated to the root 
 system $\Delta$ is the matrix with entries 
 $a_{ij}=2(\alpha_i,\alpha_j)/(\alpha_i,\alpha_i)$. (The reader is 
 referred to  [H, Chapters II
 and III] for
 additional information on  semisimple Lie algebras and root systems.)
 
 Let $q$ be an
 indeterminate and set $q_i=q^{(\alpha_i,\alpha_i)/2}$.  
 Let $[m]_q $ denote the $q$ number $(q^{m}-q^{-m})/(q-q^{-1})$ and
 $[m]_q!$ denote the $q$ factorial $[m]_q[m-1]_q\cdots [1]_q$.  The 
 $q$ binomial coefficients are defined by $$\left[\eqalign{&m\cr 
 &j\cr}\right]_q={{[m]!_q}\over{[j]!_q[m-j]!_q}}.$$
   The quantized
enveloping algebra
 $U=U_q({\bf g})$ is generated by $x_1,\dots, x_n$, $t_1^{\pm
 1},\dots, t_n^{\pm 1}, y_1,\dots, y_n$ over ${\bf C}(q)$ with 
 the relations listed below (see for example [Jo, 3.2.9] or [DK, 
 Section 1]).  
 \begin{enumerate}
 \item[(1.4)]$x_iy_j-y_jx_i=\delta_{ij}(t_i-t_i^{-1})/(q_i-q_i^{-1})$ for each 
 $1\leq i\leq n$.
 \item[(1.5)] The $t_1^{\pm
 1},\dots, t_n^{\pm 1}$ generate a free abelian group $T$ of rank $n$.
 \item[(1.6)] $t_ix_j=q^{(\alpha_i,\alpha_j)}x_jt_i$ and
 $t_iy_j=q^{-(\alpha_i,\alpha_j)}y_jt_i$  for all $1\leq i,j\leq n$.
 \item[(1.7)] The quantum Serre relations:
 $$\sum_{m=0}^{1-a_{ij}}(-1)^m\left[\eqalign{1&-a_{ij}\cr 
 &m\cr}\right]_{q_i}x_i^{1-a_{ij}-m}x_jx_i^{m}=0$$
 and $$
  \sum_{m=0}^{1-a_{ij}}(-1)^m\left[\eqalign{1&-a_{ij}\cr 
 &m\cr}\right]_{q_i}y_i^{1-a_{ij}-m}y_jy_i^{m}=0$$ for all 
  $1\leq i,j\leq n$ with  
 $i\neq j$.
 \end{enumerate}

 The algebra $U$ is a Hopf algebra with comultiplication
${\it\Delta}$, antipode $\sigma$, and counit $\epsilon$  defined on generators as follows. 
 \begin{enumerate} 
\item[(1.8)]${\it \Delta} (t) =
t\otimes t\qquad \epsilon(t)=1\qquad
\sigma(t )=t ^{-1}$ 
for all $t$ in  $T$
\item[(1.9)]${\it \Delta}(x_i)= x_i\otimes 1 + t_i\otimes x_i\qquad
\epsilon(x_i)=0\qquad \sigma(x_i)= -t_i^{-1}x_i$ 
\item[(1.10)]${\it \Delta}(y_i)= y_i\otimes t_i^{-1} + 1\otimes y_i\qquad
\epsilon(y_i)=0\qquad \sigma(y_i)= -y_it_i$
\end{enumerate}
for $1\leq i\leq n$.

It is well known that the algebra $U$ specializes to $U({\bf g})$ as 
$q$ goes to $1$.   This can be made more precise as follows. Set $A$ 
equal to ${\bf C}[q,q^{-1}]_{(q-1)}$.  Let 
$\hat U$ be the $A$ subalgebra of $U$ generated by $x_i,y_i,t_i^{\pm 
1},$ 
and $(t_i-1)/(q-1)$ for $1\leq i\leq n$.   Then $\hat U\otimes_A{\bf 
C}$ is isomorphic to $U({\bf g})$. (See for example [L2, beginning of 
Section 2]).  Given a subalgebra $S$ 
of $U$, set $\hat S=S\cap \hat U$.   We say that $S$ specializes to 
the subalgebra $\bar S$ of $U({\bf g})$ if the image of $\hat S$ in  $\hat U\otimes_A{\bf 
C}$ is $\bar S$.
    
    Set
$Q(\pi)$ equal to the integral lattice generated by $\pi$. Let
$Q^+(\pi)$ (resp. $Q^-(\pi)$) be the subset of $Q(\pi)$ consisting
of nonnegative (resp. nonpositive) integer linear combinations of
elements in $\pi$. The standard partial
  ordering  on the root lattice $Q(\pi)$ is defined by 
  $\lambda\geq\mu$ provided $\lambda-\mu$ is in $Q^+(\pi)$. Let $P^+(\pi)$ denote the set of dominant integral
weights associated to $\pi$.  In particular, $\lambda\in {\bf h}^*$
is an element of $P^+(\pi)$ if and only if 
$2(\lambda,\alpha_i)/(\alpha_i,\alpha_i)$ is a nonnegative integer 
for all $1\leq i\leq n$. There is an isomorphism $\tau$ of
abelian groups from $Q(\pi)$ to $T$ defined by
$\tau(\alpha_i)=t_i$, for $1\leq i\leq n$. Using this isomorphism,
we can replace condition (1.6) with
$$\tau(\lambda)x_i\tau(\lambda)^{-1}=q^{(\lambda,\alpha_i)}x_i{\rm
\ \ and\ }
\tau(\lambda)y_i\tau(\lambda)^{-1}=q^{-(\lambda,\alpha_i)}y_i\ 
\leqno{(1.11)}$$
for all $\tau(\lambda)\in T$ and $1\leq i\leq n$.

Let $M$ be a $U$ module. A nonzero vector $v\in U$ has 
weight $\gamma\in {\bf h}^*$ provided that $\tau(\lambda)\cdot 
v=q^{(\lambda,\gamma)}v$ for all $\tau(\lambda)\in T$. 
Given a subspace $V\subset M$, the subspace of $V$ spanned by
the 
$\gamma$ weight vectors  is called 
the $\gamma$ weight space of $V$ and denoted 
by $V_{\gamma}$.   Now $U$ can be given the structure 
of a $U$ module using the quantum adjoint action (1.2).  Let  $v$ be 
an element of $ U$. We say that $v$ has weight $\gamma$
provided that it is a $\gamma$ weight vector in terms of this adjoint 
action.   In particular, $v$ has weight $\gamma$ if 
$\tau(\lambda)v\tau(\lambda)^{-1}=q^{(\lambda,\gamma)}v$ for all
$\tau(\lambda)\in T$.

  Let $G^+$ be the subalgebra of $U$ generated by $x_1t_1^{-1},\dots,
 x_nt_n^{-1}$.  Similarly, let  $U^-$ be the subalgebra of $U$ generated by
 $y_1,\dots, y_n$.  Let $U^o$ be the group algebra of $T$. It is well known
that both
 $U^-$ and $G^+$ are a direct sum of their weight spaces.
  The quantized enveloping algebra
   $U$ admits a triangular decomposition.  More precisely, there is 
   an  isomorphism of vector spaces using the multiplication map ([R]):
 $$U\cong U^-\otimes U^o\otimes G^+.\leqno{(1.12)}$$
  It follows that  there is a direct sum decomposition
 $$U=\bigoplus_{t\in T}U^-G^+t.\leqno{(1.13)}$$
 Let $G^+_+$ (resp. $U^-_+$) denote the augmentation ideal of $G^+$
 (resp.  $U^-$)  and set $Y$ equal to the vector space
 $(U^-_+G^+U^o+ U^-G^+_+U^o)$.
  The  direct sum decomposition (1.13) implies that
 $$U=U^o\oplus Y.\leqno{(1.14)}$$ Using the definition of  the comultiplication of
 $U$, it is straightforward to check that for any $b\in U^-G^+t$
 $${\it\Delta}(b)\in t\otimes b+Y\otimes U.$$  Thus
 the projection of ${\it\Delta}(U)$ onto $U^o\otimes U$ makes $U$ into a left $U^o$ comodule
 with $U_t=U^-G^+t$. Hence we have the following version of Lemma 1.1 
 for quantized enveloping algebras.
 
 \begin{lemma}
 Let $I$ be   left coideal   of $U$. Then   $$I=\bigoplus_{t\in 
 T}(I\cap U^-G^+t).$$
 \end{lemma}

 \section{Harish-Chandra modules}

Consider a Lie subalgebra ${\bf k}$ of the semisimple Lie algebra 
${\bf g}$.  A Harish-Chandra module associated to the
pair ${\bf g},{\bf k}$ is a ${\bf g}$ module which can be written
as a direct sum of finite-dimensional simple ${\bf k}$ modules.
Harish-Chandra modules are an important tool in classical 
representation theory.  This is especially true when 
  ${\bf g},{\bf k}$ is a  symmetric pair.  Harish-Chandra 
  modules associated to symmetric pairs provide an algebraic approach
  to the
  representation theory of real reductive Lie groups.

There is a nice introduction 
to the theory of Harish-Chandra modules presented
 in [D, Chapter 9] from an algebraic point of view.    
The first section of [D, Chapter 9] only assumes that  
 ${\bf k}$ is  reductive in ${\bf g}$. 
 A basic result which is used repeatedly in this chapter of [D] is the following.
 
 \medskip
\noindent
\begin{theorem} ([D, Proposition 1.7.9]) The direct sum  of all the
finite-di-

\noindent mensional simple ${\bf k}$ modules inside a ${\bf
g}$ module is a Harish-Chandra module for the pair ${\bf g}, {\bf
k}$.
\end{theorem}

\noindent
This theorem allows one to find large Harish-Chandra 
modules inside of infinite-dimensional ${\bf g}$ modules.  Its proof
uses the fact that   ${\bf k}$ is  reductive in  ${\bf g}$ and 
that $U({\bf g})$ is a 
locally finite $\ad U({\bf g})$ module. 
 
   In the quantum 
  setting, $U_q({\bf k})$ is not always a subalgebra of $U_q({\bf g})$
  when ${\bf k}$ is a Lie subalgebra  of ${\bf g}$.   However, 
  one   often finds a quantum analog of $U({\bf k})$ which is a
   one-sided 
 coideal subalgebra  of $U_q({\bf g})$.  
   Thus any good 
 theory of quantum Harish-Chandra modules must work for pairs 
 $U_q({\bf g})$, $I$ where $I$ is a one-sided coideal subalgebra of 
 $U_q({\bf g})$. In order to begin such a theory, it is necessary to 
 have an analog of Theorem 2.1.   This presents two 
 difficulties.   The first is that $U$, in contrast to the classical 
 situation,  is not a locally finite $\ad U$ 
 module.  (We will return to this obstruction in Section 4.) The second is: what does it mean for a  coideal subalgebra
 to be reductive in $U$?

  In this section, we present a quantum version of Theorem 2.1   using 
  the locally finite part $F(U)$ of $U$ and  a certain condition on 
  coideal subalgebras which substitutes for reductivity.      
  The
 material of this section is based on  [L1, Section 4] and [L3, Section 3].
This result sets the stage for the development of a quantum Harish-Chandra 
module theory.   Indeed, the author has checked that 
many  of the results of [D, Section 9.1] and their 
 proofs carry over to  coideal subalgebras satisfying this quantum 
 version of Theorem 2.1.   
Some properties of 
quantum principal series modules  analogous to those in  [D, Section 9.3] 
are proved   in 
[L3, Section 6].  Spherical modules (see [D, 9.5.4]) have been classified in the 
quantum case (see Section 7, Theorem 7.7  and [L3, Section 4]).  This 
is discussed further in Section 7.

Recall the definition of the adjoint action (1.2) and define
$$F(U)=\{v\in U| \dim(\ad U)v <\infty\}.\leqno{(2.1)}$$  By [JL1, 
Corollary 2.3, Theorem 5.12, and Theorem 6.4],
$F(U)$ is an algebra, it can be written as a direct sum of
finite-dimensional simple $U$ modules,  and it is ``large'' in $U$.
It is also true that  $F(U)$ is   a left
coideal of $U$, a subject we will return to in Section 5.

Fix a left coideal subalgebra $I$ of $U$.  Note that
$$F(U)I=\{\sum f_ir_i|f_i\in F(U), r_i\in I\}$$ is also a
subalgebra of $U$. This follows from the fact  that
$rf=\sum r_{(1)}\epsilon(r_{(2)})f=\sum ((\ad r_{(1)})f)r_{(2)}$  for any
$r\in I$ and $f\in F(U)$.   Since $I$ is a left  coideal,
each $r_{(2)}\in I$. Furthermore the ad-invariance of   $F(U)$ implies
that $(\ad r_{(1)})f$ is in $F(U)$. We use $F(U)I$ to define
Harish-Chandra modules.  

\begin{definition}   A {\it Harish-Chandra module} for the pair $U,I$
is an $F(U)I$ module which is a direct sum of finite-dimensional
simple $I$ modules.
\end{definition}

Let us take a closer look at the condition that ${\bf k}$ is 
 reductive in ${\bf g}$.   Reductivity means  that
$(\ad {\bf k})$ acts semisimply on ${\bf g}$.   This assumption is
enough to prove that ${\bf k}$ is itself reductive and that   the center of
${\bf k}$ can be extended to a Cartan subalgebra of ${\bf g}$. It
is unclear what  the corresponding assumption in the quantum case,
namely that $(\ad I)$ acts semisimply on $F(U)$, implies. It seems
unlikely that this assumption alone will yield an analog of
Theorem 2.1.

 Of course, there would be no problem if   $I$   acted semisimply
 on
all finite-dimensional $I$ modules. When $I$ turns out to be a
Hopf subalgebra of $U$ isomorphic to a quantized enveloping
algebra of a semisimple Lie subalgebra of ${\bf g}$, this is certainly true.
However, complete reducibility does not hold in general for the
large class of coideal subalgebras considered in Section 7.   Thus
we need a replacement
 for the notion of reductive in ${\bf g}$.  This substitute is 
 invariance under the action of a certain conjugate linear 
 antiautomorphism of $U$.

Let  $\kappa$ denote  the conjugate linear form of 
the quantum Chevalley
antiautomorphism.  In particular, let $U_{{\bf R}(q)}$ denote the 
${\bf R}(q)$ subalgebra of $U$ generated by $x_i,y_i,t_i^{\pm 1}$,
for $1\leq i\leq n$.  The 
antiautomorphism $\kappa$ of $U_{{\bf R}(q)}$ is
 defined by $\kappa(x_i)=y_it_i$,
$\kappa(y_i)=t_i^{-1}x_i$ and $\kappa(t)=t$ for all $t\in T$. We then 
extend $\kappa$ to $U$  using conjugation.   More precisely, given 
$a\in{\bf C}$, write $\bar a$ for the complex conjugate of $a$.   
Set $\bar q=q$.   Then $\kappa(au)=\bar a\kappa(u)$ for all $u\in 
U_{{\bf R}(q)}.$

It is straightforward to check using (1.8), (1.9), and (1.10) that
$${\it\Delta}(\kappa(b))= (\kappa\otimes
\kappa){\it\Delta}(b)\leqno{(2.2)}$$ for all $b\in U$.  
Moreover $\kappa$ 
 gives $U$ the structure 
of a Hopf $*$ algebra where $*=\kappa$ ([CP, Section 4.1F ]).

 For the remainder of this section, we assume
that $I$ is a left coideal subalgebra such that $\kappa(I)=I$.
Thus one can think of $I$ as a $*$ subalgebra of $U$.

The field ${\bf R}(q)$ can be made into a  real ordered field ([J, 
Section 11.1]) where the positive 
elements are defined as follows.  Write a polynomial $f(q)$ in the 
form 
$(q-1)^s(f_m(q-1)^m+\cdots +f_1(q-1)+f_0)$ where each $f_i\in {\bf R}$ 
and  $m,s\in{\bf N}$.  Then $f(q)$ is positive if and only 
if $f_0>0$.  An element $h\in {\bf R}(q)$ is positive 
if and only if $h$ can be written as a quotient of positive 
polynomials. This induces a total order on ${\bf R}(q)$.

 We next specify  a class of ``nice''
finite-dimensional  $I$ modules. 

\begin{definition}   An $I$ module $W$ is
called {\bf unitary}  if it admits a sesquilinear form $S_W$ (i.e. 
 linear in the 
first variable and conjugate linear in the second variable) such that 
\begin{enumerate}
\item[{(i)}] $S_W(av, w)=S_W(v,\kappa(a)w)$ for all $a\in I$ and $v, 
w$ in $W$
\item[{(ii)}]  $S_W(v,v)$ is a positive element of  ${\bf R}(q)$ 
for each nonzero vector $v\in W$
\item[(iii)] $S_W(v,w)=\overline{S_W(w,v)}$ for all $v$ and $w$ in $W$.
\end{enumerate}
\end{definition}

Let $W$ be a finite-dimensional  unitary $I$ module. Choose a nonzero 
vector $v\in W$ such that $S_W(v,v)=1$.   Now suppose that $w\in W$ 
such that $S_W(w,v)=0$.  By Definition 2.3(iii), it follows that 
$S_W(v,w)$ also equals zero.
Hence one can show using induction    that $W$
has an orthonormal basis with respect to $S_W$.

The following result and its corollary show  that $I$ has an
extensive family of unitary modules, namely the finite-dimensional
simple $I$ submodules of any finite-dimensional simple $U$ module.

\begin{theorem}
Every finite-dimensional   unitary $I$ module can be written as a
direct sum of finite-dimensional simple unitary $I$ modules.
\end{theorem}

\noindent {\bf Proof:} Let $W$ be a finite-dimensional unitary $I$ module
with sesquilinear form $S=S_W$ as in Definition 2.3. Let $V$
be a finite-dimensional simple $I$ submodule inside of
$W$. By Definition 2.3, the restriction of $S$ to $W$ again 
satisfies conditions (i), (ii), and (iii). Furthermore, Definition 2.3(i) implies that the
orthogonal complement $W^{\perp}$ of $W$ with respect to $S$ is an
$I$ module. Hence $V\cong W\oplus W^{\perp}$,
  a direct sum of unitary $I$ modules with smaller dimension.   The proof
  follows by induction on $\dim V$.
  $\Box$

\begin{corollary}
Every finite-dimensional simple $U$ module $V$ is a Harish-Chandra
module for the pair $U$, $I$.
\end{corollary}

\noindent {\bf Proof:} Let $V$ be a finite-dimensional simple 
$U$ module.  It is well known that finite-dimensional 
$U$ modules are a direct sum of their weight spaces.   Moreover, the 
weight space of maximal weight is one dimensional.
Let $v$ be a basis vector for this highest weight space   and 
note that $v$ generates $V$ as a $U$ module.  The vector $v$ is called 
a highest weight generating vector of $V$. Recall that $U^-_+$ denotes 
the augmentation ideal of $U^-$. Note that $U^-_+v$ is the subspace of $V$ 
spanned by those weight vectors whose weights are    strictly less than 
that 
 of $v$. Furthermore,   $V$ is the 
direct sum of ${\bf C}(q)v$ and   $U^-_+v$. 

Let $\varphi$ 
be  the projection of $U$ 
onto $U^o$ using the direct sum decomposition $U=U^o\oplus Y$ (1.14). 
  Define a  sesquilinear form $S$ on $V$
by
$S(v,v)=1$ and $S(av,bv)=S(v,\varphi(\kappa(a)b)v)$ for all 
$a,b\in U$. Since $\varphi(b)=0$ for all $b$ in $U^-_+$, it follows 
that $S(v,U^-_+v)=0$.

Note  that $S$ satisfies Definition 2.3(i).  As in say ([L1, Lemma 4.2]), $S$
specializes to the classical positive definite Shapovalov  form of [Ka, 11.5 and Theorem 11.7].  Thus ([L1, Lemma 4.2])
 $S(w,w)\neq 0$ for any nonzero 
vector $w\in V$. It is straightforward to check that $S$ restricts to a 
${\bf R}(q)$ bilinear form on $U_{{\bf R}(q)}v$. Moreover, $S$ restricted 
to $\hat U_{{\bf R}(q)}$ takes values in ${\bf R}[q,q^{-1}]_{(q-1)}$.
  Let 
$w$ be an element in $\hat U_{{\bf R}(q)}v$. We can write $S(w,w)=f(q)$ 
with $f(q)$ in ${\bf R}[q,q^{-1}]_{(q-1)}$.   Since 
the specialization of $S$ is positive definite,  we must have that 
$f(1)\geq 0$.   It follows that $f(q)\geq 0$.  This fact and the 
nondegeneracy property    implies  that $S$ 
satisfies Definition 2.3(ii).  

Recall the direct sum decomposition 
(1.14) of $U$.   Note that $\kappa(Y)=Y$ and $\kappa(a)=a$ for all 
$a\in U^o\cap U_{{\bf R}(q)}$.   Therefore $\varphi(\kappa(b))=\varphi(b)$ for 
all $b\in U_{{\bf R}(q)}$.  It follows that $S$ is symmetric when
 restricted to  $U_{{\bf R}(q)}v$.  In particular, $S$ satisfies 
 Definition 2.3(iii). Thus   $V$ is a  unitary $I$ module.  The 
result now follows from Theorem 2.4.   $\Box$

\medskip 
 Note that we cannot expect  
  ${\it\Delta}(I)$ to be a subset of $I\otimes I$. Hence
 the tensor product of two $I$ modules does not necessarily admit an 
 action of  $I$ via the comultiplication of $U$. 
 However, since $I$ is a left 
coideal, the tensor product $V\otimes W$ of a $U$ module $V$ with an 
$I$ module
$W$ is an 
$I$ module. In particular, $a(v\otimes w)=\sum a_{(1)}v\otimes a_{(2)}w$ 
for all $v\otimes w\in V\otimes W$ and $a\in I$.  
The next lemma shows that the notion of unitary behaves well with
respect to the tensor product  of a $U$ module  with an $I$ module.

\begin{lemma}
Let $V$ be a finite-dimensional   $U$ module and let $W$ be a
finite-dimensional   unitary $I$ module.   Then $V\otimes W$ is a
finite-dimensional unitary $I$ module.
\end{lemma}

\noindent {\bf Proof:} Let $S_V$ (resp.  $S_W$) denote the
sesquilinear form on $V$ (resp. $W$) satisfying the conditions
of Definition 2.3.   Define a sesquilinear form $S=S_{V\otimes
W}$ on $V\otimes W$ by setting $S(a\otimes b,a'\otimes b')=
  S_V(a,a')S_W(b,b')$. It is easy to check   Definition  2.3(iii) 
  holds for $S$.  Let $\{v_i\}$ be an 
   orthonormal basis for $V$ with respect to $S_V$ and let $\{w_i\}$ be 
   an orthonormal basis for $W$ with respect to $S_W$.   Then $S(\sum 
   b_{ij}v_i\otimes w_j,\sum b_{ij}v_i\otimes w_j)=\sum b_{ij}\bar 
   b_{ij}$.   Thus Definition 2.3(ii) holds for $S$.  Condition (2.2) on $\kappa$ ensures that $S$
satisfies
  Definition 2.3(i).  In particular, for  $c\in I$,
  we have $S(c(a\otimes b),a'\otimes b')=S(\sum c_{(1)}a\otimes c_{(2)} b,
a'\otimes b')
  =S( a\otimes b,\sum \kappa( c_{(1)})a'\otimes \kappa(c_{(2)}) b')=S(a\otimes
b,\kappa(c)(a'\otimes b')).$
  $\Box$

  \medskip

  We now obtain a quantum analog of Theorem 2.1.

  \begin{theorem}  The sum of all the finite-dimensional simple
  unitary $I$ modules inside of the $F(U)I$ module $M$ is a
  Harish-Chandra module for the pair $U,I$.
  \end{theorem}

  \noindent
  {\bf Proof:} Assume that $W$ is a finite-dimensional simple unitary
  $I$ module contained in $M$. It suffices to show that the $F(U)I$ module
  generated by $W$ is a direct sum of finite-dimensional simple
  unitary
  modules.  Note that  
  $F(U)IW=F(U)W=IF(U)W$ is 
  an
$I$ module.  The vector space $F(U)\otimes W$ is also an $I$ module 
where the action is given  by   $$a\cdot (f\otimes w)=\sum( \ad
  a_{(1)})f\otimes a_{(2)} w$$ for all $f\in F(U)$, $w\in W$, and  $a\in I$.
  Furthermore, $F(U)W$ is  a homomorphic image of the $I$ module
$F(U)\otimes W$.
  Recall that
  $F(U)$ is a direct sum of finite-dimensional simple $(\ad U)$
  modules.  By Corollary 2.5, each finite-dimensional simple $(\ad U)$ module
  is a unitary $I$ module.
  Thus by Lemma 2.6, $F(U)\otimes W$, and hence $F(U)W$,
   splits into a direct sum of finite-dimensional simple unitary $I$ modules.
  $\Box$

\medskip
 Let ${\cal H}_{\bf R}$ be the set of all Hopf algebra 
automorphisms of $U$ which restrict to a Hopf algebra automorphism 
of $U_{{\bf R}(q)}$. Let $\Upsilon\in {\cal H}_{\bf R}$.  Suppose that $I$ is a left coideal subalgebra such that 
$\Upsilon^{-1}\kappa\Upsilon(I)=I$.   Then the results of this section 
hold for $I$ where we define unitary $I$ modules using 
$\Upsilon^{-1}\kappa\Upsilon$ instead of $\kappa$.

  \section{The Dual Picture}

In this section, we consider the connection between coideal
subalgebras of $U$ and coideal subalgebras inside the Hopf dual of
$U$.  The results presented in this section are well known and are
related to the theory of quantum homogeneous spaces. A good
reference for most of the material presented here  and for other
basic results about quantum homogeneous spaces is  [KS, Chapter
11] (see also [Jo, 1.4.15]). 

 Let $R_q[G]$ denote the quantized function algebra of
 the connected, simply connected algebraic
Lie group $G$ with Lie algebra ${\bf g}$. (See [Jo, Section 9.1]
for a precise definition.)  Note that up to a finite group,
$R_q[G]$ is the Hopf dual of $U$. Furthermore, $R_q[G]$ satisfies
a Peter-Weyl theorem ([Jo, 9.1.1 and 1.4.13]). That is, there is an isomorphism of $U$
bimodules $$R_q[G]\cong\bigoplus_{\lambda\in P^+(\pi)}
L(\lambda)\otimes L(\lambda)^*.\leqno{(3.1)}$$ Here $L(\lambda)$
is the (left) finite-dimensional simple $U$ module with highest
weight $\lambda$ contained in the set  $P^+(\pi)$   of dominant
integral weights. Moreover, $L(\lambda)^*$ is thought of as a
right $U$ module. Thus, the right $U$ module action on $R_q[G]$
comes from the action of $U$ 
 on $L(\lambda)^*$, while the left action comes from the action of
 $U$ on $L(\lambda)$.

Given a left coideal  $I$ of $U$ and a (left) $U$ module $M$, a (left) 
invariant is an $m\in M$ such that $am=\epsilon(a)m$ for all $a\in I$.
Write   $M_l^I$ for the collection of all left invariants in $M$. Equivalently, $M^I_l$
is equal to the elements of $M$ annihilated (on the left) by the
augmentation ideal of $I$. Consider the special case
where $I$ is the
quantum analog of the enveloping algebra of a Lie subalgebra of
${\bf g}$ corresponding to a subgroup $H$ of $G$.  Then $R_q[G]_l^I$ is
often written as $R_q[G/H]$.  In particular, $R_q[G/H]$ is thought
of  as the quantized function algebra on the quotient space $G/H$.

\begin{theorem}  For any left coideal $I$ of $U$,
$R_q[G]_l^I$ is a left coideal subalgebra of $R_q[G]$.
\end{theorem}

\noindent {\bf Proof:}  Let $\phi,\phi'$ be elements of
$R_q[G]_l^I$ and $r$ an element of $I$.   We first show that
$\phi\phi'$ is also in $R_q[G]_l^I$. To see this, consider
$$\eqalign{r\cdot (\phi\phi')=&\sum (r_{(1)}\cdot \phi)(r_{(2)}\cdot
\phi')\cr &=\sum (r_{(1)}\cdot \phi)\epsilon(r_{(2)})\phi'= (r\cdot
\phi)\phi'=\epsilon(r)\phi\phi'\cr}$$

We now check the coideal condition.  One can show using the
precise definition of the action of $U$ on $R_q[G]$ and the
coalgebra structure of $R_q[G]$ that  $${\it\Delta}(r\cdot
\phi)=(1\otimes r){\it\Delta}(\phi)=\sum\phi_{(1)}\otimes r\cdot
\phi_{(2)}.\leqno{(3.2)}$$  On the other hand, $${\it\Delta}(r\cdot
\phi)={\it\Delta}(\epsilon(r)\phi)=\sum \phi_{(1)}\otimes
\epsilon(r)\phi_{(2)}.\leqno{(3.3)}$$  Since we can choose the $\phi_{(1)}$ to
be linearly independent, (3.2) and (3.3) force  $r\cdot
\phi_{(2)}=\epsilon(r)\phi_{(2)}$.   Thus each $\phi_{(2)}\in
R_q[G]_l^I$. $\Box$

\medskip

   In [KS, Chapter 11.6],
 a quantum homogeneous space associated to $R_q[G]$ is defined up 
to isomorphism as  a one-sided coideal subalgebra of $R_q[G]$. (Note that quantum 
homogeneous spaces are actually quantum analogs of the algebra of regular 
functions on  classical homogeneous spaces.) Thus the 
theorem above shows there is a  left quantum
homogeneous space, $R_q[G]_I^I$, associated to each left coideal subalgebra $I$.  
Using the Peter-Weyl decomposition   (3.1), we  obtain the following nice description of 
$R_q[G]_l^I$.     $$R_q[G]_l^I\cong\bigoplus_{\lambda\in P^+(\pi)}
L(\lambda)^I_l\otimes L(\lambda)^*.\leqno{(3.4)}$$

Now $R_q[G]_l^I$ is the set of   left $I$ invariants of $R_q[G]$.
We may similarly define the set of right $I$ invariants
$R_q[G]^I_r$. Using the fact that the right action satisfies
${\it\Delta}(\phi\cdot r)={\it\Delta}(\phi)\cdot(r\otimes  1)$,
the same argument as in the proof of Theorem 3.1 shows that
$R_q[G]^I_r$ is a right coideal subalgebra of $R_q[G].$  One may
also study the set of bi-invariants $R_q[G]^I_{bi}=R_q[G]^I_l\cap
R_q[G]^I_r$. As above,  $R_q[G]^I_{bi}$ is a
subalgebra of $R_q[G]^I$. However, it is not a  coideal.

\section{Generators and Filtrations} 

Consider the  Hopf algebra $U(L)$, the universal 
enveloping algebra of a complex Lie algebra $L$. Since $U(L)$ is 
cocommutative, the one-sided coideal subalgebras of $U(L)$ are 
exactly the subbialgebras of $U(L)$.  It is very easy to 
understand the coideal subalgebras of $U(L)$.  Indeed,the next 
observation is   well known.  It follows from say 
 [Mo, Theorem 5.6.5] 
  and 
the fact that every subcoalgebra of $U(L)$ is connected ([Mo,  Definition 
5.1.5 and Lemma 5.1.9]). (Theorem 5.6.5 of [Mo] is stated for Hopf 
algebras, however, the proof also works for bialgebras.)

\begin{theorem}
The set of (left) coideal subalgebras of $U(L)$ is the set of  
enveloping 
algebras $U(L')$ of Lie subalgebras $L'$ of $L$.  
\end{theorem}

An immediate consequence of the above result is that any coideal 
subalgebra of $U(L)$ is generated by elements of the underlying Lie 
algebra $L$.   We would like to obtain a similar result for coideal 
subalgebras of the quantized enveloping algebra.
However, passing to the quantum case, the situation becomes more complicated.  
Indeed the coalgebra structure is not cocommutative for quantized 
enveloping algebras.  So the set of coideal subalgebras of the 
quantized enveloping algebra is much larger than the set of 
subbialgebras. By analyzing and deepening  Lemma 1.3 and studying the 
comultiplication of $U$, 
we are able to obtain detailed information 
about coideal subalgebras and their generators.

 The next result is known as well. It 
 describes the coideal subalgebras of a group algebra.
 
 \begin{lemma}  Let $I$ be a (left) coideal subalgebra of the group 
 algebra of the group $G$.   Then 
 $I\cap G$ is a semigroup and $I\cap kG$  is spanned by $I\cap G$ as a vector 
 space. 
 \end{lemma}

\medskip
We introduce two subalgebras of $U$ which are similar to $U^-$ and 
$G^+$.  
   Let $U^+$ be the subalgebra of $U$ generated by $x_1,\dots,
 x_n$ and $G^-$ be the subalgebra of $U$ generated by
 $y_1t_1,\dots, y_nt_n$. Once again, we have that
 $U^+$ and $G^-$ are a direct sum of their weight spaces.  
We may replace $U^-$ by $G^-$ and $G^+$ 
  by $U^+$ to obtain the following version of the triangular 
  decomposition.
 $$U\cong G^-\otimes U^o\otimes U^+.\leqno{(4.1)}$$
    In this section, we  show how to break up a coideal subalgebra 
 into   three parts corresponding to coideal subalgebras of $G^-$, $U^o$, and 
 $U^+$ respectively.   First, however, we obtain basic properties of  
 coideal subalgebras 
 of $G^-$ and $U^+$.  

 Using the formulas for comultiplicaton  (1.8), (1.9), and (1.10), it is straightforward 
 to check that $G^-$ and $U^+$ are  
 left coideal subalgebras of $U$. 
 Consider now an 
arbitrary coideal subalgebra $J$ of $U$ which is either a subset of 
$G^-$ or $U^+$.  Note that if $J$ is also an
 $\ad T$ module, then $J$ can be written as a direct sum of its weight 
 spaces.  We obtain a nice result on the generators of $J$ analogous 
to Theorem 4.1 when $J$ is an $\ad T$ module.
 
\begin{lemma}  Let $J$ be an $\ad T$ submodule and a coideal subalgebra of 
$G^-$ (resp. $U^+$).  Then there exists a subset 
$\Delta'$  of $ \Delta^+$   and weight 
vectors $\tilde f_{-\gamma}$ of weight $-\gamma$, $\gamma\in {\Delta'}$ (resp. $\tilde 
e_{\gamma}$ of weight $\gamma\in {\Delta'}$)  which generate $J$ as 
an algebra.  Moreover, each  $\tilde f_{-\gamma}$ (resp. $\tilde e_{\gamma}$) specializes to the root 
vector $f_{-\gamma}$ (resp. $e_{\gamma}$) as $q$ goes to $1$. 
\end{lemma}

\noindent
{\bf Proof:}   
Note that the  weight spaces of $G^-$ are finite-dimensional.  Hence 
$J$ has 
finite-dimensional weight spaces.  Let $\bar J$ denote the 
specialization of $J$ as $q$ goes to $1$.
Consider a weight space $J_{\mu}$ 
of $J$.   We have that $\hat J_{\mu}=\hat U\cap J_{\mu}=\hat G^-\cap 
J_{\mu}$. Also, 
$G^-$ is a free $A$ module and $A$ is a principal ideal domain
with unique maximal ideal generated by $(q-1)$. Hence one 
can find a basis for $J_{\mu}$ which is a subset of $\hat J_{\mu}$ 
and remains linearly independent modulo $(q-1)\hat U$.   In 
particular, the specialization of this basis as $q$ goes to $1$ is a 
basis  for $\bar J_{\mu}$.  Hence  
 the weight spaces of $\bar J$ have the same dimension as the 
weight spaces of $J$. 

 Note that the 
comultiplication of $U$ specializes to the comultiplication of $U({\bf 
g})$.  Hence $\bar J$ is a coideal subalgebra
of $U({\bf n}^-)$.  By Theorem 4.1,   $\bar J$ is an enveloping algebra of a Lie 
subalgebra, say ${\bf a}$, of ${\bf n}^-$.  Now  $\bar J$ is a direct 
sum of its weight spaces.   Hence there exists a 
subset $\Delta'$ of $\Delta^+$ such that the set 
$\{f_{-\gamma}|\gamma\in \Delta'\}$ is a basis of ${\bf a}$. 
Thus for each $\gamma\in \Delta'$,
we can find a vector $\tilde f_{-\gamma}$ of weight $-\gamma$ in $J$ 
such the image of $\tilde f_{-\gamma}$ under specialization is $f_{-\gamma}$. 
Write $\Delta'=\{\gamma_1,\dots, \gamma_m\}$.   A standard argument 
shows that the set $${\cal B}_{\eta}=\{f_{-\gamma_1}^{i_1}\dots 
f_{-\gamma_m}^{i_m}|i_j\in{\bf N}  {\rm \ for\ }  1\leq j\leq m 
{\rm \ and \ }i_1\gamma_{1}+\dots +i_m\gamma_m=\eta \}$$ is a basis for the 
$-\eta$ weight space of 
$U({\bf a})$. Furthermore ${\cal B}=\cup_{\eta}{\cal B}_{\eta}$  
is a basis for  $U({\bf a})$.   
Since the dimensions of the $-\eta$ weight spaces of $U({\bf a})$ and 
$J$ agree,  the corresponding set $\tilde{\cal B}_{\eta}$ with 
$\tilde f$ playing the role of $f$ is a basis of $J_{\eta}$.
  Thus the set $\tilde{\cal B}=\cup_{\eta} {\tilde {\cal  
B}}_{\eta}$ is a basis for $J$. 
It follows that  the $\tilde f_{-\gamma}$, $\gamma\in \Delta'$  
generate $J$ as an algebra. 

The same analysis applies to coideal 
subalgebras of $U^+$. 
$\Box$

\medskip

 The  direct sum decomposition (1.13) can be made finer using weight
spaces. It is well known that  the
set of weights of
 $G^+$ and $U^+$  equals $Q^+(\pi)$ and the set of weights of 
 $G^-$ and $U^-$ equals $Q^-(\pi)$.   
 Thus there are   direct sum decompositions 
$$U=\oplus_{\lambda, \mu}U^-_{-\lambda}G^+_{\mu}U^o{\rm \ and \ }
U=\oplus_{\lambda, \mu,t}U^-_{-\lambda}G^+_{\mu}t\leqno{(4.2)}$$ where
$\lambda$ and $\mu$ run over elements of $Q^+(\pi)$ and $t$ runs over 
elements in $ T$.
Let $\pi_{\lambda,\mu}$ be the projection of $U$ onto the subspace
$U^-_{-\lambda}G^+_{\mu}U^o$. Write $[\lambda,\mu]$ for a typical
element in $Q(\pi)\times Q(\pi)$ (so as to avoid confusion with
the Cartan inner product).

Consider elements $c\in U^-_{-\lambda}$ and $d\in G^+_{\mu}$. 
It follows from the definition of the comultiplication map on the 
generators of $U$ ((1.8), (1.9), and (1.10)) that  
\begin{enumerate}
\item[(4.3)]$(\pi_{\lambda,\mu}\otimes Id){\it\Delta}(cd)=cd\otimes 
\tau(-\lambda-\mu)$
\item[(4.4)]$(\pi_{\lambda,0}\otimes 
Id){\it\Delta}(cd)=\sum c\otimes \tau(-\lambda)d$
\item[(4.5)]$(\pi_{0,\mu}\otimes 
Id){\it\Delta}(cd)=\sum d\otimes c\tau(-\mu).$
\end{enumerate}

In the next theorem, we consider coideal subalgebras of $U$ which 
behave rather nicely in terms of the second decomposition in (4.2).

\begin{theorem} Let $I$ be a left coideal subalgebra of $U$ such that
$$I=\sum_{\lambda, \mu,t}(I\cap  U^-_{-\lambda}G^+_{\mu}t)\leqno{(4.6)}$$ and $I\cap T$ 
is a group.   Then 
$I\cap G^-$, $I\cap U^o$, and $I\cap U^+$ are $\ad T$ submodules and 
left coideal subalgebras of $I$.   Moreover, the 
multiplication map induces an isomorphism
$$I\cong (I\cap G^-)\otimes (I\cap U^+)\otimes (I\cap 
U^o)$$ of vector 
spaces.
\end{theorem}

\noindent
{\bf Proof:}  Since $I,G^-, U^o, $ and $U^+$ are all left coideal 
subalgebras, so are $I\cap G^-$, $I\cap U^o$, and $I\cap U^+$.
Note that every element in 
 $U^-_{-\lambda}G^+_{\mu}t$ is a weight vector of weight $-\lambda+\mu$.
 Thus   $I$ is a 
direct sum of its weight spaces and $I$ is an $\ad T$ module.   
It follows that $I\cap G^-$, $I\cap U^o$, and $I\cap U^+$ are all $\ad 
T$ modules.  

The triangular decomposition of $U$ (4.1) ensures that the 
multiplication map induces an injection 
$$ (I\cap G^-)\otimes (I\cap U^+)\otimes (I\cap U^o)\rightarrow I$$
of vector spaces.  We obtain an isomorphism by showing that each 
element of $I$ is contained in $(I\cap G^-)(I\cap U^+)(I\cap U^o)$.

Recall the direct sum decomposition of $I$ given in Lemma 1.3.
 Let $b$ be an element of $I\cap  U^-_{-\lambda}G^+_{\mu}t$
where $t\in T$.
There exists $c_i\in U^-_{-\lambda}$ and $d_i\in G^+_{\mu}$
  so that $b=\sum_ic_id_it$. We may
further assume that  $\{c_i\}$ and $\{d_i\}$ are each linearly
independent sets. By (1.8) and (4.3),
$(\pi_{\lambda,\mu}\otimes Id){\it\Delta}(b)=b\otimes 
\tau(-\lambda-\mu)t$.   Hence $\tau(-\lambda-\mu)t$ is an element of 
$I\cap T$.   Since $I\cap T$ is a group, $\tau(\lambda+\mu)t^{-1}$ is 
also contained in $I\cap T$.

Equation (4.4)  implies that $$(\pi_{\lambda,0}\otimes 
Id){\it\Delta}(b)=\sum_i c_it\otimes \tau(-\lambda)d_it.$$
Hence each $ \tau(-\lambda)d_it\in I$. Recall that $d_i$ is a weight 
vector of weight $\mu$ in $G^+$.   Thus, multiplying  $ \tau(-\lambda)d_it$
by $\tau(\lambda+\mu)t^{-1}$ yields that $d_i\tau(\mu)$ is an element 
of $U^+_{\mu}\cap I$.  
Similarly,   (4.5) ensures that
$$(\pi_{0,\mu}\otimes 
Id){\it\Delta}(b)=\sum_i d_it\otimes c_i\tau(-\mu)t.$$
So $c_i\tau(-\mu)t\in I$ and hence 
$c_{i}\tau(-\mu)t\tau(\lambda+\mu)t^{-1}=c_i\tau(\lambda)$ is an 
element of $I\cap G^-$.   Therefore,
$$\eqalign{
b=&\sum_ic_id_it=\sum_iq^{(-\lambda,\mu)}(c_i\tau(\lambda))(d_i\tau(\mu))\tau(-\lambda-\mu)t\cr
\in&(I\cap G^-)(I\cap U^+)(I\cap U^o).\Box\cr}  $$

\medskip

Let $I$ be a left coideal subalgebra of $U$  such that $I\cap T$ is a 
group and $I$ satisfies (4.6).  Then   Theorem 4.4 combined with 
Lemmas 4.2 and 4.3 imply that $I$
is generated by $I\cap T$ and quantum analogs of root vectors in $G^-$ 
and $U^+$.   This description of the generators can be thought of as 
an
analog of Theorem 4.1 for these special coideal subalgebras.  Below, 
we generalize these results to other coideal subalgebras by introducing
filtrations and associated gradings of $U$.

\medskip
\noindent
{\bf    Filtration I}
 
\medskip
\noindent 
 Define the
  filtration ${\cal F}$ on $U$   using the degree function:
  $$\deg x_it_i^{-1} = \deg y_i =1 \qquad \deg t_i=-1$$  for all $1\leq
  i\leq n$.   Write $\gr_{\cal F} U$ for the associated
  graded algebra of this filtration.      This filtration is
  invariant under the adjoint action and  used to understand the
  locally finite part of $U$ (see [JL2, Section 2.2]).   
  (It should be noted that the quantized enveloping 
  algebra is defined in a different though equivalent manner in [JL2].   So
  the $x_i$ (resp. $t_i$) in this paper corresponds to $x_it_i$ 
  (resp. $t_i^2$) in [JL2].  Furthermore the   degree of 
  an element as defined in [JL2] is twice the   degree  of 
  the corresponding element  given here.)
  
 Given
$\gamma=\sum_{\alpha_i\in\pi}m_i\alpha_i$, set
$\ht(\gamma)=\sum_im_i$. 
 Note that any nonzero element of 
$U^-_{-\lambda}G^+_{\mu}$ has degree ${\rm ht}(\lambda+\mu)$. 
Let $x\in U$ and set 
$\supp(x)=\{[\lambda,\mu]|\pi_{\lambda,\mu}(x)\neq 0\}$. 
Further, for   $x$   an element of $U^-G^+t$ for some $t\in T$, set
 $${\rm max}_{\ht}(x)=\{[\lambda,\mu]| 
[\lambda,\mu]\in \supp(x){\rm \ and \ }\ht(\lambda+\mu)=\deg(x)-\deg(t)\}.$$   
 The 
next lemma connects the filtration $\cal F$ with the height function.

\begin{lemma}  Let $I$ be a left coideal   of $U$ and let
$b$ be an element of $I\cap U^-G^+t$ for some $t\in T$.     Then 
$$b=\sum_{\{[\lambda,\mu]|\
[\lambda,\mu]\in\max_{\ht}(b)\}}\pi_{\lambda,\mu}(b)+{\rm \ lower \ degree\ 
terms}.$$  
\end{lemma}

\noindent
{\bf Proof}:    The lemma follows from the fact that 
$\deg\pi_{\lambda,\mu}(b)=\deg b$ if and only if $[\lambda,\mu]\in 
\max_{\ht}(b)$. $\Box$

\medskip

By induction on $\ht(\lambda+\mu)$
and  the definition of the comultiplication (1.8), (1.9), and (1.10), we have the 
following:
 $${\it\Delta}(U^-_{-\lambda}G^+_{\mu})\subset 
\sum_{\gamma+\beta=\lambda,\alpha+\xi=\mu}U_{-\gamma}^-G^+_\alpha\otimes 
U_{-\beta}^-G_{\xi}^+\tau(-\gamma-\alpha).\leqno{(4.7)}$$  

Consider a subset $S$ of $Q^+(\pi)\times Q^+(\pi)$. Set $|S|$ equal 
to the number of elements in $S$.  We call $S$ ${\it 
transversal}$ if whenever both $[\lambda,\mu]$ and $[\lambda',\mu']$ are in $S$ 
and $[\lambda,\mu]\neq [\lambda',\mu']$ then $\lambda\neq \lambda'$ and $\mu\neq \mu'$. 
Now assume that $b\in U^-G^+t$ and that $\max_{\ht}(b)$ is transversal. 
Given $[\lambda,\mu]\in {\max}_{\ht}(b)$, find
 $c_i\in U^-_{-\lambda}$ and $d_i\in G^+_{\mu}$
  such that $\pi_{\lambda,\mu}(b)=\sum_ic_id_it$.  As in the proof of 
  Theorem 4.4, we may 
  further assume that $\{c_i\}$ and $\{d_i\}$ are each linearly
independent sets.  It follows from (4.7), (4.3), (4.4), and (4.5) 
that
\begin{enumerate}
\item[(4.8)]$(\pi_{\lambda,\mu}\otimes Id){\it\Delta}(b)=\sum_ic_id_it\otimes 
\tau(-\lambda-\mu)t$ 
\item[(4.9)]$(\pi_{\lambda,0}\otimes 
Id){\it\Delta}(b)=\sum_i c_it\otimes (\tau(-\lambda)d_it+$ terms of 
lower degree)
\item[(4.10)]$(\pi_{0,\mu}\otimes 
Id){\it\Delta}(b)=\sum d_it\otimes (c_i\tau(-\mu)t$ + terms of lower 
degree)
\end{enumerate}

 A consequence of    the next lemma  is that any left coideal subalgebra which is also an $\ad T$ 
module has a basis ${\cal B}$ such that $\max_{\ht}(b)$ is transversal 
for each $b\in {\cal B}$.  This 
in turn is used to generalize Theorem 4.4. 

\begin{lemma} Let $b\in U$ be a weight vector.   Then $\max_{\ht}(b)$ is 
transveral.
\end{lemma}

\noindent
{\bf Proof:} 
 Fix $\eta$ and  let 
 $b$ 
 be an element in $U$ of weight $\eta$.   Note 
 that $\pi_{\lambda,\mu}(b)\neq 0{\rm\ implies \ that
 }-\lambda+\mu=\eta.$ Now assume that both $[\lambda,\mu]$
 and $[\lambda',\mu']$ are in $\supp(b)$.   Hence
 $-\lambda+\mu=-\lambda'+\mu'$.   Thus $\lambda=\lambda'$ if and only 
 if $\mu=\mu'$.  In particular, $\supp(b)$ is transversal.   The lemma 
 now follows 
 from the fact that $\max_{\ht}(b)$ is a subset of $\supp(b)$.
$\Box$

\medskip
Given a left coideal subalgebra $I$ of $U$, set $I^-_{\eta}$ equal to the 
subset of $G^-$ such that $I\cap G^-\tau(\eta)=I^-_{\eta}\tau(\eta)$.  
Similarly, set $I^+_{\eta}$ equal to the 
subset of $U^+$ such that $I\cap U^+\tau(\eta)=I^+_{\eta}\tau(\eta)$.
The following result can be thought of as an analog of Theorem 4.4 for coideal 
subalgebras which admit an $\ad T$ module structure.

\begin{theorem} Let $I$ be a left coideal subalgebra and $\ad 
T$ submodule of $U$. Then $$\gr_{\cal F}I\subset \sum_{\{\eta|\tau(\eta)\in I\cap 
T\}}\gr_{\cal F}I^-_{\eta}I^+_{\eta}\tau(\eta).$$ 
\end{theorem}

\noindent
{\bf Proof:}  
 Let $b$ be a weight vector of $I$ which is 
also contained in $I\cap U^-G^+\tau(\beta)$ for some $\tau(\beta)\in T$.
  By Lemma 4.6,  ${\rm max}_{\rm ht}(b)$ is transversal. We prove the 
  theorem when ${\rm max}_{\rm ht}(b)$ contains exactly one element 
  $[\lambda,\mu]$.   The same argument works in general. We argue as 
  in the proof of Theorem 4.4. 
Find $c_{i}\in U^-_{-\lambda}$ and $d_{i}\in G^+_{\mu}$
  so that $\pi_{\lambda,\mu}(b)=\sum_{i}c_{i}d_{i}\tau(\beta)$. We may
further assume that  $\{c_{i}\}$ and $\{d_{i}\}$ are each linearly
independent sets. By 
  our assumption on $\max_{\ht}(b)$ and Lemma 4.5,
$$
b=\sum_{i}c_{i}d_{i}\tau(\beta)+{\rm\ lower\ degree\ terms.\ }$$ 

Set $\eta=-\lambda-\mu+\beta$. By (4.8), 
 $\tau(\eta)$ is in $I$.   Now 
(4.9) implies that there exist  elements 
$\tau(-\lambda)D_{i}\tau(\beta)\in I$ such that
$D_{i}=d_{i}+$ lower degree terms and 
$$(\pi_{\lambda,0}\otimes Id){\it\Delta}(b)=\sum_ic_{i}\otimes 
\tau(-\lambda)D_{i}\tau(\beta).$$   Note that 
(4.7) ensures that 
$D_{i}-d_{i}$ is an element of $U^-G^+\tau(-\lambda+\beta)$.  Also, 
$d_i$ is in $G^+_{\mu}\tau(-\lambda+\beta)$.   Thus $d_i$ has degree 
$\ht(\mu+\lambda-\beta)$.  By Lemma 4.5, 
$[\xi,\gamma]\in \supp(D_{i}-d_{i})$ implies that 
$\ht(\xi+\gamma)<\ht(\mu)$. Since $\xi$ is in $Q^+(\pi)$, we also 
have  $\ht(-\xi+\gamma)<\ht(\mu)$ and thus
$-\xi+\gamma$ is not equal to $\mu$.  Therefore, for each 
$[\xi,\gamma]\in \supp(D_i-d_i)$, $\pi_{\xi,\gamma}(D_i-d_i)$ 
has weight $-\xi+\gamma$ which is different from the weight $\mu   $ 
of  
$d_i$. Since $I$ is an $\ad T$ module, it  follows that the $\mu$ 
weight term of $\tau(-\lambda)D_{i}\tau(\beta)$, namely  
$\tau(-\lambda)d_{i}\tau(\beta)$, is contained 
in $I$. Hence
$d_{i}\tau(\mu)\tau(\eta)\in I\cap 
U^+\tau(\eta)$ and $d_{i}\tau(\mu)\in I^+_{\eta}$.    A 
similar argument shows that $c_{i}\tau(\lambda)\in I^-_{\eta}$. 
Therefore
$${\gr_{\cal F}b=\gr_{\cal F}\sum_{i}c_{i}d_{i}\tau(\eta) =
\gr_{\cal F}\sum_iq^{-(\lambda,\mu)}(c_i\tau(\lambda))(d_i\tau(\mu))\tau(\eta)}$$ is 
an element of
$\gr_{\cal F}I^-_{\eta}I^+_{\eta}\tau(\eta)$.   
$\Box$

\bigskip
\noindent
{\bf Filtration II}

\medskip
\noindent
Order the  set ${\bf N}\times {\bf N}$ lexicographically from left to 
right. 
Define a filtration on $U$ by $${\cal G}_{m,n}(U)=\{u\in
U|\ (\ht(\lambda),\ht(\mu))\leq (m,n){\rm \ for \ all\ }[\lambda,\mu]\in \supp(u)\}.$$ 
The    associated graded algebra for this 
filtration is defined by setting $$\gr_{\cal
G}^{m,n}(U) ={\cal
G}_{m,n}(U)/\sum_{(m',n')<(m,n)}{\cal
G}_{m,n}(U)$$ and $$\gr_{\cal 
G}(U)=\bigoplus_{m,n}\gr_{\cal G}^{m,n}(U).$$

Given a subset $S$ of $Q^+(\pi)\times Q^+(\pi)$, set $||S||_1={\rm 
max}_{[\lambda,\mu]\in S}\{\ht(\lambda)\}$.
 We can define a bidegree: for $x$ in $U$, we say that
$\bideg(x)= (m,n)$ if $(m,n)$ is the 
smallest element of ${\bf N}\times {\bf N}$ such that $x\in {\cal 
G}_{m,n}(U)$. Set ${\max}(x)=\{[\lambda,\mu]|[\lambda,\mu]\in 
\supp(x)$ and $\bideg(x)=(\ht(\lambda),\ht(\mu))\}$. 
Now consider an element $b\in U^-G^+t$ for some $t\in T$.    The inclusion 
(4.7) 
implies the following variation of (4.3): 
$$(\pi_{\lambda,\mu}\otimes Id)(b)=\pi_{\lambda,\mu}(b)\otimes 
t\tau(-\lambda-\mu){\rm \ for \ all\ }[\lambda,\mu]\in 
\max(b).\leqno{(4.11)}$$

 \begin{lemma}  Let $I$ be a left coideal subalgebra such that $I\cap T$ 
 is a group. Then $I$ has a basis ${\cal B}$ such that for each $b\in 
{\cal  B}$, $\max(b)$ is transversal.
 \end{lemma}
 
 \noindent
 {\bf Proof:}
Recall that $I$ is a direct sum of the subspaces $I\cap U^-G^+t$.  
Let $C=\{x\in I|\max(x)$ is transversal$\}$.  
It is enough to show that for each $t\in T$,  every element of $I\cap 
U^-G^+t$ is contained in the 
span of $C$.  Consider  $b\in I\cap U^-G^+t$.  We prove this result 
under the additional assumption  that $\max(b)$ 
  consists of exactly two elements $[\lambda,\mu]$ and $[\lambda,\mu']$. 
  A similar argument works in the general case using induction on $|{\max}(b)|$ and 
  $||{\max}(b)||_1$. 
 Note that $$b=\pi_{\lambda,\mu}(b)+\pi_{\lambda,\mu'}(b)+{\rm \ terms 
  \ of \ lower \ bidegree}.$$ By (4.11),   $ t \tau(-\lambda-\mu)$
 and  $ t \tau(-\lambda-\mu')$
are both elements of   the group $I\cap T$. Hence
 $\tau(\mu-\mu')$ is 
  contained in $I\cap T$.   Consider the element 
  $$\eqalign{b'=\tau(\mu-\mu')b\tau(\mu-\mu')^{-1}
  =&q^{(-\lambda+\mu,\mu-\mu')}\pi_{\lambda,\mu}(b)+
  q^{(-\lambda+\mu',\mu-\mu')}\pi_{\lambda,\mu'}(b)\cr+&{\rm \ terms 
  \ of \ lower \ bidegree}.\cr}$$   Now $(\mu-\mu',\mu-\mu')$ is positive since
  $(\ ,\ )$ is positive definite on $Q(\pi)$.  Hence $q^{(-\lambda+\mu,\mu-\mu')} 
  \neq q^{(-\lambda+\mu',\mu-\mu')}$.  Taking linear combinations of 
  $b$ and $b'$, it follows that there exists  $b_1$ and $b_2$ in $U^-G^+t\cap I$
  such that $\{[\lambda,\mu]\}={\rm max}(b_1)$ and  
  $\{[\lambda,\mu']\}=
  {\rm max}(b_2)$. In particular, both ${\rm max}(b_1)$ and ${\rm 
  max}(b_2)$ are transversal and $b$ is a linear combination of $b_1$ 
  and $b_2$. $\Box$
  
 \medskip
  Now assume that $b\in U^-G^+t$ and that $\max(b)$ is transversal.
  We have versions of (4.4) and (4.5) similar to (4.9) and (4.10) 
  in the discussion of the first filtration.   
Given $[\lambda,\mu]\in {\max}(b)$, find
 $c_i\in U^-_{-\lambda}$ and $d_i\in G^+_{\mu}$
  so that $\pi_{\lambda,\mu}(b)=\sum_ic_id_it$  and that $\{c_i\}$ and $\{d_i\}$ are each linearly
independent sets.  It follows from (4.7), (4.4), and (4.5)
that
\begin{enumerate}
\item[(4.12)]$(\pi_{\lambda,0}\otimes 
Id){\it\Delta}(b)=\sum_i c_it\otimes (\tau(-\lambda)d_it+$ terms of 
lower bidegree)
\item[(4.13)]$(\pi_{0,\mu}\otimes 
Id){\it\Delta}(b)=\sum d_it\otimes (c_i\tau(-\mu)t$ + terms of lower 
bidegree).
\end{enumerate}

The filtration ${\cal G}$ restricts to filtrations on the 
subalgebras $G^+$, $U^-$ and $U^o$.  Indeed, 
$U^o={\cal G}_{0,0}(U)$ and so $\gr_{\cal G}U^o\cong U^o$
as algebras.
Upon restriction to  $G^+$, the filtration ${\cal G}$ becomes 
filtration by the degree function associated to the first 
filtration ${\cal F}$.  
   The subalgebra of $G^+$ satisfies exactly the same relations as 
   $U^+$.   In particular,  the only relations satisfied 
by  the elements of $G^+$ are the quantum Serre relations (1.7) (see 
the discussion in Section 7 concerning (7.18))
 which are
homogeneous with respect to degree.   Hence we have an algebra isomorphism    
$\gr_{\cal G}G^+\cong G^+$.  A similar argument shows that $\gr_{\cal G}U^-\cong U^-$.  
Since the elements in $T$ have bidegree $(0,0)$, we further have that
$\gr_{\cal G}G^-\cong G^-$ and  $\gr_{\cal G}U^+\cong U^+$. 
For the rest of the paper, we will often identify $\gr_{\cal G}G^-$
with $ G^-$,  $\gr_{\cal G}U^+$ with $U^+$, and $\gr_{\cal G}U^o$ with 
$U^o$.

    Now the images of $ x_i$ and $y_j $
commute with each other inside the associated graded algebra of $U$ with 
respect to ${\cal G}$. (See 
relation (1.4) of $U$.)    It
follows that the image of $U^-U^+$ in the associated graded
algebra is isomorphic to the tensor product $U^-\otimes U^+$ as an
algebra. If we replace $U^-$ by $G^-$, the images of the  elements
$x_i$ and $y_jt_j$ do not commute.   However, they do commute up
to a power of $q$. Thus the image of $G^-U^+$ in the associated
graded algebra can be thought of as a $q$ form of the tensor
product which we write as $G^-\otimes_q U^+$.  

The group algebra $U^o$ 
acts  on weight 
vectors  by $\tau(\lambda)\cdot 
a_{\mu}=\tau(\lambda)a_{\mu}\tau(\lambda)^{-1}=q^{(\lambda,\mu)}a_{\mu} $
for $a_{\mu}\in U_{\mu}$.   Thus we obtain the following algebra 
isomorphism using a smash product construction:
$$\gr_{\cal G}(U)\cong U^o\#(G^-\otimes_q 
U^+).\leqno{(4.14)}$$  (Compare this with a similar result for a different 
filtration in [Jo, 7.4.7].)

\begin{theorem}  Let $I$ be a left coideal subalgebra such that 
$I\cap T$ is a group.  Then  
$$\gr_{\cal G}(I)\cong (I\cap U^o)\#((\gr_{\cal G}(I)\cap G^-)\otimes_q 
((\gr_{\cal G}(I)\cap U^+)).$$  
\end{theorem}

\noindent{\bf Proof:}  The proof is a graded version of the proof of Theorem 4.4. 
Using Lemma 1.3 and Lemma 4.8, we can find a basis ${\cal B}$ of $I$ such 
that ${\cal B}=\bigcup_t({\cal B}\cap U^-G^+t)$ and $\max(b)$ is 
transversal for each $b\in {\cal B}$.  By (4.14), $$(I\cap U^o)\#((\gr_{\cal G}(I)\cap G)^-\otimes_q 
((\gr_{\cal G}(I)\cap U^+))$$ is isomorphic to a subalgebra of 
$\gr_{\cal G}I$.  To show this subalgebra is all of $\gr_{\cal G}I$ it 
is sufficient to show that each element of ${\cal B}$ is 
contained in  $((\gr_{\cal G}(I)\cap G^-)
((\gr_{\cal G}(I)\cap U^+))\gr_{\cal G}(I\cap U^o)$. 

Let $t$ be an element of $ T$ and let $b$ be an element of ${\cal B}\cap 
U^-G^+t$.  Choose $[\lambda,\mu]\in \max(b)$. 
There exists $c_i\in U^-_{-\lambda}$ and $d_i\in G^+_{\mu}$
  so that $\pi_{\lambda,\mu}(b)=\sum_ic_id_it$ and the $\{c_i\}$ and $\{d_i\}$ are each linearly
independent sets. Using (4.11), (4.12), and (4.13) and arguing as in 
the proofs of Theorem 4.4 and Theorem 4.7, $I$ contains
 $\tau(-\lambda-\mu)t$ and $\tau(\lambda+\mu)t^{-1}$ and elements 
 $\tilde d_i$ and $\tilde c_i$ such that
 $$\tilde d_i= d_i\tau(\mu)+{\rm \ terms \ of \ lower\ bidegree}$$
 and 
 $$\tilde c_i=c_i\tau(\lambda)+{\rm \ terms \ of \ lower \ bidegree}.$$
Note that $\gr_{\cal G}\tilde d_i\in \gr_{\cal G} (I)\cap U^+$ and $\gr_{\cal 
G}\tilde c_i\in \gr_{\cal G}(I)\cap G^-$. Set 
$$
b'=b-\sum_iq^{-(\lambda,\mu)}\tilde c_i\tilde 
d_i\tau(-\mu-\lambda)t.$$ Note that $b'$ is in $I$.  By construction,
$\pi_{\lambda,\mu}(b')=0$.  Thus
either $\max(b')=\max(b)-\{[\lambda,\mu]\}$ or the bidegree of $b'$ is 
strictly smaller than the bidegree of $b$.   The theorem now follows 
by induction on 
$|\max(b)|$ and the bidegree of $b$.
$\Box$

\medskip

Consider a left coideal 
subalgebra $I$ such that $I\cap T$ is a group. Given $x$ in $U$, set $\tip(x)=\sum_{[\lambda,\mu]\in 
\max(x)}\pi_{\lambda,\mu}(x)$.   The element $\tip(x)$ can be thought 
of as the highest bidegree term of $x$. Note that $\gr_{\cal G}(I)\cap \gr_{\cal G}(G^-)$ 
identifies with  $\tip(I)\cap G^-$  under the isomorphism $G^-\cong 
\gr_{\cal G}(G^-)$. Thus $\tip(I)\cap G^-$ is a subalgebra of $G^-$.  
Consider the elements $\tilde c_i$ in the proof of Theorem 4.9.   Note that 
each $\tip(\tilde c_i)$ is a weight vector.   In particular, it follows 
implicitly from the proof of Theorem 
4.9 
that $\tip(I)\cap G^-$ is spanned by weight vectors and hence is an $\ad 
T$ module.  One can further show using (4.7) that $\tip(I)\cap G^-$ 
is a left coideal of $G^-$.   Thus $\tip(I)\cap G^-$ is a left coideal 
subalgebra and $\ad T$ submodule of $G^-$.   Similarly, $\tip(I)\cap 
U^+$ is a left coideal subalgebra and $\ad T$ submodule of $U^+$.  
Thus combining Theorem 4.9 with Lemmas 4.2 and 4.3 yields the following.

\begin{corollary} Let $I$ be a left coideal subalgebra of $U$ such 
that  $I\cap T$ is a subgroup of $T$.   
Then there exists  subsets $\Delta'$ and $\Delta''$ of $\Delta^+$ such that $I$ is generated by elements  $
c_{-\gamma},\gamma\in\Delta'$;  $d_{\beta},\beta\in\Delta''$; and 
$I\cap T$.   Moreover each 
$\tip(c_{-\gamma})$ (resp. $\tip( d_{\beta})$) is a weight vector of 
weight $-\gamma$ (resp. $\beta$) which specializes to the root vector $f_{-\gamma}$
(resp. $e_{\beta}$) of $U({\bf 
g})$.
\end{corollary}

 \section{The locally finite part of $U$}

  One of the most important coideal subalgebras contained in the
  quantized enveloping algebra is the locally finite part,
 $F(U)$, defined by (2.1).  This subalgebra is studied extensively in
[JL1]  and [JL2] (see also [Jo]). Here we present some of the known
results about this algebra by directly showing that $F(U)$ is a 
coideal subalgebra of $U$. We will see some of the implications of Section 4 on
the structure of $F(U)$. 

 Recall that $F(U)$ is defined using the
quantum adjoint action in Section 2.
   It is helpful
to see
  how the generators of $U$ act via the adjoint action. In particular
$$( \ad \ y_i)b = y_i bt_i -by_it_i \quad (\ad
x_i)b=x_ib-t_ibt_i^{-1}x_i\quad   (\ad t_i)b=t_ibt_i^{-1}$$ for
all $1\leq i\leq n$ and $b\in U$.

\begin{theorem}
$F(U)$ is a left coideal subalgebra of $U$.
\end{theorem}

\noindent {\bf Proof:}   Let $b\in F(U)$. A straightforward
computation shows $$\eqalign{{\it\Delta}((\ad x_i)b) &= \sum x_ib_{(1)}\otimes
b_{(2)} -\sum t_ib_{(1)}t_i^{-1}x_i\otimes t_ib_{(2)}t_i^{-1}
\cr &+\sum t_ib_{(1)}\otimes (\ad x_i)b_{(2)}\cr}\leqno{(5.1)}$$
for each $1\leq i\leq n$.

We may write ${\it\Delta}(b)$ as a sum $\sum_{j=1}^sc_j\otimes b_j$ 
where the $b_j$, $1\leq j\leq s$, are linearly independent weight 
vectors in $U$ of weight $\lambda_j$ respectively.   Extend the 
standard partial ordering on the integral lattice $Q(\pi)$ to a total 
archimedean ordering. (This can be done by embedding $Q^+(\pi)$ in 
the nonnegative real numbers.)  We may further 
suppose that $\lambda_1\geq \lambda_2\geq \cdots\geq \lambda_s$. 
For  sake of simplicity,  assume that these inequalities are all 
strict.  (A similar argument works in general.) 

 Let $U_i$ be the subalgebra of $U$ 
generated by $x_i,y_i,t_i^{\pm 1}$.  Note that $U_i$ is isomorphic to 
$U_q({\bf sl}\, 2)$. We show below that each $b_j$ generates a 
finite-dimensional $\ad U_i$ module for $1\leq i\leq n$.  By [JL1, 
Theorem 5.9], this 
forces each $b_j$ to be an element of $F(U)$ (see also the proof of [JL1, 
Proposition 6.5]).
 
Suppose that $(\ad x_i)^mb=0$.  Using (5.1) and induction,  we obtain
$${\it\Delta}((\ad x_i)^mb)\in t_i^mc_1\otimes (\ad 
x_i)^mb_1+\sum_{\beta<\lambda_1+m\alpha_i}U\otimes U_{\beta}.$$
Hence $(\ad 
x_i)^mb_1 =0$. Choose $r$ such that 
$(m-1)\alpha_i+\lambda_1<(m+r)\alpha_i+\lambda_2$.  We further have that 
$${{\it\Delta}((\ad x_i)^{m+r}b)\in
t_i^{m+r}c_{2}\otimes (\ad x_i)^{m+r}b_{2}
+\sum_{\beta<\lambda_{2}+({m+r})\alpha_i}U\otimes 
U_{\beta}}.$$
Thus $(\ad x_i)^mb=0$ also implies that $(\ad 
x_i)^{m+r}b_2 =0$. By induction, it follows that there exists $M>0$ 
such that
$(\ad x_i)^{M}b_j=0$ for all $1\leq j\leq s$ and $1\leq i\leq n$.
    One obtains a similar property for the action of
each $\ad y_i$ on $b$. In particular, for each $1\leq i\leq n$ and 
each $1\leq j\leq s$, both $\ad y_i$ and $\ad x_i$ act nilpotently on 
$b_j$.   Since $b_j$ is a weight vector, it further follows that $\ad 
t_i$ acts semisimply on $b_j$.  Thus   $b_j$   
generates a finite-dimensional $\ad U_i$ module for all $1\leq i\leq n$ 
and all $1\leq j\leq s$. Therefore  each
$b_j\in F(U)$. $\Box$

\medskip
Set $T_{F}=T\cap F(U)$. 
It follows from Lemma 4.2 that
the  algebra generated by $T_F$ is equal to the intersection of
$U^o$ with $F(U)$.  By [JL1, 6.2],
$\tau(\lambda)\in F(U)$ if and only if $(\ad x_i)$ and $(\ad y_i)$ 
act nilpotently on $\tau(\lambda)$ for $1\leq i\leq n$.    
Furthermore, ([JL1, the proof of Lemma 6.1]) $s$ is the least positive integer such that 
$(\ad x_i)^s\tau(\lambda) =0$ and $(\ad y_i)^s\tau(\lambda) =0$    if and only if
$${{(\lambda,\alpha_i)}\over
{(\alpha_i,\alpha_i)}}=-s+1.$$ 
Thus, $\tau(\lambda)\in F(U)$ if and only if $(\lambda,\alpha_i)/
(\alpha_i,\alpha_i)$ is a nonpositive integer   for all $1\leq i\leq n$. 
In particular, $-\lambda/2$ is a dominant integral weight.   So  
$\tau(\lambda)$ is in $F(U)$ if and only if  
$\lambda$ is in $R(\pi):= 
Q(\pi)\cap -2P^+(\pi)$. (This is [JL1, Lemma 6.1. 
 Note that the notation in [JL1] is different than in this paper.   
In particular, $t_i$ here corresponds to $t_i^2$ in [JL1]. Thus 
divisibility by $4$ in [JL1, Lemma 6.1] corresponds to divisibility by 
$2$ in this paper.) For example, when ${\bf g}$
is ${\bf sl\ 2}$, then $T_{F}$ is just the set $$\{t^{-m}| m \in{\bf N}\}.$$   Note that this set
is a semigroup but is not a group.   This is true in general for 
$T_{F}$.

 Recall ${\cal F}$, the first filtration discussed in Section 4.
For each $\xi\in R(\pi)$, set $K_{\xi}^-$ equal to the subspace of $G^-$ such that 
 $F(U)\cap G^-\tau(\xi)=K_{\xi}^-\tau(\xi)$.  
 Similarly, set $K_{\xi}^+$ equal to the subspace of $U^+$ such that 
 $F(U)\cap U^+\tau(\xi)=K_{\xi}^+\tau(\xi)$.  It is shown in [JL2, 
 Section 4.9, 4.10] that 
 $$\gr_{\cal F}(F(U))=\oplus_{\xi\in R(\pi)}
 \gr_{\cal F}(K_{\xi}^-K_{\xi}^+\tau(\xi)).\leqno{(5.2)}$$  Note that 
  the 
inclusion of the  left hand side of (5.2) inside the right hand side
is just Theorem 4.7 applied to $F(U)$.  In particular, Theorem 4.7 gives a new 
proof of this inclusion.  Moreover, Theorem 4.7 can be thought of 
as
a  generalization of  this part of (5.2) 
to other left coideal subalgebras which admit an $\ad T$ module 
structure.

The
analysis in [JL2, Section 4, see Section 4.10],  shows that
$$\gr_{\cal F}(K_{\xi}^-K_{\xi}^+\tau(\xi))=(\ad U)\gr_{\cal F} 
\tau(\xi).$$ Moreover, 
$$(\ad U)\gr_{\cal F} \tau(\xi) \cong (\ad U)\tau(\xi)$$ as $\ad U$ 
modules  for each
$\tau(\xi)\in T_F$. 
 Thus one has the
direct sum decomposition in the nongraded case [JL2, Corollary 4.11]): $$ F(U)=
\oplus_{t\in T_F}(\ad U)  t.\leqno{(5.3)}$$
 Now (5.3) implies that 
$F(U)$ is the $\ad U$ module generated by the  algebra $F(U)\cap 
U^o$. The fact that $F(U)$ is a left coideal   was originally proved using
this fact and a weakened version of Lemma 1.2 ( [JL3, Lemma 5.3]). 
  
Note that $(\ad U)t$ is an {\it ad-invariant  left coideal} of $U$ for 
each $t\in T$.
 On the other hand, (4.11) guarantees that any left
coideal of $U$ contains an element of $T$.  Thus by (5.3) the minimal
ad-invariant left coideals of $U$ contained in $F(U)$ are exactly the vector
subspaces $(\ad U)t $ where $t\in T_F$.  This argument and result
is due to [HS, Theorem 3.9] where it is actually proved in the
dual setting. The description of the minimal ad-invariant left
coideals is, in turn, a crucial step in the classification of
bicovariant differential calculi on the quantized function algebra
$R_q[G]$.

The algebra $F(U)$ can be localized  by the normal elements $T_{F}$ to
obtain the larger coideal subalgebra ${\bf F}=F(U)T_{F}^{-1}$. Now
${\bf F}\cap T$ is just the subgroup  generated by
$T_{F}$. It is 
straightforward to show that ${\bf F}$ is generated by $x_i, y_it_i$,
and ${\bf F}\cap T$ for $1\leq i\leq n$.  In particular,  given $i$, 
there exists some $t\in T_{F}$
such that $(\ad x_i)t$ is a nonzero multiple of $x_it$. Hence
 $x_it\in F(U)$ and $x_i\in {\bf F}$.  A similar argument shows that  $y_it_i\in {\bf F}$ for all
$1\leq i\leq n$.  Thus ${\bf F}$ contains  ${\bf F}\cap T$, 
$x_i$, and $y_it_i$, for $1\leq i\leq n$. In the notation of  Corollary 
4.10, we get that ${\Delta}'={\Delta}''={\Delta^+}$, and, moreover,  ${\bf F}\cap T$, 
$x_i$, and $y_it_i$, $1\leq i\leq n$, generate
${\bf F}$.  It further follows that $G^-$ and 
$U^+$ are subalgebras of ${\bf F}$ 
 and that  $G^-U^+t$ is a subset 
of ${\bf F}$ for each $t\in {\bf F}\cap T$. Recall the notation of 
Theorem 4.7.   Note that ${\bf F}^-_{\eta}=G^-\subset {\bf F}$ and ${\bf 
F}^+_{\eta}=U^+\subset {\bf F}$ for all $\eta\in Q(\pi)$. Hence Theorem 4.7 implies 
the following direct sum decomposition of ${\bf F}$: $${\bf
F}=\oplus_{t\in {\bf F}\cap T}G^-U^+t.$$  Since ${\bf
F}\cap T$ is a subgroup of finite index in $T$, we see that  ${\bf F}$, and
hence $F(U)$,  is ``large'' in $U$ (For a stronger version of this, 
see [JL1, Theorem 6.4]).

A particular type of quantum Harish-Chandra module, defined
differently (and earlier) than those of Section 2, was introduced
in [JL3] in order to classify the
 primitive ideals of $U$. These modules were originally specified as
a subcategory of the  $F(U)$ bimodules with a ``compatible" $\ad
U$ action (see [JL3, 5.4] or [Jo, 8.2.3 and 8.4.1]). In [JL3, 5.4],
an  
$F(U)$ bimodule $M$ has a compatible $\ad U$ module structure provided that
$$\sum ((\ad a)(b\cdot m\cdot c)=\sum (\ad a_{(1)})b\cdot (\ad 
a_{(2)})m\cdot (\ad a_{(3)})c\leqno{(5.5)}$$ and 
$$(\ad t)m\cdot t=t\cdot m\leqno{(5.6)}$$
for all $a\in U$, $b$ and $c$ in $ F(U)$, $m\in M$, and $t\in F(U)\cap T$.   A 
different  definition of compatible is given in [Jo, 8.2.3].   In 
particular, the $\ad U$ action must satisfy the following condition 
in [Jo, 8.2.3]:
$$\sum ((\ad a_{(1)})m)\cdot a_{(2)}=a\cdot m\leqno{(5.7)}$$ for all
$a\in F(U)$ and $m\in M$. Note that (5.6) follows from (5.7) by 
setting $a=t$.  

The purpose of introducing the compatibility conditions (5.5) and (5.6) was to study the 
 specific
 Harish-Chandra module 
category ${\cal H}_{\chi}$ associated to a dominant regular 
weight $\Lambda$ defined in    [JL3, Section 5.7].  By 
[JL3, 5.12] and [Jo, 8.4.11], 
this category  is the same as the one described in 
[Jo, 8.4.1] using condition (5.7).   
Hence this category consists of   modules
with an $F(U)$ bimodule structure and $(\ad U)$ 
module action which satisfy both (5.5) and (5.7). In this paper, we 
say that $F(U)$ has a compatible $\ad U$ module action if both (5.5) 
and (5.7) hold. We show here that $F(U)$ bimodules with a 
compatible $\ad U$ module action fit exactly into the framework of Section 2.

Let $U^{op}$ denote the Hopf algebra with underlying vector space $U$, the opposite
multiplication, the same comultiplication and counit as $U$, and with
antipode $\sigma^{-1}$ ([Jo, 1.1.12]).  Note that the
 algebra $U\otimes U^{op}$ can be made into a Hopf algebra with
 comultiplication ${\it\Delta}(a\otimes b)=(Id\otimes {\bf tw}\otimes
 Id)({\it\Delta}\otimes{\it\Delta})(a\otimes b)$ where ${\bf tw}$ denotes the
twist map
 sending $a\otimes b$ to $b\otimes a$. The other Hopf
 operations can be defined similarly.  Observe that $U\otimes U^{op}$ is 
 isomorphic to $U_q({\bf g}\oplus {\bf g}^*)$ as a Hopf algebra.
 There is an algebra embedding
$\psi$ of  $U$ into $U\otimes U^{op}$ which sends an
element $u$ to $\sum u_{(1)}\otimes \sigma(u_{(2)})$.  
The image of $U$ in $U\otimes U^{op}$ under $\psi$ is not a Hopf 
subalgebra of $U\otimes U^{op}$.   However, by the next lemma it is a 
coideal subalgebra.

\begin{lemma}  The algebra $\psi(U)$ is a left coideal of $U\otimes 
U^{op}$.  
\end{lemma}

\noindent{\bf Proof:} By (1.3),
$${\it\Delta}(\sum u_{(1)}\otimes\sigma(u_{(2)}))=\sum (u_{(1)}\otimes\sigma(u_{(4)}))
\otimes (u_{(2)}\otimes\sigma(u_{(3)})).$$ Thus $\psi(U)$ is a
left coideal since ${\it\Delta}(u_{(2)})=\sum u_{(2)}\otimes 
u_{(3)}$.  $\Box$

\medskip

Let $F(U\otimes U^{op})$ denote the locally finite part of $U\otimes 
 U^{op}$.
We show that $F(U)$ modules with compatible $\ad U$ module action 
are $F(U\otimes U^{op})\psi(U)$ modules. 
The next lemma relates $F(U\otimes U^{op})$ 
to the locally finite part $F(U)$ of $U$.

\begin{lemma} $F(U\otimes U^{op})=F(U)\otimes F(U)^{op}$
\end{lemma}

\noindent
{\bf Proof:}
 Let ${\rm ad}^{op}$ denote the (left) adjoint action of $U^{op}$.  
 With sufficient care to indentification of elements in $U$ and 
 $U^{op}$, one checks using (1.3) that   
  $({\rm ad}^{op} \sigma(a))b= (\ad a) b$.  Thus $F(U^{op})=F(U)^{op}$ 
  as algebras.  $\Box$

\medskip
 Recall that since $\psi(U)$ is a left 
coideal and $F(U\otimes U^{op})$ is an $\ad U\otimes U^{op}$ module, 
we have that $F(U\otimes U^{op})\psi(U)=\psi(U)F(U\otimes U^{op})$.  
The next lemma shows that $F(U)$ is a free as a left $\psi(U)$ module.

\medskip
\begin{lemma}
The multiplication map induces an isomorphism $\phi$ of vector spaces
$$\psi(U)\otimes (1\otimes F(U)^{op})\rightarrow \psi(U)F(U\otimes 
U^{op}).$$
\end{lemma}

\noindent
{\bf Proof:}
Let 
$a\in F(U)$. Note that
$$\eqalign{a\otimes 1 &=\sum a_{(1)}\epsilon(a_{(2)})\otimes 1\cr
&=\sum a_{(1)}\otimes \epsilon(a_{(2)})\cr
&=\sum a_{(1)}\otimes 
a_{(3)}\sigma(a_{(2)})
\cr &=\sum \psi(a_{(1)})(1\otimes a_{(2)})\cr}
\leqno{(5.8)}$$ for all $a\in U$.  It follows that $$\psi(U)F(U\otimes 
U^{op})=
\psi(U)(1\otimes F(U)^{op}).$$ This proves that $\phi$ is surjective.

Suppose that 
$\sum_i \psi(c_i)(1\otimes b_i)=0 $ where the set $\{b_i\}$ is a 
linearly independent subset of $F(U)^{op}$.  We argue that each 
$\psi(c_i)=0$. This in turn implies that $\phi$ is injective.

 There is a version of Lemma 
1.3 for right coideal subalgebras.  In particular, $U=\oplus_tG^-U^+t$
and each $G^-U^+t$ is a right coideal of $U$. We may write
$c_i=\sum_tc_{it}$ where each $c_{it}\in G^-U^+t$.  It follows that 
$\psi(c_{it})(1\otimes b_i)\in G^-U^+t\otimes U$ for each $i$ and 
$t\in T$. Thus $\sum_{i}\psi(c_{it})(1\otimes b_i)=0$.   This allows us 
to reduce to the case where there exists $t\in T$ such that $c_i$ 
is in $G^-U^+t$ for all $i$. 

Recall the notation of Section 4, Filtration II. Let $(M,N)$ be the 
maximum  value of the set of bidegrees of the $c_i$.  Reordering if 
necessary, we may assume that $c_1$ has bidegree $(M,N)$. Choose 
$[\lambda,\mu]\in \max(c_1)$.   By say (4.11), we have $$(\pi_{\lambda,\mu}\otimes 
Id )(\sum_i\psi(c_i)(1\otimes b_i))=\sum_i\pi_{\lambda,\mu}(c_i)\otimes 
\sigma(t)b_i.$$  Note that $\sigma(t)=t^{-1}$.  Since the set $\{b_i\}$ is linearly independent, the 
set $\{\sigma(t)b_i\}$ is also linearly independent.  Hence 
$\pi_{\lambda,\mu}(c_i)=0$ for each $i$.   The choice of 
$[\lambda,\mu]$ now forces $c_i=0$ for each $i$.
$\Box$

\medskip

The next result shows that $F(U)$ modules with compatible $\ad U$ 
modules are just $F(U\otimes U^{op})\psi(U)$ modules.

\begin{theorem} The  set of $F(U)$ bimodules   with  compatible $\ad 
U$  module action can be identified with the set of     $F(U\otimes U^{op})\psi(U)$ modules.
\end{theorem}

\noindent
{\bf Proof:} 
Let $M$ be a $F(U\otimes U^{op})\psi(U)$ module.   Note that $M$ is a 
$F(U)$ bimodule in a natural way.   In particular, set $a\cdot m\cdot 
b=(a\otimes b)m$ for all $a,b\in F(U)$ and $m\in M$.   We define an 
action of $\ad U$ on $M$ by setting $(\ad c)m=\psi(c)m$ for all $c\in 
U$.  By (5.8), it follows that this $\ad U$ action satisfies (5.7).
A straightforward computation shows that this action satisfies (5.5) 
as well.

Now  let $M$ be an
$F(U)$ bimodule with compatible $\ad U$ module action. Make $M$  into a
$F(U\otimes U^{op})\psi(U)$ module by setting $$(a\otimes b)m 
=a\cdot 
m\cdot b{\rm \ and \ }(\psi(c))m=(\ad c)m\leqno{(5.9)}$$ for all $a\otimes b\in F(U\otimes 
U^{op})$,
$c\in U$, and $m\in M$. 

One checks that  $\psi(c)(1\otimes b)=\sum(1\otimes (\ad c_{(2)})b)\psi(c_{(1)})$.  By (5.5), 
$(\ad c)(m\cdot b)=(\ad c_{(1)}\cdot m)\cdot ((\ad c_{(2)}) b)$.  
Hence the action of $\psi(c)$ on $(1\otimes b)m$ described in (5.9) agrees with the
action of $(1\otimes (\ad c_{(2)})b)$ on $\psi(c_{(1)})m$.   
Therefore, to 
show that the action in (5.9) is well defined it is sufficient to show 
that the action of an element $x\in  F(U\otimes U^{op})\psi(U)$ on $M$ 
is independent of the way $x$ is written as a sum of terms of the 
form $bu$ where $b\in  F(U\otimes U^{op})$ and $u\in \psi(U)$.
 
 The compatibility condition (5.7) ensures that  $$(\ad c)(a\cdot m\cdot 
b)=\sum(\ad ca_{(1)})(m\cdot ba_{(2)}) $$ for all 
$a\otimes b\in F(U\otimes 
U^{op})$,
$c\in U$, and $m\in M$. 
   Thus using (5.9) formally, we see that 
$$\eqalign{\psi(c) ((a\otimes b)m) &=\psi(c) ((a\cdot m\cdot  b)) 
\cr &=\sum\psi(ca_{(1)})(m\cdot 
ba_{(2)})\cr &=\sum\psi(ca_{(1)})((1\otimes 
ba_{(2)})m).\cr}$$  In particular the action of 
$\psi(c)(a\otimes b)$  agrees with the action of $\sum\psi(ca_{(1)})(1\otimes 
ba_{(2)})$ on $M$.  By  Lemma 5.4,  every element in $F(U\otimes U^{op})\psi(U)$ can be expressed 
uniquely in the form $\sum_i\psi(a_i)(1\otimes b_i)$ where $\{b_i\}$ is 
a basis of $F(U)^{op}$.   The theorem now follows.
$\Box$

\medskip

One can  apply the results of Section 2 to the study of 
Harish-Chandra modules for the pair $U\otimes
U^{op}$, $\psi(U)$.  Identify the algebra $U\otimes 
U^{op}$ with $U_q({\bf g}\oplus {\bf g}^*)$. 
Let $\tilde\kappa$ denote the 
conjugate linear Chevalley antiautomorphism of Section 2 associated 
here to the quantized enveloping algebra $U_q({\bf g}\oplus {\bf g}^*)$. One can find a Hopf algebra 
automorphism $\Upsilon\in {\cal H}_{\bf R}$ such that $\Upsilon(\psi(U))$ is invariant
under $\tilde\kappa$.  Thus the results in
Section 2 apply here.  However, the main results of Section 2, such as Theorem 
2.7, can be proved  easily in this case since $\psi(U)$ is 
isomorphic as an algebra to $U_q({\bf g})$.  Thus it acts completely 
reducibly on all finite-dimensional $\psi(U)$ modules.  Furthermore, 
one checks that all 
 finite-dimensional 
$\psi(U)$ modules are unitary  using the fact that this is 
true for  
$U_q({\bf g})$.

For an example of a Harish-Chandra module associated to 
the pair $U\otimes
U^{op}$, $\psi(U)$, consider two left  
$U$ modules $M$ and $N$. Define the $U$ bimodule ${\rm Hom}(M,N)$  by $ (a\cdot f\cdot b)(m)=af(bm).$ As
explained in [JL3, 5.4] and [Jo, 8.2.3], ${\rm Hom}(M,N)$ has a compatible $(\ad U)$
module structure in the sense of (5.5) and (5.7) given by $(\ad
a)f=\sum a_{(1)}\cdot f\cdot \sigma( a_{(2)})$.  Thus from the above
Theorem 5.5, we see that ${\rm Hom}(M,N)$ is a $(F(U )\otimes
F(U)^{op})\psi(U)$ module.  By Theorem 2.7, the sum of all 
finite-dimensional $\ad U$ modules $F(M,N)$ inside of ${\rm Hom}(M,N)$
is a Harish-Chandra module for the pair $U\otimes U^{op}$,
$\psi(U)$. 

In [JL3, Theorem 5.13] (see also [Jo, Chapter 8]),  the theory of Harish-Chandra
modules associated to the pair $U\otimes U^{op}$, $\psi(U)$ is
used to prove an equivalence of categories   between certain
Harish-Chandra modules
  and various
 category ${\cal O}$ modules. This is  critical in
obtaining the quantum version of Duflo's theorem: every primitive
ideal of $U$ is the annihilator of a highest weight simple module
([JL3, Corollary 6.4] or [Jo, 8.4.17]).

\section{Nilpotent and parabolic coideal subalgebras}

Yet another   left coideal subalgebra of $U$ is $G^{-}$,
an obvious quantum analog of $U({\bf n}^-)$.  In
this section, we consider coideal subalgebras of $G^-$ which
correspond to classical enveloping algebras of  Lie subalgebras of
${\bf n}^-$ and related Lie subalgebras of ${\bf g}$.   Most of
the results presented here are from [Ke].

Let $\pi'$ be a subset of the simple roots $\pi$ of ${\bf g}$.
There are a number of Lie subalgebras of ${\bf g}$ which can be
associated to $\pi'$.   The most obvious is the semisimple Lie
subalgebra ${\bf m}$ of ${\bf g}$ generated by the $e_i,f_i$, $h_i$,
for those $i$ with $\alpha_i\in \pi'$.
Since the simple roots $\pi'$ associated to the root system of 
${\bf m}$ are contained in the simple roots $\pi$ of ${\bf g}$, the
entire picture can be lifted to the quantum setting. In
particular, $U_q({\bf g})$ contains a Hopf subalgebra ${\cal M}$
isomorphic to $U_q({\bf m})$ and generated by the
$x_i,y_i,t_i,t_i^{-1}$ for the same $i$.
 Set ${\cal M}^-={\cal M}\cap G^-$ and ${\cal M}^+={\cal M}\cap U^+$.

Let $\Delta'$ denote the set of positive roots associated to the simple
roots $\pi'$.  
The vector space ${\bf n}^-_{\pi'}$ spanned by
the root vectors $f_{-\gamma}$, $\gamma$
 in $\Delta^+-\Delta'$, is a second Lie subalgebra of ${\bf n}^-$.
Let ${\bf m}^-$ denote the Lie subalgebra of ${\bf m}$ generated by 
the $f_i$ for $\alpha_i\in \pi'$.  Then $${\bf n}^-={\bf
n}^-_{\pi'}\oplus {\bf m}^-.$$  Thus the multiplication map
defines a vector space isomorphism: $$U({\bf n}^-)\cong U({\bf
n}^-_{\pi'})\otimes U({\bf m}^-).\leqno{(6.1)}$$  We shall see that the algebra
$U({\bf n}^-_{\pi'})$ can be lifted to the quantum setting using a
coideal subalgebra.

Let $G_{\pi-\pi'}^-$ be the subalgebra of $G^-$ generated by the
$y_it_i$ such that $\alpha_i$ is in $\pi-\pi'$. Note that
$G_{\pi-\pi'}^{-}$ is a left coideal subalgebra of $G^-$. Now
$G_{\pi-\pi'}^{-}$ is generated by weight vectors and in
particular, $(\ad T)G_{\pi-\pi'}^{-}=G_{{\pi-\pi'}}^{-}$. Also,
$(\ad x_i)y_jt_j=0$ for all $i\neq j$. Thus $(\ad ({\cal
M}^+T))G_{\pi-\pi'}^{-}\subset G_{\pi-\pi'}^{-}$.
 Recall that ${\cal M}$ is equal to the    quantized
enveloping $U_q({\bf m})$.  Hence the triangular decomposition (1.12)
  implies that ${\cal M}T={\cal M}^-{\cal M}^+T$.
Hence $(\ad
 {\cal M}^-)G_{\pi-\pi'}^{-}$ equals $(\ad {\cal M})G_{\pi-\pi'}^{-}$.
By Lemma 1.2, $(\ad
 {\cal M}^-)G_{\pi-\pi'}^{-}$ is a left coideal. Let $N^-_{\pi'}$ be the
subalgebra of $G^-$ generated by  $(\ad {\cal
M}^-)G_{\pi-\pi'}^{-}$. It is a left coideal subalgebra since
$(\ad
 {\cal M}^-)G_{\pi-\pi'}^{-}$ is a left
coideal.

By [Ke], one has  a quantum analog of (6.1).  Namely there is an
isomorphism of vector spaces $$G^-\cong N^-_{\pi'}\otimes
 {\cal M}^-.\leqno{(6.2)}$$   K\'eb\'e
  actually proves a stronger result 
 with this as a consequence, namely,  $G^-$ is
isomorphic to the smash product of $N^-_{\pi'}$ and ${\cal M}^-$.

By construction, $N^-_{\pi'}$ is generated by weight vectors and hence 
is a direct sum of its weight spaces.   By Lemma 4.3, we can find a 
subset $\Delta_1$ of $\Delta^+$ such that $N^-_{\pi'}$ is generated by 
weight vectors 
 $\tilde f_{-\gamma}$ of weight $\gamma\in \Delta_1$  
 which specialize to root vectors in 
$U({\bf n}^+)$.  By [L3, proof of Proposition 2.2], $\Delta_1$ consists 
of those
 positive roots which are not linear combinations of roots 
in $\pi'$.  In particular, $N^-_{\pi'}$ specializes to  
 $U({\bf n}^-_{\pi'})$ as $q$ goes to $1$ ([L3, proof of Proposition 2.2]).
Thus the left coideal subalgebra  $N^-_{\pi'}$ is a natural choice
of quantum analog of $U({\bf n}^-_{\pi'})$ inside of
 $U({\bf g})$.

It is instructive to look at the generators of $N^-_{\pi'}$.  Let
$I$ be a tuple $(i_1,\dots,i_r)$ of (arbitrary) length $r$ and suppose 
that
     $\alpha_{i_s}$ is in $\pi-\pi'$ for $1\leq s\leq r$. By the argument in
     [L3, Proposition 2.2], the algebra $N^-_{\pi'}$ is
generated by elements of the form $$Y_{I,j}=(\ad y_{i_1}\cdots
y_{i_r})y_jt_j $$ where   $\alpha_{j}\notin \pi'$. 

Now each $Y_{I,j}$ is an element of  the subcoideal  $(\ad {\cal M}^-) y_jt_j$ of  $N^-_{\pi'}$ as well as an element of $G^-$.  Hence $$(Id\otimes
\pi_{0,0}){\it\Delta}(Y_{I,j}) =Y_{I,j}\otimes 1.$$   Thus 
$${\it\Delta}(Y_{I,j})=Y_{I,j}\otimes 1 + \sum Y_i\otimes Y'_i $$
where $Y_i$ is in $U$ and $Y'_i$ is in  $(\ad{\cal
M}^-)y_jt_j.$ We can actually say more about
the  $Y_i$. First recall that  $Y_{I,j}$ is in  $G^-$.
Set $\lambda=\alpha_{i_1}+\cdots +\alpha_{i_r}+\alpha_j$ and note that the 
weight of $Y_{I,j}$ is $-\lambda$.   Set $\mu=0$. We may apply (4.7) to 
$Y_{I,j}\tau(-\lambda)$ using this 
$\lambda$ and $\mu.$   By  (4.7) and weight space considerations, each $Y_i$ is 
in ${\cal M}\cap G^-U^o$.   Furthermore, (4.7) implies that  each $Y_i \in
U^-\tau(\lambda)$. Since $\tau(\lambda)\in{\cal M}t_j$, it follows 
that  each $Y_i$ is an element of $({\cal M}\cap 
U^-U^o)t_j$.  In particular, we  get
that   (see [AJS, 
Proposition C.5]) 
$${\it\Delta}(Y_{I,j})\in Y_{I,j}\otimes 1 + ({\cal M}\cap 
U^-U^o)t_j\otimes (\ad
{\cal M}^-)(y_jt_j).\leqno{(6.3)} $$ The elements $Y_{I,j}$ also satisfy a
uniqueness property. In particular, by [L2, Proposition 4.1], if
$Y$ is an element of $G^-$ of weight $-\lambda$ such that 
$${\it\Delta}(Y)\in Y\otimes 1 +
({\cal M}\cap 
U^-U^o)t_j\otimes (\ad {\cal M}^-)(y_jt_j)$$ then $Y$ is a
nonzero scalar multiple of $Y_{I,j}$.  This uniqueness property will be 
used in the uniqueness result Theorem 7.5 concerning quantum symmetric pairs.

Let ${\bf
n}^+_{\pi'}$ be the Lie subalgebra of ${\bf n}^+$ spanned by the
root vectors $e_{\gamma}$, where $\gamma$ runs over 
 $\Delta^+-\Delta'$. One can similarly define left coideal subalgebras $N^+_{\pi'}$ of
$U^+$ which are analogs of $U({\bf n}^+_{\pi'})$.
These can be constructed directly using the
same methods described above for $N^-_{\pi'}$. 

Of course, one could take the perspective of right coideal
subalgebras instead of left coideal subalgebras.  This will be
useful in the next section.  For example, right coideal analogs of
$U({\bf n} ^+_{\pi'})$ are subalgebras of $G^+$  defined using the
right adjoint action,
$$({\rm ad}_r a)b=\sum \sigma(a_{(1)})ba_{(2)}\leqno{(6.4)}$$ for all
$a$ and $b$ in $U$. Let $G^+_{\pi-\pi'}$ be the subalgebra of
$G^+$ generated by the $x_it_i^{-1}$ for all $i$ such that $
\alpha_i\in \pi-\pi'$. Then the subalgebra $N^{+}_{\pi',r}$
generated by $({\rm ad}_r {\cal M}^+)G^+_{\pi-\pi'}$ is a right
coideal subalgebra of $G^+$ and an analog of $U({\bf n} ^+_{\pi'}).$  The algebra
$N^{+}_{\pi',r}$ is generated by elements of the form
$$X_{I,j}=({\rm ad}_r x_{i_1}\cdots x_{i_r})x_jt^{-1}_j $$ where
each $\alpha_{i_s}\in {\pi'}$ and 
$\alpha_{j}\notin \pi'$.   Moreover the comultiplication of these
elements is similar to that of the $Y_{I,j}$, e.g.,
 $${\it\Delta}(X_{I,j})\in 1\otimes X_{I,j} + ({\rm ad}_r {\cal
M}^+)(x_jt^{-1}_j)\otimes({\cal M}\cap G^+U^o)t_j^{-1}.\leqno{(6.5)} $$

Using $N^-_{\pi'}$, $N^+_{\pi'}$,  and ${\cal M}^{-}$, one can
construct what are called generalized Verma modules.  Let ${\cal P}$ be the
subalgebra of $U$ generated by ${\cal M}$, $U^o$, and
$N^+_{\pi'}$. Note that ${\cal P}$ is a left coideal subalgebra
since it is generated by left coideal subalgebras. It is an analog
of the enveloping algebra of the parabolic Lie subalgebra $({\bf
m}+{\bf h})\oplus{\bf n}^+_{\pi'}$. Using (6.2), one obtains
 an isomorphism of vector spaces via the multiplication
map $$U\cong N^-_{\pi'}\otimes {\cal P} .\leqno{(6.6)}$$

Let  $W$ be a finite-dimensional simple ${\cal M}$ module. Extend
the action of ${\cal M}$ on $W$ to $U^o$ by insisting that the
highest weight generating vector of $W$ is a weight vector of say
weight $\Lambda$ with respect to the action of $T$. Extend further
the action on $W$ to $N^+_{\pi'}$ by insisting that the
augmentation ideal of $N^+_{\pi'}$ acts as zero on all vectors in
$W$. These extensions make $W$ into a ${\cal P}$ module. The
generalized Verma module $M_{\pi'}(\Lambda)$ is defined to be
$U\otimes_{\cal P} W$.  In particular,  elements of $U$
 act by left multiplication and $p u\otimes w=\sum (\ad p_{(1)})u\otimes p_{(2)}w$
 for all $p\in {\cal P}$, $u\in U$, and $w\in W$. As a left
$N^-_{\pi'}$ module,
$U\otimes_{\cal P}W\cong N^-_{\pi'}\otimes W$.  Furthermore,
the action of ${\cal M}$ on ${\bf N}^-_{\pi'}$
is  both locally finite and semisimple. Hence the generalized Verma module
$M_{\pi'}(\Lambda)$ is a Harish-Chandra module for the pair $U$,
${\cal M}$.

Using the coideal subalgebras discussed in this section, one can
form  quantized homogenous spaces as in Section 3.  For example,
the    homogeneous space  associated to $G^-$,
$R_q[G/N]=R_q[G]^{G^-}_l$ is studied in [Jo, Chapter 9]    where
it is used to obtain the complete description of the prime and primitive
spectra of the quantized function algebra $R_q[G]$.

  \section{Quantum symmetric pairs}
  We turn now to the theory of
  quantum symmetric pairs.  First, we present the construction and
characterization of the
  coideal subalgebras used to form such pairs.   The results are drawn from
[L2] and [L3],
  but the methods in this paper are often simpler. 
The involutions
 used to construct these algebras are  given in a   concrete
 fashion here.  The relations for the coideal subalgebras as algebras 
 are also presented more explicitly.  
  Moreover, using the results of Section 4,
  we give a new, less intricate, proof of
  the uniqueness characterization for the subalgebras used to
  form quantum symmetric pairs (see Theorem 7.5 below.)  The
  Harish-Chandra module and symmetric space theory associated to
  these pairs is also described with the aid of Sections 2 and 3.
  
  A  symmetric pair is defined for each Lie 
algebra involution (equivalently, a Lie algebra automorphism of 
order $2$) of ${\bf g}$. More precisely, let $\theta$ be a
 Lie algebra involution of $\bf g$.
Write ${\bf g}^{\theta}$ for the Lie subalgebra of ${\bf g}$ 
consisting of elements fixed by $\theta$.   The pair ${\bf g},{\bf 
g}^{\theta}$ is a classical symmetric pair. A classification of 
involutions and
classical symmetric pairs up to isomorphism can be found in 
[He1, Chapter 10, Sections 2, 5, and 6] and [OV, Section 4.1.4]. 

Let ${\bf p}=\{v\in {\bf g}|\theta(v)=-v\}$. A commutative Lie subalgebra of ${\bf g}$ 
which is reductive 
in ${\bf g}$ and is equal to its centralizer in ${\bf p}$
is called a Cartan subspace of ${\bf p}$ (see [D, 1.13.5].)
A Cartan subalgebra ${\bf h'}$ of ${\bf g}$ is called 
maximally split ([V, Section 0.4.1]) with respect to $\theta$ provided that ${\bf h'}\cap 
{\bf p}$ is a Cartan subspace of ${\bf p}$.
By [D, 1.13.6, 1.13.7], ${\bf p}$ contains Cartan subspaces and moreover
each Cartan subspace can be extended to a Cartan subalgebra of ${\bf 
g}$. 

Recall that we have already  specified  a Cartan subalgebra ${\bf h}$
of ${\bf g}$. Let $\theta$ be
an involution of ${\bf g}$ such that ${\bf h}$ is maximally split with 
respect to $\theta$.
Let ${\cal L}$ be the set of Lie algebra automorphisms $\psi$ of ${\bf g}$ such that
  $\psi({\bf p}\cap {\bf h})$ is a subset of ${\bf h}$.
If $\psi\in {\cal L}$ then 
 ${\bf h}$ is also maximally split with respect to the involution   
$\psi\theta\psi^{-1}.$   By [D, 1.13.7 and 1.13.8], 
one can replace $\theta$ by  $\psi\theta\psi^{-1}$ for some  
$\psi\in {\cal L}$    so that 
$\theta$ also satisfies the following conditions:
\begin{enumerate}
\item[(7.1)] $\theta({\bf h})={\bf h};$ 
\item[(7.2)] if $\theta(h_i)=h_i$ then  $\theta(e_i)=  e_i$ and 
$\theta(f_i)=  f_i$;
\item[(7.3)]if $\theta(h_i)\neq h_i$ then $\theta(e_i)$ (resp. $\theta(f_i)$) is a nonzero  root vector in
${\bf n}^-$ (resp. ${\bf n}^+$).
\end{enumerate}
By [D, 1.13.8], $\theta$ also induces an 
automorphism $\Theta$ of the root system $\Delta$.

Now consider an arbitrary involution $\theta'$ of ${\bf g}$.  One can 
find a  Lie algebra automorphism $\Upsilon$ of ${\bf g}$ so that ${\bf 
h}$ is maximally split with respect to the involution 
$\Upsilon\theta'\Upsilon^{-1}$.   In the quantum case,
we do not have as much flexibility 
in ``moving'' involutions around using an automorphism of $U$. 
In particular, there is only one choice of quantum Cartan subalgebra,  
since the only invertible elements of $U$ are the nonzero scalars and 
the elements of $T$.  Hence any automorphism of $U$ restricts to an 
automorphism of $T$. Thus  the relationship between an 
 involution of ${\bf g}$ and the particular Cartan 
subalgebra ${\bf h}$ is important in lifting the involution to the 
quantum case. In this section, we call an involution $\theta$ of 
${\bf g}$ a maximally split involution if ${\bf h}$ is maximally 
split with respect to $\theta$ and $\theta$ satisfies (7.1), (7.2),
and (7.3).  (Similar terminology was introduced in [Di, Section 5].)  
We discuss lifts of 
maximally split involutions and the associated quantum symmetric 
pairs.  
   There are also   
a few scattered results on quantum symmetric pairs
when the involution is not
 maximally split.   The reader is referred to [G] and [BF] for more 
 information.

For the remainder of this section, let
 $\theta$ be a maximally split involution with respect to the 
 fixed Cartan subalgebra ${\bf h}$. 
 Consider the Cartan subspace ${\bf a}={\bf p}\cap {\bf h}$ of ${\bf 
 p}$. Since ${\bf a}$ is subset of ${\bf h}$, the action of $\ad {\bf a}$ 
 on ${\bf g}$ is semisimple.   Given $\lambda\in {\bf a}^*$, 
 set $${\bf g}_{\lambda}=\{x\in {\bf g}|(\ad a)x=\lambda(a)x{\rm \ 
 for\  all \ }a\in {\bf a}\}.$$  Let $$\Sigma=\{\lambda\in 
 {\bf a}^*|\ {\bf g}_{\lambda}\neq 0\}.$$ We can write
 ${\bf g}=\oplus_{\lambda\in \Sigma}{\bf g}_{\lambda}$.  Furthermore, 
 by [OV, Theorem 3.4.2], $\Sigma$ is an abstract root system 
called the restricted root system associated to
$\theta$ (or  more precisely, to ${\bf g},{\bf a}$.)  A classification 
of restricted root systems associated to involutions can be found in 
[Kn, Chapter VI, Section 11] (see also [He1, Chapter X, Section F 
under Exercises and Further Results]).  Note that an abstract root system is slightly 
more general than an   ordinary 
root system  (often called a reduced root system) described in 
[H, Chapter III]. Good references for abstract root systems 
are   [Kn, Chapter II, Section 5] and  [OV, Chapter 3, 
Section 1.1]. The abstract root systems have been 
classified as the set  of  
 reduced root systems  and one additional nonreduced family referred to as 
type $BC$ ([Kn, Chapter II, Section 8]).  

Before discussing the quantum case, we further describe the action of 
$\theta$ on the generators of ${\bf g}$.
 Set
$\Delta_{\Theta}=\{\alpha\in \Delta|\Theta(\alpha)=\alpha\}$ where 
$\Theta$ is the associated root system automorphism. This
is the root system for the semisimple Lie subalgebra ${\bf m}$ of
${\bf g}$ generated by the $e_i,f_i,h_i$ with $\theta(h_i)=h_i$.
Write ${\bf m}= {\bf m}^-\oplus  {\bf m}^o\oplus  {\bf m}^+$ for
the obvious triangular decomposition of ${\bf m}$. Set 
$\pi_{\Theta}=\Delta_{\Theta}\cap \pi.$ 
  Note that
$\pi_{\Theta}$ is  a set   of positive simple roots for the root system
$\Delta_{\Theta}$. Write $Q(\pi_{\Theta})$ for the 
lattice of integral linear combinations of the simple roots in
$\pi_{\Theta}$.  Let $Q^+(\pi_{\Theta})$ be the set   of nonnegative integral
linear combinations of the elements in $\pi_{\Theta}$.

Note that $\pi_{\Theta} =\Theta(\pi)\cap \pi$. 
Also, $\Theta(-\alpha_i)\in \Delta^+$ for all 
$\alpha_i\notin\pi_{\Theta}$ by (7.3).   It follows that $$\Theta(-\alpha_i)\in 
\sum_{\alpha_j\notin\pi_{\Theta}}{\bf 
N}\alpha_j+Q^+(\pi_{\Theta})\leqno{(7.4)}$$
for each $\alpha_i\notin\pi_{\Theta}$. Since $\Theta$ is a root system automorphism, every element of $\Delta$
can be written as an integral linear 
combination of roots in $\{\Theta(\alpha_i)|\alpha_i\in \pi\}$
where either all the coefficients are positive or all the coefficients are 
negative.  
 Hence  each $\alpha_i\notin\pi_{\Theta}$
 can be written as a linear 
combination of elements in 
$\{\Theta(\alpha_i)|\alpha_i\notin\pi_{\Theta}\}\cup \pi_{\Theta}$
with just negative integers as   coefficients.
Observation (7.4) thus implies that there exists
 a permutation $p$ on the set $\{i\
|\alpha_i\in\pi-\pi_{\Theta}\}$
 such that
for each $\alpha_i\in\pi- \pi_{\Theta}$, $$
\Theta(-\alpha_i)-\alpha_{p(i)}\in
Q^+(\pi_{\Theta}).\leqno{(7.5)}$$ 

Choose a maximal subset $\pi^{*}$ of $\pi-\pi_{\Theta}$ such that if 
$j= p(j)$ then $\alpha_j\in \pi^{*}$ and if $j\neq p(j)$, then exactly 
one of the pair $\alpha_j,\alpha_{p(j)}$ is in $\pi^{*}$.
Consider $i$ such that $\alpha_i\in\pi^{*}$.  The root vector
$e_{p(i)}$ associated to the   simple root $\alpha_{p(i)}$ satisfies 
$(\ad f_j)e_{p(i)}=[f_j,e_{p(i)}]=0$ for all 
$\alpha_j\in\pi_{\Theta}$.  Thus  $e_{p(i)}$ is a lowest weight
vector for the action of $\ad {\bf m}^-$. Let $V$ be the
corresponding simple $\ad {\bf m}$ module generated by $e_{p(i)}$.
By (7.3), $\theta(f_i)$ is a   root vector in 
${\bf n}^+$.  Furthermore (7.5) implies that
the weight of this root vector is 
$\alpha_{p(i)}$ plus some element in $Q^+(\pi_{\Theta})$.  Thus $\theta(f_i)$ can be written as a
 bracket $[a_1[a_2,\dots,[a_{s-1},a_s]\dots]$ where exactly one of the 
$a_j$ equals  $e_{p(i)}$ and the others are elements of  
${\bf m}^+$.  Using the Jacobi identity, we see that 
$\theta(f_i)$ is an element of $(\ad {\bf m}^+)e_{p(i)}$.   In 
particular $\theta(f_i)$ is   an element of $V$. Furthermore,
since elements of ${\bf m}^+$ commute with $f_i$ and thus with
$\theta(f_i)$, we see that $\theta(f_i)$ must be a highest weight vector of
$V$. 
 Thus we can find a sequence of elements
$\alpha_{i_1},\dots, \alpha_{i_r}$ in $\pi_{\Theta}$ and
a sequence of positive integers $m_{ 1},\dots, m_{ r}$ such that (up to a slight 
adjustment of $\theta$)
$$\theta(f_i)=(\ad e_{i_1}^{(m_1)}\cdots
e_{i_r}^{(m_r)})e_{p(i)}.\leqno{(7.6)} $$ Here  $e_j^{(m)}=e_j^m/m!$.
 We may further assume that both the
sequence  of roots and the sequence of  integers are chosen
so that each $(\ad e_{i_s}^{(m_s)}\cdots
e_{i_r}^{(m_r)})e_{p(i)}$, $1\leq s\leq r$, is an extreme vector
of $V$.  (In particular, $(\ad e_{i_{s}}^{(m_s)}\cdots
e_{i_r}^{(m_r)})e_{p(i)}$ is a highest weight vector for the action 
of $\ad e_{i_s}$ and  $(\ad e_{i_{s-1}}^{(m_{s-1})}\cdots
e_{i_r}^{(m_r)})e_{p(i)}$ is a lowest weight vector for the action of 
 $\ad f_{i_s}$.) Suppose that the sequence $\alpha_{j_1},\dots, 
 \alpha_{j_s}$ of elements in $\pi_{\Theta}$ and
the positive integers $n_{ 1},\dots, n_{s}$ also satisfy this 
condition on extreme vectors and that 
$\sum_km_k\alpha_{i_k}=\sum_kn_k\alpha_{j_k}$.  By [Ve], $(\ad e_{i_1}^{(m_1)}\cdots
e_{i_r}^{(m_r)})e_{p(i)}=(\ad e_{j_1}^{(n_1)}\cdots
e_{j_s}^{(n_s)})e_{p(i)}$.  Thus (7.6) is independent of the choice 
of such sequences. 

Using lowest weight vectors instead of 
highest weight vectors, we  obtain 
$$\theta(e_{p(i)})=(\ad f_{i_r}^{(m_r)}\cdots
f_{i_1}^{(m_1)})f_{i}\leqno{(7.7)} $$ 
up to a nonzero scalar.
 A straightforward ${\bf sl}\ 2$ computation shows that $$(\ad
e_{i_1}^{(m_1)}\cdots e_{i_r}^{(m_r)})[ (\ad f_{i_r}^{(m_r)}\cdots
f_{i_1}^{(m_1)})f_{i}]=f_i$$ and $$(\ad f_{i_r}^{(m_r)}\cdots
f_{i_1}^{(m_1)})[( \ad e_{i_1}^{(m_1)}\cdots
e_{i_r}^{(m_r)})e_{p(i)}]=e_{p(i)}.\leqno{(7.8)}$$
Since $\theta^2$ is the identity, the scalar in (7.7) must be $1$.

Set $m(i)= m_1+\dots +m_r $.
Now  $[\theta(e_i),\theta(f_i)]=\theta(h_i)$ is an element of ${\bf 
h}$ by (7.1).  Furthermore,   by (7.2) and
(7.3), $\theta(h_i)$ must be the coroot  $h_{\Theta(\alpha_i)}$ associated to the root
$\Theta(\alpha_i)$. 
The description of the Chevalley basis for ${\bf g}$ given in [H,
Proposition 25.2 and Theorem 25.2] ensures that  both $\theta(e_i)$ and $\theta(f_i)$
are Chevalley basis vectors up to a sign. Furthermore,   
by [H, Proposition 25.2(b)] 
and (7.6),we must have $$\theta(e_{i})=(-1)^{m(i)}(\ad f_{i_1}^{(m_1)}\cdots
f_{i_r}^{(m_r)})f_{p(i)}.$$  Similarly, by [H, Proposition 25.2(b)] 
and (7.7) $$\theta(f_{p(i)})=(-1)^{m(i)}(\ad e_{i_r}^{(m_r)}\cdots
e_{i_1}^{(m_1)})e_{i}.$$  Note that when $p(i)=i$, we have  $$(\ad e_{i_r}^{(m_r)}\cdots
e_{i_1}^{(m_1)})e_{i}=(\ad e_{i_1}^{(m_1)}\cdots
e_{i_r}^{(m_r)})e_{i}.$$  Hence $m(i)$ is even in this case.

The above analysis allows us to better describe the root space automorphism
$\Theta$.  Let $W'$ denote the Weyl group associated to the root
system $\Delta_{\Theta}$ of ${\bf m}$ considered as a subgroup of
the Weyl group  of $\Delta$.  Let $w_o$ denote the longest element
of $W'$. Note that $w_o$ is a product of reflections in $W'$ but
can also be considered as an element of $W$. Let $d$ be the  diagram
automorphism on $\pi_{\Theta}$ such that $d=-w_o$ when restricted
to $\pi_{\Theta}$. Note that $d$ induces a permutation on  the set 
$\{i|\alpha_i\in\pi_{\Theta}\}$ which we also denote by $d$.   In particular, 
given $\alpha_i\in\pi_{\Theta},$ $d(\alpha_i)=\alpha_{d(i)}$.
 Extend $d$ to a function on $\pi$, and thus to
$\Delta$, by setting $d(\alpha_i)=\alpha_{p(i)}$ for
$\alpha_i\notin\pi_{\Theta}$. It follows that $\Theta=-w_od$. Note
that this forces $d$ to be a diagram automorphism of the larger
root system $\Delta$.  

Before lifting   $\theta$ to the quantum case, we  recall and introduce 
more notation.
The  right
adjoint action is defined by (6.4). This action on the generators 
  of  $U$ is given by: $$( {\rm ad}_r
y_i)b =   by_i - y_it_ibt_i^{-1} \quad ({\rm ad}_r
x_i)b=t_i^{-1}bx_i-t_i^{-1}x_ib \quad ( {\rm ad}_r
t_i)b=t_i^{-1}bt_i$$ for $1\leq i\leq n$. 
Recall the definitions of $[m]_q$ and $ q_i$ used to define the quantized enveloping algebra 
((1.4)-(1.7)).  The divided powers of $x_i$ 
 and $y_i$ are defined by 
 $x_i^{(m )}=x_i^{m}/[m]_{q_i}!$ and $y_i^{(m 
 )}=y_i^{m}/[m]_{q_i}!$.  (Note that these are quantum analogs of the 
 divided power $e_i^{(m)}$.)
Let ${\cal M}$ denote the  subalgebra of $U$ generated by the
corresponding elements $x_i,y_i,t_i,t_i^{-1}$ where
$\theta(h_i)=h_i$.   Note that ${\cal M}$ is just a copy of the
quantized enveloping algebra $U_q({\bf m})$ so this notation is 
consistent with that of Section 6. 
Let $\iota$ be the ${\bf C}$
algebra automorphism of $U$ fixing $x_it_i^{-1}$ and  $t_iy_i$ for $ 1\leq i\leq n$,
sending $t$ to $t^{-1}$ for all $t\in T$ and $q$ to
$q^{- 1}$.
Recall the sequences $\alpha_{i_1},\dots, \alpha_{i_r}$ and
 $m_{ 1},\dots, m_{ r}$
used in (7.6) and (7.7).  (As in the classical case, using [Lu, 
Proposition 39.3.7], the 
description of $\tilde\theta(y_i)$ in (7.12) below is independent of 
the choice of such sequences.)

In the next theorem, we lift $\theta$ to a ${\bf C}$ algebra 
automorphism of $U$.   This is in the spirit of [L2, Theorem 3.1].
The main difference here is that we do not insist that $\tilde\theta$ 
is a ${\bf C}$ algebra involution on all of $U$.

\begin{theorem} There exists a ${\bf C}$ algebra automorphism    
$\tilde\theta$ on $U$ such that:
 \begin{enumerate}
\item[(7.9)] $ \tilde\theta(x_i)=x_i$ and $\tilde\theta(y_i) =y_i
$
for  all $ \alpha_i\in \pi_{\Theta}.$
\item[(7.10)] $\tilde\theta(\tau(\lambda))=\tau(\Theta(-\lambda))$
for  all $\tau(\lambda)\in T.$
\item[(7.11)]  $\tilde\theta(q)=q^{-1}.$
\item[(7.12)]
 $\tilde\theta(y_i)=
 [( {\rm ad}_r
x_{i_1}^{(m_1)}\cdots x_{i_r}^{(m_r)})t_{p(i)}^{-1}x_{p(i)}]$

\noindent
and $\tilde\theta(y_{p(i)})=
(-1)^{m(i)} [( {\rm ad}_r
x_{i_r}^{(m_r)}\cdots x_{i_1}^{(m_1)})t_{i}^{-1}x_{i}]$
 for $ \alpha_i\in \pi^{*}$.
 \end{enumerate}
Furthermore, $\tilde\theta^2$ is  the identity when restricted to  ${\cal 
M}$ and to $T$.   Finally, $\tilde\theta$ 
specializes to $\theta$ as $q$ goes to $1$.
  \end{theorem}

\medskip
\noindent
{\bf Proof:} To show that $\tilde\theta$ extends to a ${\bf C}$ algebra
automorphism of $U$, we relate it to Lusztig's automorphisms.
Let $T_{w_o}$ be Lusztig's automorphism
associated to $w_o$, the longest element of $W'$.  We follow the notation of [DK, Section 1.6]. Fix
$\alpha_i\in \pi_{\Theta}$. Recall that $-w_o(\alpha_i)=\alpha_{d(i)}$. By [DK, 
Section 1.6 and Proposition 1.6], $T_{w_o}$ sends $y_i$ to a 
nonzero scalar multiple of $x_{d(i)}t_{d(i)}^{-1}$, sends $x_i$ to a nonzero scalar 
multiple of $y_{d(i)}t_{d(i)}$,  and sends $t_i$ to $t_{d(i)}^{-1}$.
Furthermore, one checks using [DK, Remark 1.6]
that for each $\alpha_i\notin\pi_{\Theta}$  the composition $$(\iota\circ T_{w_o})(t_{p(i)}^{-1}x_{p(i)})=u_i[( {\rm ad}_r
x_{i_1}^{(m_1)}\cdots x_{i_r}^{(m_r)})t_{p(i)}^{-1}x_{p(i)}]$$ for 
some nonzero scalar $u_i$.   

Define a function $\tilde\theta$ on the generators of $U$ using 
(7.9), (7.10), (7.11), (7.12), and setting 
$$\tilde\theta(x_i)= u_i^{-1}(\iota\circ 
T_{w_o})(y_{p(i)}t_{p(i)}){\rm\ \ and\ \ }
\tilde\theta(x_{p(i)})= (-1)^{m(i)}u_{p(i)}^{-1}(\iota\circ 
T_{w_o})(y_{i}t_{i})$$ for each  $\alpha_i\in\pi^{*}$.
It is clear from (7.9) and (7.10)  that $\tilde\theta$ extends to a ${\bf 
C}$ algebra automorphism on both ${\cal M}$ and $T$.  Now
$\tilde\theta^2$ is clearly the identity on ${\cal M}$.  
Since $\Theta$ is an involution on the root system of ${\bf g}$,
condition (7.10) ensures that $\tilde\theta$ also restricts to   an
involution on the group $T$. 
 
We check that $\tilde\theta$ extends to a ${\bf C}$ algebra 
automorphism of $U$. In particular, 
$\tilde\theta(y_i)\tilde\theta(x_i)-\tilde\theta(x_i)\tilde\theta(y_i)
=(\iota\circ T_{w_o})(y_{p(i)}x_{p(i)}-x_{p(i)}y_{p(i)})=
\tilde\theta(y_ix_i-x_iy_i)$ for $\alpha_i\notin \pi_{\Theta}$. Furthermore, for $\alpha_i\in \pi_{\Theta}$,
$(\iota\circ T_{w_o})(t_{d(i)}x_{d(i)})=y_{i}=\tilde\theta(y_{i})$ up
to some nonzero scalar. Hence the
$\tilde\theta(y_i),1\leq i\leq n$ satisfy the quantum Serre
relations (1.7).  Similarly, $(\iota\circ
T_{w_o})(y_{i}t_i)=x_{d(i)}=\tilde\theta(x_{d(i)})$ up
to a nonzero scalar when $\alpha_i\in \pi_{\Theta}$.  It follows
that the $\tilde\theta(x_i)$ for $1\leq i\leq n$ 
 satisfy the quantum Serre relations (1.7).
Moreover, $\tilde\theta$ preserves the relations between the
$x_i$ and the $y_j$ for  $1\leq i,j\leq n$.  Thus $\tilde\theta$ extends to  a ${\bf C}$ algebra
automorphism $\tilde\theta$ of $U$.

Now consider an element $b$ in ${\cal M}\cup T\cup \{y_i|1\leq i\leq n\}$
and write $\bar b$ for its specialization as $q$ goes to $1$.   Note that the specialization of 
$\tilde\theta(b)$ is just $\theta(\bar b)$.   This is enough to force 
 $\tilde\theta$ to specialize to $\theta$. $\Box$

\medskip
We are now ready to introduce the quantum analog of $U({\bf g}^{\theta})$.
Set $$T_{\Theta}=\{\tau(\lambda)|\Theta(\lambda)=\lambda\},$$ a 
subgroup  of $T$.
Let $B=B_{\tilde\theta}$ be the subalgebra of $U$ generated by
${\cal M}$, $T_{\Theta}$, and the elements
$$B_i=y_it_i+\tilde\theta(y_i)t_i$$ for $\alpha_i\notin\pi_{\Theta}$. 
The next result shows that $B$ is a coideal subalgebra of $U$.   This 
fact combined with the results of Sections 1 and 4 is used below to describe 
the relations satisfied by these generators.   As a consequence, we 
show below that $B$ specializes to $U({\bf g}^{\theta})$ as $q$ goes 
to $1$.

\begin{theorem}
$B$ is a left coideal subalgebra of $U$.
\end{theorem}
{\bf Proof:} We need to check that $${\it\Delta}(b)\in U\otimes 
B\leqno{(7.13)}$$ for 
all $b\in B$. Since $\it\Delta$ is an algebra homomorphism from $U$ 
to $U\otimes U$, it is sufficient to check (7.13) for a set of 
generators of $B$.   Now  $B$ is generated 
by the elements $B_i$,  for $\alpha_i\notin\pi_{\Theta}$,
and  two Hopf algebras:  ${\cal M}$ and the group algebra generated 
by $T_{\Theta}$.
In particular, each $b\in {\cal M}$ and each $b\in 
 T_{\Theta}$ satisfies (7.13).  Hence it is sufficient to check (7.13) holds for the 
 remaining generators, namely when $b=B_i$ for 
 $\alpha_i\notin\pi_{\Theta}$. 

Set ${\cal M}^+=U^+\cap {\cal M}$. Note that $t_it_{p(i)}^{-1}$ is in 
$T_{\Theta}$ for all $i$ with 
$\alpha_i\notin\pi_{\Theta}$. Thus by (6.5) and the definition of
$\tilde\theta$,  the element $\tilde\theta(y_i)t_i$ satisfies the
following nice property with respect to the comultiplication of
$U$:
 $$\eqalign{{\it\Delta}(\tilde\theta(y_i)t_i)&\in
t_i\otimes \tilde\theta(y_i)t_{i}
 + U\otimes ({\cal M}\cap G^+U^o)t_{p(i)}^{-1}t_i
 \cr &\subset t_i\otimes \tilde\theta(y_i)t_{i}
 +U\otimes  {\cal M}^+T_{\Theta}.\cr}\leqno{(7.14)}$$
This combined with the formula for ${\it\Delta}(y_it_i)$ (see (1.8) 
and (1.10) )
yields  $${\it\Delta}(B_i)\in t_i\otimes
 B_i+ U\otimes {\cal M}^+T_{\Theta}\subset U\otimes B\leqno{(7.15)}$$ and the 
 theorem follows. $\Box$
 \medskip

We turn now to understanding the relations satisfied by the 
generators of $B$.  The elements $B_i$ have already been defined 
when $\alpha_i\notin\pi_{\Theta}$.
 Set $B_i=y_it_i$ for $\alpha_i\in \pi_{\Theta}$.  Given a tuple
$I=(i_1,\dots, i_r)$, set $|I|=r$, ${\rm wt}(I)=\alpha_{i_1}+\dots 
+\alpha_{i_r}$, $B_I=B_{i_1} \cdots B_{i_r}$, and 
$Y_I=y_{i_1}t_{i_1}\cdots y_{i_r}t_{i_r}$.

Recall ((1.4) and (1.11)) that  $x_iy_jt_j=q^{(-\alpha_i,\alpha_j)}y_jt_jx_i$  
whenever $j\neq i$.  
By Theorem 7.1, $\tilde\theta(x_i)=x_i$ whenever $\alpha_i\in 
\pi_{\Theta}$. Furthermore,   $\tilde\theta(y_j)$ and $y_j$ have the 
same weight with respect to the adjoint action of $T_{\Theta}$. Hence  
 $$x_iB_j=q^{(-\alpha_i,\alpha_j)}B_jx_i{\rm \ and \ }
\tau(\lambda)B_j=q^{(\lambda,-\alpha_j)}B_j\tau(\lambda)\leqno{(7.16)}$$  for
all $\alpha_i\in \pi_{\Theta}$ with $\alpha_j\notin \pi_{\Theta}$, and 
$\tau(\lambda)\in T_{\Theta}$.
It follows that  $$B=\sum_IB_I{\cal M}^+T_{\Theta}.\leqno{(7.17)}$$ 
Let ${\cal J}$ be a set such that $\{Y_{ J}|J\in {\cal J}\}$ is a
basis for $G^-$.  Note that $B_J=Y_J$ + (terms of higher weight) for each 
tuple $J$.   The triangular decomposition (4.1) of $U$ implies that 
the subspaces $\{Y_J{\cal M}^+T_{\Theta}|J\in {\cal J}\}$, and hence 
the subspaces $\{B_J{\cal M}^+T_{\Theta}|J\in {\cal J}\}$, are 
linearly independent.

Let $F_{ij}$ be the function in two
variables $X_1$ and $X_2$ defined by 
$$F_{ij}(X_1,X_2)=\sum_{m=0}^{1-a_{ij}}(-1)^m\left[\eqalign{1&-a_{ij}\cr 
 &m\cr}\right]_{q_i} X_1^{1-a_{ij}-m} X_2 X_1^{m}.$$ 
 The 
 quantum Serre relations (1.7) are the set of equations 
 $F_{ij}(y_i,y_j)=0$ for $i\neq j$.   A straightforward 
 computation shows that if $(\lambda_i,\alpha_j)=(\lambda_j,\alpha_i)$ 
 then
 $F_{ij}(y_i\tau(\lambda_i),y_j\tau(\lambda_j))=0$. 
Hence $$F_{ij}(y_it_i,y_jt_j)=0.\leqno{(7.18)}$$ 
It follows that the generators $y_it_i$  of $G^-$ satisfy the 
same relations as the generators of $U^-$.  Furthermore, since 
$(\Theta(-\alpha_i),\alpha_j)=(\Theta(-\alpha_j),\alpha_i)$, we have 
$$F_{ij}(\tilde\theta(y_i)t_i,\tilde\theta(y_j)t_j)=0.\leqno{(7.19)}$$

We show below that the $B_i$ for  $1\leq i\leq n$ satisfy relations which come 
from the quantum Serre relations on $G^-$.
First, we consider the evaluation of the function $F_{ij}$ at $B_i$, 
$B_j$ in a few special cases.   

If both $\alpha_i$ and $\alpha_j$ are in $\pi_{\Theta}$, then 
$F_{ij}(B_{i},B_j)=F_{ij}(y_it_i,y_jt_j)$. Similarly, if  $\alpha_i\in\pi_{\Theta}$
and $\alpha_j\notin\pi_{\Theta}$, then 
$F_{ij}(B_{i},B_j)=F_{ij}(y_it_i,y_jt_j)+F_{ij}(\tilde\theta(y_i)t_i,\tilde\theta(y_j)t_j).$   
Hence (7.18) and (7.19) imply that $${\rm if\ }\alpha_i\in 
\pi_{\Theta} {\rm  \ then \ } F_{ij}(B_i,B_j)=0.\leqno{(7.20)}$$

Now suppose that $i$ and $j$ are chosen such that $\pi_{0,0}(Y_{ij})$ 
is nonzero.   It follows that $Y_{ij}$ must have a zero weight 
summand.   Checking the possibilities for the quantum Serre 
relations, we must have  $a_{ij}=0$ and
$\Theta(\alpha_i)=-\alpha_j$. In particular,  
$B_i=y_it_i+q_i^{-2}x_jt_j^{-1}t_i$ and
$B_j=y_jt_j+q_i^{-2}x_it_i^{-1}t_j$.   A
straightforward computation shows that
$$\eqalign{{\rm if\ } &a_{ij}=0 {\rm\ and\ }\Theta(\alpha_i)=-\alpha_j {\rm 
\ then\ }\cr  
&F_{ij}(B_i,B_j)=B_iB_j-B_jB_i=(t_i^{-1}t_j-t_j^{-1}t_i)/(q_i-q_i^{-1}).\cr}\leqno{(7.21)}$$

Given $\lambda\in Q(\pi)$, let $P_{\lambda}$ be the projection of $B$ onto $U^-G^+\tau(\lambda)$
with respect to  the direct sum decomposition of Lemma 1.3 applied to the coideal $B$. 
The next lemma provides more detailed information about 
$F_{ij}(B_i,B_j)$.

\begin{lemma} Let 
$Y_{ij}=F_{ij}(B_i,B_j)$ for $i\neq j$ and 
$\lambda_{ij}=(1-a_{ij})\alpha_i+\alpha_j$.
 If $(\pi_{\beta,\gamma}\circ P_{\lambda_{ij}})(Y_{ij})\neq 0$ then 
 $[\beta,\gamma]\neq 0$, 
  $\tau(\lambda_{ij}-\beta)\notin T_{\Theta}$, and 
 $\tau(\lambda_{ij}-\gamma)\notin T_{\Theta}$.
 \end{lemma}

\noindent
{\bf Proof:}
Set $P_{ij}=P_{\lambda_{ij}}$.  Suppose that
$(\pi_{0,0}\circ P_{{ij}})(Y_{ij})\neq 0$.   It follows that 
$\pi_{0,0}(Y_{ij})\neq 0$.  
 Hence
$a_{ij}=0$ and 
$\Theta(\alpha_i)=-\alpha_j$.  Now $\lambda_{ij}=\alpha_i+\alpha_j$ 
in this case. By (7.21)
$P_{ij}(Y_{ij})=P_{ij}(t_i^{-1}t_j-t_j^{-1}t_i)=0$.  Therefore, $\pi_{0,0}(Y_{ij})=0$ for all 
choices of $i$ and $j$.

By (7.20), we  may assume that  $\alpha_i$  is not 
in $ \pi_{\Theta}$. 
Assume that $\beta$ and $\gamma$ are chosen so 
that
$\pi_{\beta,\gamma}(Y_{ij})\neq 0$. Note that $Y_{ij}$ can be written as a sum of 
monomials in $2-{a_{ij}}$ terms where
$1-a_{ij}$ of those terms are from the set 
$\{y_it_i,\tilde\theta(y_i)t_i\}$ and 
the other term is from the set $\{y_jt_j,\tilde\theta(y_j)t_j\}$.
It follows that
$\gamma=
s_1\alpha_{p(i)}+s_2\alpha_{p(j)}+\eta$ for some  $\eta\in 
Q^+(\pi_{\Theta})$ and nonnegative integers $s_1$ and $s_2$ such that
$ s_1\leq 1-a_{ij}$ and $ s_2\leq 1$.
Set $\gamma'=s_1\alpha_i+s_2\alpha_j$  and note that 
$\tau(\gamma-\gamma')$ is in $T_{\Theta}$. The above description of 
the monomials which add to  $Y_{ij}$ further implies  that 
$\gamma'+\beta\leq\lambda_{ij}$.  
Moreover,
   by (7.18), $0\leq \beta<\lambda_{ij}$ and by (7.19),
 $0\leq \gamma'< \lambda_{ij}$.  Now $\beta$ and $\gamma'$ are both linear 
 combinations of $\alpha_i$ and $\alpha_j$.  Thus the lemma follows if neither $\alpha_i$ nor 
 $\alpha_j$ are elements of $Q^+(\pi_{\Theta})$.   
 In the 
 case when $\alpha_j\in \pi_{\Theta}$, (7.19)  further implies that 
that  $0\leq \gamma'<\lambda_{ij}-\alpha_j$ and
 $0\leq \beta<\lambda_{ij}-\alpha_i$.  The lemma thus  follows in this case as well.
$\Box$

\medskip
The next result gives a description of the generators and relations 
of $B$.

\begin{theorem}  Let ${\tilde B}$ be the algebra freely generated over 
${\cal M}^+T_{\Theta}$ by the elements   $\tilde B_i$, $1\leq i\leq n$.
Then there exist elements $c^{ij}_J\in {\cal M}^+T_{\Theta}$
such that  $B\cong \tilde B/L$ where $L$ is the ideal generated by the 
following elements:
\begin{enumerate}
\item[(i)] $\tau(\lambda)\tilde
B_i\tau(-\lambda)-q^{-(\lambda,\alpha_i)}\tilde B_i$ for all
$\tau(\lambda)\in T_{\Theta}$ and $\alpha_i\notin \pi_{\Theta}$.
\item[(ii)]$ t_j^{-1} x_j\tilde B_i-\tilde B_i
t_j^{-1} x_j-\delta_{ij}(t_j-t_j^{-1})/(q_j-q_j^{-1})$ for all 
$\alpha_j\in \pi_{\Theta}$ and $1\leq i\leq n.$
\item[(iii)] 
$$\sum_{m=0}^{1-a_{ij}}(-1)^m\left[\eqalign{1&-a_{ij}\cr 
 &m\cr}\right]_{q_i}\tilde B_i^{1-a_{ij}-m}\tilde B_j\tilde B_i^{m}-
 \sum_{\{J\in {\cal J}|{\rm wt}(J)<(1-a_{ij})\alpha_i+\alpha_j\}}
 \tilde B_Jc^{ij}_J$$  for each $i\neq j$, $1\leq i,j\leq n$.
\end{enumerate}

\end{theorem}

\noindent
{\bf Proof:}  Relations (i) and (ii) follow from (7.16) and (1.4). 
We now show that the $B_i$, $1\leq i\leq n$, satisfy the relations 
described in  (iii).
 Fix a quantum Serre relation  
$Y=F_{ij}(B_i,B_j) $
 for   given $\alpha_i,\alpha_j$ with $i\neq j$.  Set 
 $\lambda=(1-a_{ij})\alpha_i+\alpha_j$ and $Z=P_{\lambda}(Y)$. By (4.7), 
$((P_{\lambda}\circ \pi_{0,0})\otimes Id){\it\Delta}(Y)=
(\pi_{0,0}\otimes Id){\it\Delta}(Z)$. Moreover, (4.7) ensures that 
$$(\pi_{0,0}\otimes Id){\it\Delta}(Z)=
 \tau(\lambda)\otimes Z.$$ 

By (7.15) and (7.17), we have  $${\it\Delta}(Y)\in\tau(\lambda)\otimes Y+ \sum_{\{J|\
{\rm wt}(J)<\lambda\}}U\otimes B_J{\cal M}^+T_{\Theta} .\leqno{(7.22)}$$
Now if $J$ has weight less than $\lambda$, one checks from (1.7) that 
there is no quantum Serre relation of weight greater than or equal to 
$-\lambda$.
Hence if  ${\rm wt}(J)<\lambda$ then $J$ is an element of the set ${\cal 
J}$. 
Now, (7.22)
 implies that 
 $$
 ((P_{\lambda}\circ\pi_{0,0})\otimes Id){\it\Delta}(Y)\in
  \tau(\lambda)\otimes 
( Y+ \sum_{\{J\in {\cal J}|\
wt(J)<\lambda\}} B_J{\cal M}^+T_{\Theta})$$ Thus we can
find $X\in \sum_{\{J\in {\cal J}|\wt(J)<\lambda\}}B_J{\cal M}^{+}T_{\Theta}$ such
that $Y+X=Z$.  We obtain a relation
of the form described in (iii) by proving $Z=0$.

Recall the notation of Section 4, Filtration II. Assume that $Z$ is 
nonzero and hence $\max(Z)$ is nonempty. Choose 
$[\beta,\gamma]\in{\max}(Z)$.   It follows that 
$\pi_{\beta,\gamma}(Z)\neq 0$.   By Lemma 7.3, $[\beta,\gamma]\neq 
[0,0]$, 
and neither $\tau(\lambda-\beta)$ nor $\tau(\lambda-\gamma)$ is an 
element of $T_{\Theta}$. Write $(\pi_{\beta,0}\otimes 
Id)({\it\Delta}( Z))=\sum v_i\otimes u_i$ where the $v_i\in 
U^-_{\beta}\tau(\lambda)$ and the $u_i\in G^+T$.  We may assume that 
the $v_i$ are linearly independent elements of 
$U^-_{-\beta}\tau(\lambda)$. Note that at
least one of the $u_i$ has a (nonzero)   summand of 
weight $\gamma$ in $G^+_{\gamma}\tau(\lambda-\beta)$.   By (7.19), the maximality of 
$[\beta,\gamma]$,  and the fact that $\beta\neq 0$, each $u_i$ is in $G^+T\cap 
\sum_{\{J\in {\cal J}| \wt(J)<\lambda\}}B_J{\cal 
M}^{+}T_{\Theta}$.  This intersection is just ${\cal 
M}^{+}T_{\Theta}$.  Hence $\tau(\lambda-\beta)\in T_{\Theta}$, a 
contradiction. This forces $\beta=0$.   It follows that 
$(\pi_{0,\gamma}\otimes Id)({\it\Delta}(Z))\in U\otimes \tau(\lambda-\gamma)$.   Again 
$\tau(\lambda-\gamma)$ must be in $T_{\Theta}$.   This contradiction 
forces $\max(Z)$ to be empty.   In particular, $Z=0$.

We have shown that $B$ is isomorphic to a homomorphic image of 
$\tilde B/L$.
A consequence of  relations (i), (ii), and (iii) is that ${\cal 
M}^+T_{\Theta}\tilde B_I\subset\sum_{J\in {\cal J}}\tilde B_J{\cal M}^+T_{\Theta}+L$ for each tuple $I$.
Thus  $$\tilde B/L=\bigoplus_{J\in {\cal J}}(\tilde B_J{\cal M}^+T_{\Theta}+L).$$ 
Since the elements $B_i$ in $B$ satisfy the 
relations (i), (ii), (iii), we also have the following direct sum 
decomposition:
$$B=\bigoplus_{J\in {\cal J}}(B_J{\cal M}^+T_{\Theta}).$$  Therefore $B\cong 
\tilde B/L$.  $\Box$

\medskip

Note that (7.20) and (7.21) both provide examples of the relations 
described in Theorem 7.4 (iii).   
 We illustrate how  to  compute the $c^{ij}_J$ in a more complicated  example.
  Consider the case where $\Theta(\alpha_i)=-\alpha_i$, 
 $\Theta(\alpha_j)=-\alpha_j$ and $a_{ij}=-1$.  So $B_i=y_it_i+q_i^{-2}x_i$
 and $B_j=y_jt_j+q_j^{-2}x_j$ and
  $Y=B_i^2B_j-(q_i+q_i^{-1})B_iB_jB_i+B_jB_i^2.$
   Thus by (1.9) and (1.10),
  $${\it\Delta}(B_r)=B_r\otimes 1+ 
 t_r\otimes B_r$$ for $r=i,j$. It follows that $${\it\Delta}(Y)=
 t_i^2t_j\otimes Y + 
 (B_i^2t_j-(q_i+q_i^{-1})B_it_jB_i+t_jB_i^2)\otimes B_j+ W\otimes B$$
 for some $W$ which satisfies $\pi_{0,0}(W)=0$.  A straightforward computation 
 using the relations of $U$ shows that
 $ P_{\lambda}\circ\pi_{0,0}( 
 (B_i^2t_j-(q_i+q_i^{-1})B_it_jB_i+t_jB_i^2)=-q_i^{-1}t_i^2t_j.$ 
 Thus $$0= (P_{\lambda}\circ\pi_{0,0})\otimes 
 Id({\it\Delta}(Y) =  t_i^2t_j\otimes Y -q_i^{-1}t_i^2t_j\otimes B_j.$$
 It follows that 
 $$B_i^2B_j-(q_i+q_i^{-1})B_iB_jB_i+B_jB_i^2=q_i^{-1}B_j.$$  (This 
 relation is also computed in  [L1, Lemma 2.2 (2.2)].  The generators 
 for $U$ and $B$ are somewhat different in [L1].   In particular, when 
 $\alpha_i=-\Theta(\alpha_i)$,  $B_i$ in [L1] is equal to 
 $y_it_i+x_i$ in the notation of this paper. Thus using a Hopf 
 algebra automorphism of $U$, the $B_i$ in [L1] corresponds to 
 $q_i^{-1}$ times the $B_i$ defined in this paper. This explains the 
 difference in coefficient of $B_j$ found in the two papers.)   Note 
 that a similar argument shows that $c^{ij}_J=0$ whenever 
 $-\alpha_i\neq \alpha_{p(i)}$ and $-\alpha_i\neq \alpha_{p(j)}$.  
 The $c^{ij}_J$ are computed in [L1, Lemma 2.2] for the cases when 
 $a_{ij}\geq -2$ and
  $\Theta(-\alpha_i)=\alpha_{p(i)}$.
\medskip

    Note that the 
  generators of $B$ specialize to the generators of  $U({\bf g}^{\theta})$ as $q$ goes to $1$.
 Thus the specialization of  $B_J{\cal M}^+T_{\Theta}$
 is contained in $U({\bf g}^{\theta})$. Moreover the set of spaces  $\{B_J{\cal 
 M}^+T_{\Theta}$, $J\in {\cal J}\}$ remain linearly independent after 
 specialization. As $q$ goes to $1$,  since these spaces 
  span $B$, we conclude that $B$ specializes
to $U({\bf g}^{\theta})$. 

The algebra $B$ also satisfies a  maximality
condition. Indeed, suppose that $C$ is a  subalgebra of $U$ containing $B$
and that $C$ also specializes to $U({\bf g}^{\theta})$.  Then by
[L2, Theorem 4.9], $C=B$.  The proof in [L2] uses a quantum version of 
the Iwasawa decomposition. The result also follows directly from 
Theorem 7.4.   The idea is as follows. Recall the notation of Section 6.
Set  $N^+_{\Theta}= N^+_{\pi_{\Theta}}$. 
By (6.6) (interchanging the roles of 
$N^+_{\pi'}$ with $N^-_{\pi'}$), we have  $$U= \sum_{J\in{\cal J}}
Y_J{\cal M}^+T
N_{ {\Theta}}^+.$$  By induction on $|J|$ (as in [L2, Lemma 4.3]), one can show that 
$U$ is 
spanned by the spaces $B$; $Bt$, $t\notin T_{\Theta}$; and 
$BT(N_{ {\Theta}}^+)_+$ where $(N_{ {\Theta}}^+)_+$ is the 
augmentation ideal of $N_{ {\Theta}}^+$.  Let $X$ be in $C$.   Subtracting  an 
element of $B$ if necessary, we may assume that $X$ is a linear combination 
of elements in $B(t-1)/(q-1)$ for  $t\notin T_{\Theta}$, and $BT(N^+_{{\Theta}})_+$. 
Assume that $X$ is nonzero. Rescale $X$ by a power of $(q-1)$
 so that it is an element of $\hat C-(q-1)\hat C$.   It 
follows that $X$ does not specialize to an element of $U({\bf 
g}^{\theta})$.  This contradiction forces $X=0$ and thus $B=C$.

We have shown that the algebra   $B$ satisfies the following properties.

\begin{enumerate}
\item[(7.23)]$B$ is a left coideal in $U$.
\item[(7.24)]$B$ specializes to $U({\bf g}^{\theta})$.
\item[(7.25)]If $B\subset C$ and $C$ is a subalgebra of $U$ which
specializes to $U({\bf g}^{\theta})$ then $B=C$.
\end{enumerate}
We now turn to characterizing all subalgebras of $U$ which satisfy
(7.23), (7.24), and (7.25). First, we present two variations which
satisfy these conditions as well.

\medskip
\noindent
 {\bf Variation 1}:  For  sake of simplicity, we assume first that
 $\bf g$ is simple. Recall the permutation $p$ used in (7.5). Suppose  that there exists an $r\in
 \{1,2,\dots,n\}$ such that
 $\alpha_r\notin\pi_{\Theta}$
and  $p(r)\neq r$. Assume further that
$(\alpha_r,\Theta(\alpha_r))\neq 0$. 
 Recall the  Cartan 
subspace ${\bf a}=\{x\in{\bf h}|\theta(x)=-x\}$ and 
the restricted root system $\Sigma$ associated 
to $\theta$ introduced at the 
beginning of this section. Let $\beta$ be the 
  restricted root 
corresponding to $e_{r}$. Note that $\beta$ is just the 
restriction of $\alpha_r\in {\bf h}^*$ to ${\bf a}^*$.  Furthermore,   
$(\ad a)[e_{r},\Theta(f_{r})]=2\beta(a)[e_{r},\Theta(f_{
r})]$  
for all $a\in {\bf a}$.   
In particular, 
 the restricted root system $\Sigma$ contains both $\beta$ and 
 $2\beta$.  Thus $\Sigma$
   is nonreduced
 and hence must be  of type $BC$. 
 One can 
 choose the positive roots of $\Sigma$ so that each $\alpha_j$ restricted 
 to ${\bf a}^*$ is either zero or 
 a simple positive root in $\Sigma$.  Furthermore, 
 $\alpha_r$ and $\alpha_j$ restrict to the same root if and only if 
 $j=r$ or $j=p(r)$.
 It follows from [Kn, Chapter II, Section 8] that there is exactly one positive simple root 
 in $\Sigma$ such that twice this root is also in $\Sigma$.    Hence 
 $r$ and $p(r)$ are the only values of $j$ such that
 $(\alpha_j,\Theta(\alpha_j))\neq 0$.   
 
 Let $c$ be an
element in $A={\bf C}[q,q^{-1}]_{(q-1)}$  which  specializes to $1$ as $q$ goes to $1$.
Define the ${\bf C}$ algebra automorphism $\tilde\theta_{c}$ of $U$ by
$$\tilde\theta_{c}(y_r)=c^{-1}\tilde\theta(y_r)$$
$$\tilde\theta_{c}(x_{r})=c\tilde\theta(x_{r})$$ and
$\tilde\theta_c$ agrees with $\tilde\theta$ on all other
generators of $U$.
Note that
$\tilde\theta_c$  is
 also a ${\bf C}$ algebra automorphism of $U$ which specializes
to $\theta$ and restricts to $\tilde\theta$ on ${\cal
M}T_{\Theta}$. Define $B_{\tilde\theta_c}$ in the same way as $B_{\tilde\theta}$
using $\tilde\theta_c$ instead of $\tilde\theta$. Thus 
$B_{\tilde\theta_c}$ is generated by ${\cal M}$, $T_{\Theta}$, and 
elements $B_i^c=y_it_i+\tilde\theta_c(y_i)t_i$ for 
$\alpha_i\notin\pi_{\Theta}$. Moreover, $B^c_{i}=B_i$ for $i\neq r$.
Since $\tilde\theta_c(y_i)$ is a scalar 
multiple of $\tilde\theta(y_i)$ for all $i$, the proof of Theorem 7.2 
also works for
  $B_{\tilde\theta_c}$.  Hence $B_{\tilde\theta_c}$ is a left coideal 
  subalgebra of $U$.
Consider a quantum Serre relation $F_{ij}(y_i,y_j)$
where either $i$ or $j$ equals $r$. Note that if
$\Theta(\alpha_i)=-\alpha_j$  then, by the assumptions on $r$,
$\{i,j\}=\{r,p(r)\}$ and $(\alpha_i,\alpha_j)\neq 0$.   Thus as in the 
proof of Lemma 7.3, $(\pi_{0,0}\circ P_{ij})(F_{ij}(B_i^c,B_j^c))\neq 0$ whenever $i\neq 
j$.
Hence the arguments for
$B_{\tilde\theta}$ used to prove Lemma 7.3 and Theorem 7.4 work for
$B_{\tilde\theta_c}$ as well. In particular,   $B_{\tilde\theta_c}$ satisfies
conditions (7.23), (7.24), and (7.25).

 Note that $B_{\tilde\theta_c}$ is not
isomorphic to $B_{\tilde\theta}$ via a Hopf algebra automorphism
of $U$ for $c\neq 1$.   It appears unlikely   in general that two such algebras
are isomorphic using just an algebra isomorphism.  
It should be noted that the existence
of this one parameter family of analogs is implicit in the proof
of [L2, Theorem 5.8]. However, it was mistakenly concluded in the 
paragraph directly preceding [L2, Theorem 
5.8]
 that all the analogs of Variation 1 were isomorphic to
$B_{\tilde\theta}$ via a Hopf algebra automorphism.

In the general semisimple case,
the one parameter $c$ is replaced by a multiparameter ${\bf c}$. In particular, 
each parameter corresponds to a pair of roots  $\alpha_{i_j}, 
\alpha_{p(i_j)}$
such that $(\alpha_{i_j},\Theta(\alpha_{i_j}))\neq 0$.   The 
automorphism
$\tilde\theta_{\bf c}$ is defined in a similar fashion to
$\tilde\theta_{ c}$. Let $[\Theta]$ be the set of automorphisms
of the form $\tilde\theta_{\bf c}$. Following the convention in
[L3], we refer to $B_{\tilde\theta'}$, $\tilde\theta'\in [\Theta]$
as a standard analog of $U({\bf g}^{\theta})$.

\medskip
\noindent {\bf Variation 2:}  Let ${\cal S}_1$ be the subset of
$\pi-\pi_{\Theta}$ consisting of those roots $\alpha_i$
 such that $\Theta(\alpha_i)=-\alpha_i$. Let
${\cal S}$ be the subset of ${\cal S}_1$ such that if $\alpha_i\in
{\cal S}$ and $\alpha_j\in {\cal S}_1$ then
$2(\alpha_i,\alpha_j)/(\alpha_j,\alpha_j)$ is even.  Let ${\bf S}$ be the set of $n$
tuples ${\bf s}=(s_1,\dots,s_n)$ such that each $s_i$ is in $A={\bf 
C}[q,q^{-1}]_{(q-1)}$ and
$s_i\neq 0$ implies $\alpha_i\in {\cal S}$.  Given $\hat\theta\in
[\Theta]$, let $B_{\hat\theta,{\bf s}}$ be the subalgebra of $U$
generated by $T_{\Theta}, {\cal M}$, the $B_i$ for $ \alpha_i\in
\pi-{\cal S}$, and the $B_{i,{\bf s} }$ defined by
$$B_{i, {\bf s}} =y_it_i + q^{-(\alpha_i,\alpha_i)}x_i+ s_it_i$$
for $\alpha_i\in {\cal S}$.
In particular, when the entries of ${\bf s}$ are all zero, $B_{i,{\bf 
s}}$ is just equal to $B_i$.

Note that $${\it\Delta}(B_{i,{\bf s}})= t_i\otimes B_i + (y_it_i +
q^{-(\alpha_i,\alpha_i)}x_i)\otimes 1.$$  Thus by the same
arguments as in Theorem 7.2 , $B_{\hat\theta, {\bf s}}$ is a left
coideal subalgebra of $U$.
Note that if $\tau(\lambda)\in T_{\Theta}$ and $\alpha_i\in {\cal S}$ 
then $(\lambda,\alpha_i)=0$.   It follows that $aB_{i,{\bf 
s}}=B_{i,{\bf s}}a$ for all $\alpha_i\in {\cal S}$ and $a\in {\cal 
M}T_{\Theta}$.  Recall the notation of  Lemma 7.3. To show that $B_{\hat\theta, {\bf s}}$
satisfies the conditions (7.23), (7.24), and (7.25), it suffices
to check for all $i,j$, $i\neq j$, that $(\pi_{0,0}\circ 
P_{ij})(F_{ij}(B_{i,{\bf s}}, 
B_{j,{\bf s}}))=0$.  A lengthy but routine computation shows that
this holds exactly when the $n$ tuple ${\bf s}$ is in ${\bf S}$.

Following the convention in [L3], the $B_{\hat\theta,{\bf s}}$,
$\hat\theta\in [\Theta]$ are called nonstandard analogs of $U({\bf
g}^{\theta})$. A nonstandard analog $B_{\hat\theta,{\bf s}}$ is not isomorphic to
a standard analog using a Hopf algebra automorphism of $U$.
However,  ([L2, Lemma 5.7])  $B_{\hat\theta,{\bf
s}}$ is isomorphic as an algebra to $B_{\hat\theta}$.  

Nonstandard analogs were first observed (to the
suprise of the author) in [L2, Section 5]. In [L3, Section 2],
nonstandard analogs were claimed to exist when ${\bf S}$ is defined 
using the larger set ${\cal S}_1$ instead of ${\cal S}$. (See in particular the definition of ${\bf S}$ given 
following [L3, (2.11)] and Theorem 2.1].) Our analysis in Variation 2 corrects this point.

\bigskip

We are now ready to show that the only possible subalgebras of $U$
which satisfy (7.23), (7.24), and (7.25)
 are our standard and nonstandard
analogs  associated to an automorphism in
$[\Theta]$. In particular, we give a new proof of [L2, Theorem
5.8] using the approach and results of Section 4.  Note that when the 
restricted roots $\Sigma$ associated to the involution $\theta$ do
not contain a component of type $BC$, then all the analogs described 
below are isomorphic to each other as algebras.   This is precisely 
what Theorem 5.8 in [L2]  states.   On the other hand, by the 
discussion  
of Variation 1, if $\Sigma$ 
contains $m$ components of type $BC$, then
there is an $m$ parameter family of analogs up to 
algebra isomorphism. 

\begin{theorem}  A subalgebra $B$ of $U$ satisfies (7.23), (7.24), and 
(7.25)
if and only if $B$ is isomorphic as an algebra to $B_{\hat\theta}$
for some $\hat\theta\in [\Theta]$. In particular, $B$ is
isomorphic to a standard or nonstandard analog of $U({\bf
g}^{\theta})$ corresponding to an element $\hat\theta$ in
$[\Theta]$ and an element ${\bf s}$ in ${\cal S}$  via a Hopf
algebra automorphism of $U$.
\end{theorem}

\noindent {\bf Proof:} 
We use the notation of the second filtration 
introduced in Section 4.  Let $B$ be a subalgebra of $U$ which
satisfies (7.23), (7.24), and (7.25). The proof of this theorem has 
three steps:
\begin{enumerate}
\item[(i)]   $B\cap T=T_{\Theta}$ 
\item[(ii)]    $B\cap U^+={\cal M}^+$ 
\item[(iii)]   $\gr_{\cal G}B\cap G^-=G^-.$ 
\end{enumerate}
More precisely, we first prove  that $B\cap 
T$ is a subgroup of $T_{\Theta}$ and $B\cap U^oU^+$ is a coideal 
subalgebra of ${\cal M}^+T_{\Theta}$. We then  use the second 
filtration introduced in Section 4 to analyze $\gr_{\cal G}B\cap G^-$
and thus prove (iii).  This information is then used  to show that 
$B\cap U^+U^o$ specializes to $U({\bf g}^{\theta})\cap U({\bf n}^++{\bf 
h})$. Next we  obtain (i) and (ii).  The last part of the proof takes a closer look at 
the generators of $B$ whose tip is in $G^-$ and show they are of the 
desired form.  The details follow.    

\medskip
Consider the set $B\cap T$. By (7.24), $B\cap T$ is a subset
of $T_{\Theta}$. Hence $B\cap T=B\cap T_{\Theta}$. 
Note that any element of $B$ can be written as a
direct sum of weight vectors with respect to $B\cap T$.  Hence  by
(7.25), we may assume that $ B\cap T_{\Theta}$ is a group.
Condition (7.23) and Lemma 4.2 ensure that $B\cap U^o$ is the group algebra
generated by $B\cap T_{\Theta}$.  Since $T_{\Theta}$ is free 
abelian of finite rank,  $B\cap T_{\Theta}$ is free abelian of
rank at most the   rank of $ T_{\Theta}$. 

Consider the coideal subalgebra $B\cap U^+U^o$ of $B$. 
We show that $B\cap U^+U^o$ is a subalgebra of ${\cal M}^+T_{\Theta}$. 
By Lemma 1.3, $B\cap U^oU^+$ is a direct sum of the 
vector spaces $B\cap G^+\tau(\mu)$, where $\tau(\mu)\in T$.
 Suppose that $c\in B\cap G^+\tau(\mu)$.     Choose $\gamma$ maximal 
 with respect to the standard partial ordering on $Q^+(\pi)$ so that
$\pi_{0,\gamma}(c)\neq 0$ and $\gamma\in Q^+(\pi_{\Theta})$.
 Then by (4.7), $$\pi_{0,\gamma}(c)\in
G^+_{\gamma}\tau(\mu)\otimes Y $$ where
$Y\in\tau(\mu-\gamma)+
\sum_{\gamma'>\gamma}G^+_{\gamma'-\gamma}\tau(\mu-\gamma)$.  Since $B$ 
is a coideal, $Y$ is an element of $B$.
Rescaling if necessary, we may assume that $Y$ is in $\hat B-(q-1)\hat 
B$.
Hence $Y$   specializes
to a nonzero element in $U({\bf g}^{\theta})$.  The choice of $\gamma$ implies that $\gamma'-\gamma\notin
Q^+(\pi_{\Theta})$ for all $\gamma'$ which appear in the definition of $Y$.
Hence, 
$Y\in \tau(\mu-\gamma)+(q-1)\sum_{\gamma'>\gamma}\hat 
G^+_{\gamma'-\gamma}\tau(\mu-\gamma)$.  But then $(q-1)^{-1}(Y-1)$ is
also in $\hat B$ and thus specializes to an
element of $U({\bf g}^{\theta})$. This forces $\tau(\mu-\gamma)$, and 
thus $\tau(\mu)$, to be in  
$T_{\Theta}$. Now consider  $\lambda$  maximal such that 
$\pi_{0,\lambda}(c)\neq 0$.   Then by (4.7), $\tau(\mu-\lambda)\in 
T_{\Theta}$.   Hence $\lambda\in Q^+(\pi_{\Theta})$. Note that if 
$\lambda'\in Q^+(\pi)$ and $\lambda'<\lambda$ then $\lambda'$ is also 
in $Q^+(\pi_{\Theta})$. It follows that 
$c$ is a sum of weight vectors with weights in $Q^+(\pi_{\Theta})$.   
In particular, $c\in {\cal M}^+T_{\Theta}$ and $B\cap U^oU^+$ is a 
subalgebra of ${\cal M}^+T_{\Theta}$.

  We next analyze the 
part of $B$ whose top degree terms are in $G^-$.
To do this, we introduce the left $B$ module $B/N$ where  $N$ is the 
left ideal $B
(B\cap (U^+U^o)_+)$ of $B$.  (Here $(U^+U^o)_+$ is equal to the augmentation ideal of $U^+U^o$.) 
  The filtration ${\cal G}$ on $B$ induces a 
 filtration which we also denote by ${\cal G}$ on $B/N$ which makes 
 $\gr_{\cal G}B/N$ into a $\gr_{\cal G}B$ module. By Theorem 4.9, the only 
 important contributions to this graded module occur in bidegree
 $(m,0)$ for 
  $m\geq 0 $.  In particular, $\gr_{\cal G}B/N$ is spanned by 
  elements $b+N$ where $b\in B$ and  ${\rm tip}(b)\in G^-$.
   Note that the subspace of $G^-$ of elements of bidegree 
  less than or equal to 
  $(m,0)$ is finite dimensional. 
    Thus 
the filtration on $B/N$ is a finite discrete filtration.
   Moreover, $\gr_{{\cal G}}B$ is finitely generated by the 
  image of the generators of $B$ described in  Corollary 4.10. Hence
 we have equality of Gelfand Kirillov dimension: $\gkdim \gr_{{\cal G}}(B/N)=\gkdim B/N$ ([KL, Prop. 6.6]).  
Now $\gr_{\cal G}B/N$ identifies with $\gr_{\cal G}(B)\cap G^-$ as a 
left 
$\gr_{\cal G}(B)\cap G^-$ module.   It is straightforward to 
check that the GK dimension of $ \gr_{\cal G}B/N$ as a $\gr_{\cal G}B$ 
module is equal to the GK dimension of  $ \gr_{\cal G}B/N$ as a 
$\gr_{\cal G}(B)\cap G^-$ module.
Hence, the form of the generators of $B$ 
 given in Corollary 4.10 implies that $\gkdim\gr_{\cal G}B/N\leq
 \dim {\bf n}^-$. 
 
 Let ${\bf r}$  be a Lie subalgebra of ${\bf g}^{\theta}$.  A standard argument 
 similar to the argument in the previous paragraph yields that 
 the $U({\bf g}^{\theta})$ module 
 $U({\bf g}^{\theta})/(U({\bf g}^{\theta}) {\bf r})$ has GK dimension 
 equal to $\dim {\bf g}^{\theta}-\dim {\bf r}$. (This follows for 
 example from [D, Proposition 2.2.7].) 
Consider the $\hat B$ 
 module $\hat B/\hat N$.  Write $\bar N$ for the specialization of 
 $N$  at $q=1$.  By Theorem 4.1, $
 {B\cap (U^+U^o)}$ specializes to the   enveloping algebra of 
 a Lie subalgebra, say ${\bf s}$,   of ${\bf g}^{\theta}$. Note that 
 $\bar N=U({\bf g}^{\theta}){\bf s}$. Now $B\cap 
 U^+U^o\subset {\cal M}^+T_{\Theta}$.  Hence  ${\bf s}$ is a Lie 
 subalgebra of ${\bf 
 m}^++({\bf g}^{\theta}\cap {\bf h})$. The map which sends each 
 $b+\hat N$ in $\hat B/\hat N$ to $\bar b+\bar N$ in $U({\bf g}^{\theta})/\bar N$
 allows us to   specialize the left $\hat B$ module 
   $\hat B/\hat N$ to the 
   $U({\bf g}^{\theta})$ module $U({\bf g}^{\theta})/\bar N$ at $q=1$.
 We can choose  generating sets for $B$ and $B/N$ which specialize to 
 generating sets of  $U({\bf g}^{\theta})$ and $U({\bf g}^{\theta})/\bar N$
 respectively. Hence
  $\gkdim B/N\geq \gkdim U({\bf g}^{\theta})/\bar N$. Note that  
  $$\eqalign{\gkdim U({\bf g}^{\theta})/\bar N&=\dim 
 {\bf g}^{\theta}-\dim {\bf s}\cr&\geq \dim {\bf g}^{\theta}-\dim ({\bf 
 m}^++({\bf g}^{\theta}\cap {\bf h}))\cr&=\dim {\bf n}^-.\cr}$$ By the previous 
 paragraph, this inequality is an equality.   Hence  
 $$\gkdim U({\bf g}^{\theta})/\bar N=\gkdim {\bf n}^-=\gkdim G^-.$$ Moreover 
 $\dim {\bf s}= \dim ({\bf 
 m}^++({\bf g}^{\theta}\cap {\bf h}))$. Since ${\bf s}$ is a 
 subalgebra of $ {\bf 
 m}^++({\bf g}^{\theta}\cap {\bf h})$, it follows that ${\bf s}= {\bf 
 m}^++({\bf g}^{\theta}\cap {\bf h})$.
 Thus $B\cap U^+U^o$ specializes to $U({\bf 
 m}^++({\bf g}^{\theta}\cap {\bf h}))$.

Recall the set $\Delta'$  of Corollary 4.10 used to define the 
generators of $B$ whose top degree term is in $G^-$. The 
description of the generators of $B$ in Corollary 4.10 implies 
 that $\gkdim B/N$ is equal to the number of elements in 
${\Delta}'$.  Since the number of elements in ${\Delta}^+$ is just 
the dimension of ${\bf n}^-$, it follows that 
${\Delta}'={\Delta}^+$.   Hence by Corollary 4.10, $B$ contains 
elements $y_it_i+b_i$, $1\leq i\leq n$, where $b_i$ is in $U^+U^o.$ 
It follows that $\tip(B)\cap G^-=G^-$.  This proves (iii).

Let $N'$ be the left ideal of $B\cap U^+U^o$ generated by the 
augmentation ideal of $B\cap U^o$.   We can analyze the left $B\cap 
U^+U^o$ module $(B\cap U^+U^o)/N'$ in a similar fashion to the analysis 
of   $B/N$.  It follows that  
 $B\cap U^+U^o=B\cap {\cal M}^+T_{\Theta}$ contains elements $x_i+c_i\in \hat B$ for each 
 $\alpha_i\in \pi_{\Theta}$.  Furthermore, $c_i\in U^o$ and $B\cap 
 U^o$ specializes to $U({\bf g}^{\theta}\cap {\bf h})$. Now $B\cap 
 U^o$ is just the group algebra generated by $B\cap T_{\Theta}$.
  Therefore, ${\rm rank\ } B\cap
T_{\Theta}={\rm rank\ } T_{\Theta}$.    Hence  we can find generators of
$T_{\Theta}$ such that a power of each generator lies in $B$. This
in turn implies that $B$ can be written as a direct sum of
$T_{\Theta}$ weight spaces. By the maximality condition (7.25) of $B$, we
obtain $B\cap T_{\Theta}=T_{\Theta}$.  This completes the proof of 
step (i).

Since
$T_{\Theta}\subset B$,  any element in $U^+U^o\cap B={\cal
M}^+T_{\Theta}\cap B$ is a sum of $T_{\Theta}$ weight vectors
contained in $B$.  Thus $x_i+c_i\in B$ implies $x_i\in B$.  In 
particular
$B$ contains $x_i$ for all $\alpha_i\in \pi_{\Theta}$.  Hence $B\cap 
U^+U^o={\cal M}^+T_{\Theta}$ and (ii) follows.

 Fix $i$ and consider again the element $y_it_i+b_i$ in $B$
 where $b_i\in U^+U^o$. Replacing $b_i$ by another 
 element in $U^+U^o$ if necessary, we may assume that
 $y_it_i+b_i$ is a weight vector for the
action of $T_{\Theta}$. By Lemma 1.3, we may further assume that $b_i\in 
G^+t_i$.  First consider the case when $\alpha_i\in \pi_{\Theta}$.   
Choose $\beta$ maximal with respect to the standard partial ordering 
on $Q(\pi)$ such that $\pi_{0,\beta}(b_i)\neq 0$.  By (4.7), $(\pi_{0,\beta}\otimes
Id){\it\Delta}(y_it_i+b_i)$ is a nonzero element of
$G^+_{\beta}t_i\otimes \tau(-\beta)t_i$. Hence $\tau(-\beta)t_i\in 
T_{\Theta}$ and  $\beta\in Q^+(\pi_{\Theta})$. If $0<\gamma<\beta$, 
then $\gamma$ is also in $Q^+(\pi_{\Theta})$.   Thus supp($b_i$) is a 
subset of $\{0\}\times Q^+(\pi_{\Theta})$. 
This   forces $b_i$ to be an element of ${\cal M}^+T_{\Theta}$ 
and so $y_it_i\in B$.  

Now assume that $\alpha_i\notin \pi_{\Theta}$. Choose $\beta$ such  that $[0,\beta]\in \max(b)$.
 Then by (4.11), $(\pi_{0,\beta}\otimes
Id){\it\Delta}(y_it_i+b_i)$ is a nonzero element of
$G^+_{\beta}t_i\otimes \tau(-\beta)t_i$.  In particular,
$\tau(-\beta)t_i=\tau(-\beta+\alpha_i)$ is in $T_{\Theta}$.
Since $\beta\in Q^+(\pi)$, it follows that $\beta\in 
\alpha_i+Q^+(\pi_{\Theta})$ or $\beta\in 
\alpha_{p(i)}+Q^+(\pi_{\Theta})$. However, $\beta$ must also be of the same $T_{\Theta}$ weight as
$-\alpha_i$.  The only possibility is  $\beta=\Theta(-\alpha_i)$.
By the
uniqueness property of the $Y_{I,j}$ and $X_{I,j}$ discussed in Section 6 (see (6.3) 
and the following discussion), the $\beta$ weight term is a scalar multiple
of $\tilde\theta(y_i)t_i$.  Indeed this is necessary in order for 
${\it\Delta}(y_it_i+b_i)$ to be an element of $U\otimes B$.
Therefore $b_i=c\tilde\theta(y_i)t_i+dt_i$ for some scalar $c$ and 
element $d\in G^+$ of bidegree less than 
$\bideg(\tilde\theta(y_i)t_i)$. By (7.23),
$${\it\Delta}(y_it_i+b_i)\in t_i\otimes (y_it_i+b_i) +U\otimes B.$$ By 
(1.8), (1.10), and (7.14), it follows that $${\it\Delta}(dt_i)\in 
t_i\otimes dt_i+U\otimes {\cal M}^+T_{\Theta}.$$  Since 
$\alpha_i\notin \pi_{\Theta}$, this forces $dt_i$ to be a scalar 
multiple of $t_i$.   Hence,  up to a Hopf algebra automorphism of $U$,
the only possibility for $B$ is one of the standard or nonstandard
analogs of $U({\bf g}^{\theta})$.  $\Box$
\medskip

Let us return for now to our first analog $B_{\tilde\theta}$.  Recall the definition of the
antiautomorphism $\kappa$.   One checks that $\kappa(({\rm
ad}_rx_j)b)=-(({\rm ad}_ry_j)\kappa(b))$ for any $b\in U$ and
$1\leq j\leq n$. Recall  that $m(i)=m_1+\cdots +m_r$.
Hence $$\kappa[( {\rm ad}_r x_{i_1}^{(m_1)}\cdots
x_{i_r}^{(m_r)})t_{p(i)}^{-1}x_{p(i)}]= (-1)^{m(i)}( {\rm ad}_r
y_{i_1}^{(m_1)}\cdots y_{i_r}^{(m_r)})y_{p(i)}.$$   
 A straightforward 
$U_q({\bf sl}\ 2)$ 
computation as in the classical case (see (7.8)) yields
$$( {\rm ad}_r
y_{i_1}^{(m_1)}\cdots y_{i_r}^{(m_r)})( {\rm ad}_r x_{i_r}^{(m_r)}\cdots
x_{i_1}^{(m_1)})t_{i}^{-1}x_{i}=t_{i}^{-1}x_{i}.$$
Set $y_j\cdot
bt_{p(i)}=bt_{p(i)}q^{(\alpha_{p(i)},\alpha_j)}y_j-y_jt_jbt_{p(i)}t_j^{-1}$
for any $b\in U$ and $1\leq i,j\leq n$. Note  that $y_j\cdot
bt_{p(i)}=(({\rm ad}_ry_j)b)t_{p(i)}$.  
Recall the definition of $\pi^*$ immediately 
following (7.5). We have
 $$\eqalign{(-1)^{m(i)}[y_{i_1}^{(m_1)}&\cdots
y_{i_r}^{(m_r)}\cdot B_{p(i)}]t_{p(i)}^{-1}t_i \cr =\
&( [(-1)^{m(i)}( {\rm ad}_r
y_{i_1}^{(m_1)}\cdots y_{i_r}^{(m_r)})y_{p(i)}]t_{i}+
t_i^{-1}x_{i}t_{i}\cr =\
&
\kappa(\tilde\theta(y_i)t_i)+
q^{-(\alpha_i,\alpha_i)}\kappa(y_it_i)\cr}$$ is an element of
$B_{\tilde\theta}$ for each $\alpha_i\in \pi^*$. A similar argument 
shows that 
$$\kappa(\tilde\theta(y_{p(i)})t_{p(i)})+
q^{-(\alpha_i,\alpha_i)}\kappa(y_{p(i)}t_{p(i)})$$ is also in 
$B_{\tilde\theta}$ for each $\alpha_i\in \pi^*$.
Thus one  
can find a Hopf algebra automorphism $\Upsilon$ in ${\cal H}_{\bf R}$ 
such  that $\Upsilon$ restricts to the identity on ${\cal M}$ and 
$T_{\Theta}$ and
$\Upsilon(B)$ contains $\kappa(\Upsilon(B_i))$ for each 
$\alpha_i\notin\pi_{\Theta}$. Furthermore, one can show that 
$\Upsilon(B)$ is generated by ${\cal M}$, $T_{\Theta}$,  and the 
$\kappa(\Upsilon(B_i))$, 
$\alpha_i\notin\pi_{\Theta}$.  Now $\kappa(\Upsilon({\cal M}))={\cal M}$ and 
$\kappa(\Upsilon(T_{\Theta}))=T_{\Theta}$.   It follows that  
$\kappa(\Upsilon(B))=\Upsilon(B)$.  Hence 
the results of Section 2 hold for $B$. 

The same argument works for analogs of Variations 1 and 2 provided
that all entries of the tuples involved are from ${\bf R}(q)$.  In
particular, let $[\Theta]_r$ be the set $\{\theta_{\bf b}| $ all
entries of ${\bf b}$ are in ${\bf R}(q)\}.$  We
refer to analogs of the form $B_{ \theta_{\bf b},{\bf s}}$ for
$\theta_{\bf b}\in [\Theta]_r$ and all entries of $\bf s$ are in
${\bf R}(q)$ as real analogs of $U({\bf g}^{\theta})$.
Given
$\theta_{\bf b}\in [\Theta]_r$, one can find $\Upsilon\in
{\cal H}_{\bf R}$ such that $\Upsilon^{-1}\kappa\Upsilon(B_{ \theta_{\bf b}})=B_{
\theta_{\bf b}}$. Furthermore, for any  $\bf s$ such that all of
its entries are in  ${\bf R}(q)$, we also have that $\Upsilon^{-1}\kappa\Upsilon(B_{ \theta_{\bf b},{\bf s}})=
B_{ \theta_{\bf b},{\bf s}}$. 
Hence, we may apply
results of Section 2 to all real analogs of $U({\bf g}^{\Theta})$.

Consider a real analog $B$ of $U({\bf g}^{\theta})$. Given a $U$ module $M$, set $X(M)$ equal to the
sum of all the finite-dimensional unitary $B$ submodules of $M$. The 
next result on
basic Harish-Chandra modules associated to the pair $U,B$ follows 
from Section 2.

\begin{theorem}  Let $B$ be a real analog of $U({\bf
g}^{\theta})$ and let $M$ be a $U$ module.  
Then any finite-dimensional $U$ module is a $B$ unitary module
 and a Harish-Chandra module  for the pair $U,B$.  Furthermore
both $F(U)$ and $X(M)$  
 are Harish-Chandra
modules for the pair $U,B$.
\end{theorem}

We continue the assumption that $B$ is a real analog of $U({\bf
g}^{\theta})$.   Using the approach of Section 3, we can define
the quantum homogeneous space associated to $B$.   The left
invariants  $R_q[G]^B_l$ are often referred to as  $R_q[G/K]$
(or ${\cal A}_q[G/K]$) in
the literature (see for example [NS,(2.5)]). Here $K$ can be thought of
merely as a symbol or as the complexification of the compact Lie
group in $G$ with Lie algebra ${\bf g}^{\theta}$.  Thus the 
homogeneous space $G/K$ is a 
 symmetric space.  
 The notation $R_q[G/K]$
suggests that the right $B$ invariants of $R_q[G]$ is the quantum
analog of the ring of regular functions on 
$G/K$. In [L3], it is shown that $B$ is a ``good" analog  of
$U({\bf g}^{\theta})$ for constructing quantum symmetric spaces in
the sense of [Di, end of Section 3]. In particular,  $R_q[G/K]$ has the same left
$U$ module structure as its classical counterpart (see Theorem 7.8
below). We summarize this and related  results here. A good survey on 
how to construct 
quantum symmetric spaces which includes a description of the 
classical situation is [Di]. For further 
information about classical symmetric spaces, the reader is referred 
to [He1] and [He2].  

A finite-dimensional
$U$ module $V$ is called a spherical module for $B$ if the space
of invariants  $V^B$ has dimension $1$. Recall the notion of Cartan 
subspace and restricted root system introduced at the beginning of 
this section.  Let ${\bf a}$ be  the
 Cartan subspace $\{x\in{\bf h}|\theta(x)=-x\}$ 
 and let $\Sigma$ be the associated restricted root system. 
Let $P^+_{\Theta}$ be the subset of $P^+(\pi)$ containing those
$\lambda$ such that 
\begin{enumerate}
\item[(i)]$ (\lambda,\beta)=0$ for all $\beta\in Q(\pi)$ such that
$\Theta(\beta)=\beta;$
\item[(ii)]the  restriction $\tilde\lambda$ of
$\lambda$ to ${\bf a}^*$   satisfies
$(\tilde\lambda,\beta)/(\beta,\beta)$ is an integer for
every   restricted root $\beta$. 
\end{enumerate}
 The set $P^+_{\Theta}$ is
exactly the set of dominant integral weights such that the
corresponding finite-dimensional simple ${\bf g}$ module is
spherical ([Kn, Theorem 8.49]).  By [L3, Theorem  4.2 and Theorem 4.3] we have the same
classification in the quantum case.

\begin{theorem}  Let $L(\lambda)$ be a finite-dimensional $U$ module
with highest weight ${\lambda}$ up to some possible roots of
unity. Then $$\dim L(\lambda)^B\leq 1.\leqno{(i)}$$ Moreover,
$$\dim L(\lambda)^B= 1 {\rm \ if \ and \ only \ if \ } \lambda\in
P^+_{\Theta}.\leqno{(ii)} $$

\end{theorem}

We sketch the proof here and refer the reader to [L3] for full
details.  Let $v_{\lambda}$ denote the highest weight generating
vector of $L(\lambda)$. Recall that for each $y\in G^-$ there
exists a $b\in B$ such that $b=y +$ higher weight terms. Now
$L(\lambda)$ is spanned by weight vectors of the form
$yv_{\lambda}$ where $y\in G^-$.  Hence $\dim
L(\lambda)/B^+v_{\lambda}\leq 1$ where $B^+$ is the augmentation ideal 
of $B$. Statement (i) follows from the
fact that $B^+v_{\lambda}\cap L(\lambda)^B$ is  empty. A
careful analysis using the form of the generators of $B$ further
shows that $v_{\lambda}\in B^+v_{\lambda}$ if and only if
$\lambda\notin P^+_{\Theta}$. This in turn implies (ii). The
argument turns out to be much more delicate when    $B$ is a
nonstandard analog. $\Box$

\medskip

 Theorem 7.7, the Peter-Weyl decomposition of $R_q[G]$, (3.1), 
and (3.4) 
imply  the following characterization of $R_q[G]^B_l$ as a
right $U$ module. 

\begin{theorem}  There is an isomorphism of right $U$ modules
$$R_q[G]^B_l\cong \bigoplus_{\lambda\in P^+_{\Theta}} L(\lambda)^*$$
\end{theorem}

There is an analogous statement for the right $B$ invariants of
$R_q[G]$. One can also describe the $B$ bi-invariants in a nice
way.   Identifying $L(\lambda)$ with a subspace of $R_q[G]^B_l$,
set ${\cal H} (\lambda)=R_q[G]^B_{bi}\cap L(\lambda)$. Note that
${\cal H} (\lambda)$ is a trivial left and right $B$ module.
Moreover,  by Theorem 7.7, ${\cal H} (\lambda)$ is
one-dimensional if $\lambda\in P^+_{\Theta}$ and zero otherwise.  The following direct sum decomposition into
trivial one-dimensional $B$ bimodules is thus an immediate
consequence of Theorem 7.8.

$$R_q[G]^B_{bi}\cong \bigoplus_{\lambda\in P^+_{\Theta}}{\cal
H}(\lambda).$$

Let ${\cal A}$ be the subgroup of $T$ consisting of those elements 
$\tau(\lambda)$ such that $\Theta(\lambda)=-\lambda$.
Thus ${\cal A}$ can be thought of as a quantum version of ${\bf a}$.
Let $W_{\Theta}$ denote the Weyl group associated to the restricted 
root system $\Sigma$.   Since $\Sigma\subset {\bf a}^*$, ${\bf a}$ and 
hence ${\cal A}$ inherit an action of $W_{\Theta}$.
 The author has recently shown that, $R_q[G]^B_{bi}$ is commutative
 and moreover is isomorphic to ${\bf C}(q)[{\cal A}]^{W_{\Theta}}$.
Thus, the ${\cal H}(\lambda)$ are 
natural choices of quantum zonal spherical functions (see [Di, the 
discussion concerning (3.4)]).   In special
cases, these quantum zonal spherical functions have been
determined to be Macdonald polynomials or other $q$ hypergeometric
series (see for example [K], [N], [DN], [NS]). Preliminary work by
the author suggests that this should be true in general.  

 It should be noted that these papers use analogs of
$U({\bf g}^{\theta})$ whose definition differs from the definition
of the $B_{\tilde\theta}$ and its variations found in this paper.
In [NM], one-sided coideal subalgebras are used. By [L2, Section
6], using  Theorem 7.5, these are shown  to be examples of the
analogs presented here. In other papers, two-sided coideals
analogs of ${\bf g}^{\theta}$ are used. The specialization of
these two-sided coideals generate   a much larger subalgebra than
$U({\bf g}^{\theta})$. The important object in these papers, used
to define quantum symmetric spaces, is the left ideal generated by
these two-sided coideals analogs of ${\bf g}^{\theta}$. It seems
likely that these left ideals can be shown to be generated by the
augmentation ideal of one of the analogs presented here. This is
certainly true for the left coideals studied in [K] and also for those in [N] ( combine [N,
Section 2.4] with [L2, Section 6]).

\bigskip
\centerline{REFERENCES}
\bigskip
\noindent [AJS]  H.H.  Andersen, J.C. Jantzen, and W. Soergel,
Representations of quantum groups at a $p$-th root of unity and of
semisimple groups in characteristic $p$: Independence of $p$, {\it
Ast\'erisque} {\bf 220}, Soc. Math.  France, Paris
(1994).

\medskip
\noindent [BF] W. Baldoni and P.M. Frajria, The quantum analog of a 
symmetric pair: a construction in type $(C_n, A_1\times C_{n-1})$,
{\it Trans. Amer. Math. Soc.} {\bf 8} (1997), 3235-3276.

\medskip
\noindent [CP] V. Chari and A. Pressley, {\it  A Guide to Quantum
Groups}, Cambridge University Press, Cambridge, (1995).

\medskip
\noindent [DK] C. DeConcini and V.G. Kac, Representations of
quantum groups at
 roots of 1, In:
{\it Operator Algebras, Unitary Representations, Enveloping
Algebras,
 and Invariant
Theory}, Progress in Math. {\bf 92}, Birkh\"auser, Boston (1990), 471-506.

\medskip
\noindent [Di] M.S. Dijkhuizen, Some remarks on the construction
of quantum symmetric spaces, In: {\it Representations of Lie
Groups, Lie Algebras and Their Quantum Analogues, } Acta Appl.
Math. {\bf 44} (1996), no. 1-2, 59-80.

\medskip
\noindent [DN] M.S. Dijkhuizen and M. Noumi, A family of quantum
projective spaces and related $q$-hypergeometric orthogonal
polynomials, {\it Trans. Amer. Math. Soc.} {\bf 350} (1998), no.
8, 3269-3296.

\medskip
\noindent [D] J. Dixmier, {\it Alg\`ebres Enveloppantes}, Cahiers
Scientifiques, XXXVII,

\noindent Gauthier-Villars, Paris (1974).

\medskip
\noindent [G] V. Guizzi, A classification of unitary highest weight 
modules of the quantum analogue of the symmetric pair $(A_n, 
A_{n-1})$, {\it J. Algebra} {\bf 192} (1997), 102-129.

\medskip
\noindent [H] J.E. Humphreys, {\it Introduction to Lie Algebras
and Representation Theory}, Springer-Verlag, New York
(1972).

\medskip
\noindent [HS] I. Heckenberger and K. Schmudgen,  Classification
of bicovariant differential calculi on the quantum groups ${\rm
SL}\sb q(n+1)$ and ${\rm Sp}\sb q(2n)$, {\it J. Reine Angew.
Math.} {\bf 502} (1998), 141-162.

\medskip
\noindent
[He1] S. Helgason, {\it Differential Geometry, Lie Groups, and Symmetric 
Spaces,}
 Pure and Applied Mathematics {\bf 80}, Academic Press,  New York  (1978).

\medskip
\noindent [He2]
S.  Helgason, {\it Groups and Geometric Analysis,
 Integral Geometry, Invariant Differential Operators, and Spherical 
 Functions,} Pure and Applied
Mathematics {\bf 113}, Academic Press, Inc., Orlando (1984).

\medskip
\noindent [J] N. Jacobson, {\it Basic Algebra II,}
W. H. Freeman and Company, San Francisco  (1980).

\medskip \noindent [JL1] A. Joseph and G. Letzter, Local
finiteness of the adjoint action for quantized enveloping
algebras, {\it J. of Algebra} {\bf 153} (1992), 289 -318.

\medskip
\noindent [JL2] A. Joseph and G. Letzter, Separation of variables
for quantized enveloping algebras, {\it American Journal of
Math.} {\bf 116} (1994), 127-177.

\medskip
\noindent[JL3] A. Joseph and G. Letzter, Verma module
annihilators for quantized enveloping algebras, {\it Ann. Sci.
Ecole Norm. Sup.}(4) {\bf 28} (1995), no. 4, 493-526.

\medskip
\noindent [Jo] A. Joseph, {\it Quantum Groups and Their Primitive
Ideals}, Springer-Verlag, New York (1995).

\medskip
\noindent
[Ka] V.G. Kac, {\it Infinite-Dimensional Lie Algebras,} Third ed., Cambridge University Press, 
Cambridge (1990).

\medskip
\noindent [Ke] M.S. K\'eb\'e, ${\cal O}$-alg\`ebres quantiques, {\it C.
R. Acad. Sci. Paris, Ser. I Math.} {\bf 322} (1996), no. 1, 1-4.

\medskip \noindent [KS] A.  Klimyk and K. Schm\"udgen,  {\it Quantum
Groups and Their Representations},  Texts and Monographs in
Physics, Springer-Verlag, Berlin (1997).

\medskip
\noindent
[Kn] A. W.  Knapp, {\it Lie Groups Beyond an Introduction,}
 Progress in Math. {\bf 140}, Birkh\"auser,  Boston
 (1996).

\bigskip
\noindent [K] T. Koornwinder,   Askey-Wilson polynomials as zonal
spherical functions on the ${\rm SU}(2)$ quantum group. {\it SIAM
J. Math. Anal.} {\bf 24} (1993), no. 3, 795--813.

\bigskip
\noindent [KL] G.R. Krause and T.H. Lenagan, {\it
 Growth of Algebras and Gelfand-Kirillov Dimension,} 
Research Notes in Mathematics {\bf 116}, Pitman, London (1985).

\medskip
\noindent [L1] G. Letzter, Subalgebras which appear in quantum
Iwasawa decompositions, {\it Canadian Journal of Math.}
  {\bf 49} (1997), no. 6, 1206-1223.

\medskip
\noindent [L2] G. Letzter,  Symmetric pairs for quantized
enveloping algebras, {\it J. Algebra} {\bf 220} (1999), no. 2,
729-767.

\medskip \noindent [L3] G. Letzter, Harish-Chandra modules for
quantum symmetric pairs, {\it Representation Theory, An Electronic
Journal of the AMS} {\bf 4} (1999)
64-96.

\medskip\noindent[Lu] G. Lusztig, {\it Introduction to Quantum Groups,}
 Progress in Math. {\bf 110}, Birkh\"auser,  Boston
 (1994).

\medskip\noindent
[M] S. Montgomery,  {\it Hopf Algebras and Their Actions on Rings,}
 CBMS Regional Conference Series in Mathematics {\bf 82},
  American Mathematical Society, Providence (1993). 

\medskip
\noindent [N] M. Noumi, Macdonald's symmetric polynomials as zonal
spherical functions on some quantum homogeneous spaces, {\it
Advances in Mathematics} {\bf 123} (1996), no. 1, 16-77.

\medskip
\noindent [NS] M. Noumi and T. Sugitani, Quantum symmetric spaces
and related q-orthogonal polynomials, In: {\it Group Theoretical
Methods in Physics (ICGTMP)} (Toyonaka, Japan, 1994), World
Sci.~Publishing, River Edge, N.J. (1995) 28-40.

\medskip
\noindent
[OV] A. L. Onishchik and E. B. Vinberg, {\it Lie Groups and Lie 
Algebras  III:
Structure of Lie Groups and Lie Algebras,}  Springer-Verlag, Berlin (1994).

\medskip
\noindent [R]  M. Rosso, Groupes Quantiques, Repr\'esentations
Lin\'eaires et Applications, Th\`ese, Paris 7  (1990).

\medskip
\noindent
[Ve] D.N. Verma, Structure of certain induced representations of 
complex semisimple Lie algebras, {\it Bulletin of the American 
Mathematical Society} {\bf 74}
(1968), 160-166.

\bigskip
\noindent
[V] D. Vogan, {\it Representations of Real Reductive Lie Groups,}
 Progress in Math.  {\bf 15}, Birkh\"auser, Boston (1981). 
 
\bigskip
\noindent

-----------------

\noindent letzter@math.vt.edu

\end{document}